\documentclass[11pt]{article}

\usepackage{a4wide}

\usepackage{amssymb,amsfonts,amsmath}
\usepackage{pb-diagram}

\newtheorem{theorem}{Theorem}
\newtheorem{lemma}{Lemma}[section]
\newtheorem{corollary}{Corollary}[section]
\newtheorem{proposition}{Proposition}[section]
\newtheorem{remark}{Remark}[section]

\title{Relatively free nilpotent torsion-free groups and their Lie algebras}
\author{C. Kofinas, V. Metaftsis and A.I. Papistas}
\date{}

\begin{document}
\parindent=0.8in
\parskip=0in
\renewcommand{\baselinestretch}{1.5}
\baselineskip=20pt plus 1pt
\newcommand{\qed}{$\Box$}
\newcommand{\pf}{\noindent {\it Proof.}\ }

\renewcommand{\thetheorem}{\Alph{theorem}}
\renewcommand{\theenumi}{(\Roman{enumi})}

\maketitle

\begin{abstract}

Let $K$ be a field of characteristic zero. For a torsion-free
finitely generated nilpotent group $G$, we naturally associate
four finite dimensional nilpotent Lie algebras over $K$, ${\cal
L}_{K}(G)$, ${\rm grad}^{(\ell)}({\cal L}_{K}(G))$, ${\rm
grad}^{(g)}(\exp {\cal L}_{K}(G))$ and $L_{K}(G)$. Let
${\mathfrak{T}}_{c}$ be a torsion-free variety of nilpotent groups
of class at most $c$. For a positive integer $n$, with $n \geq 2$,
let $F_{n}({\mathfrak{T}}_{c})$ be the relatively free group of
rank $n$ in ${\mathfrak{T}}_{c}$. We prove that ${\cal
L}_{K}(F_{n}({\mathfrak{T}}_{c}))$ is relatively free in some
variety of nilpotent Lie algebras, and ${\cal
L}_{K}(F_{n}({\mathfrak{T}}_{c})) \cong
L_{K}(F_{n}({\mathfrak{T}}_{c})) \cong {\rm grad}^{(\ell)}({\cal
L}_{K}(F_n({\mathfrak{T}}_{c}))) \cong {\rm grad}^{(g)}(\exp {\cal
L}_{K}(F_n({\mathfrak{T}}_{c})))$ as Lie algebras in a natural
way. Furthermore $F_{n}({\mathfrak{T}}_{c})$ is a Magnus nilpotent
group. Let $G_{1}$ and $G_{2}$ be torsion-free finitely generated
nilpotent groups which are quasi-isometric. We prove that if
$G_{1}$ and $G_{2}$ are relatively free of finite rank, then they
are isomorphic. Let $L$ be a relatively free nilpotent Lie algebra
over $\mathbb{Q}$ of finite rank freely generated by a set $X$.
Give on $L$ the structure of a group $R$, say, by means of the
Baker-Campbell-Hausdorff formula, and let $H$ be the subgroup of
$R$ generated by the set $X$. We show that $H$ is relatively free
in some variety of nilpotent groups; freely generated by the set
$X$, $H$ is Magnus and $L \cong {\cal L}_{\mathbb{Q}}(H) \cong
L_{{\mathbb{Q}}}(H)$ as Lie algebras. We extend the isomorphism
between ${\cal L}_K$ and $L_K$ to relatively free residually
torsion-free nilpotent groups. We also give an example of a
finitely generated Magnus nilpotent group $G$, not relatively
free, such that ${\cal L}_{\mathbb{Q}}(G)$ is not isomorphic to
$L_{\mathbb{Q}}(G)$ as Lie algebras.

\end{abstract}

\pagebreak
\small
\tableofcontents

\normalsize

\section{Introduction and Notation}

Let ${\mathbb Z}$ and $\mathbb{Q}$ denote the ring of integers and
the field of rational numbers, respectively. Furthermore we write
${\mathbb N}$ for the set of positive integers. For a group $G$
and $i \in {\mathbb N}$, we write $\gamma_{i}(G)$ for the $i$-th
term of the lower central series of $G$. Moreover we denote
$G^{\prime} = \gamma_{2}(G)$ i.e. the commutator subgroup of $G$.
For elements $a, b$ of $G$, we write $(a, b) = a^{-1}b^{-1}ab$,
and for $c \geq 3$, and elements $a_{1}, \ldots, a_{c}$ of $G$, we
define the left-normed group commutator $(a_{1}, \ldots, a_{c-1},
a_{c}) = ((a_{1}, \ldots, a_{c-1}), a_{c})$.  We call a group $G$,
a \emph{Magnus group} if each $\gamma_{i}(G)/\gamma_{i+1}(G)$ is
torsion-free and $\cap_{i \geq 1} \gamma_{i}(G) = \{1\}$. For $n
\in {\mathbb N}$, with $n \geq 2$, let $F_{n}$ be a free group of
rank $n$ freely generated by the set $\{f_{1}, \ldots, f_{n}\}$.
For a variety of groups $\mathfrak{G}$, let
${\mathfrak{G}}(F_{n})$ denote the verbal subgroup of $F_{n}$
corresponding to $\mathfrak{G}$. Also, let $F_{n}({\mathfrak{G}})
= F_{n}/{\mathfrak{G}}(F_{n})$: thus $F_{n}({\mathfrak{G}})$ is a
relatively free group of rank $n$ in $\mathfrak{G}$ and it has a
free generating set $\{x_{1}, \ldots, x_{n}\}$, where $x_{i} =
f_{i} {\mathfrak{G}}(F_{n})$, $i = 1, \ldots, n$. Note that the
verbal subgroups of $F_{n}$ are precisely the fully invariant
subgroups of $F_{n}$ (that is, the subgroups of $F_{n}$ which are
invariant under all group endomorphisms of $F_{n}$). The same
property holds for verbal subgroups and fully invariant subgroups
of relatively free groups. (For further information concerning
relatively free groups and varieties of groups see
\cite{neumann}.) Let ${\mathfrak{N}}_{c}$ be the variety of
nilpotent groups of class at most $c$, and let
${\mathfrak{T}}_{c}$ be a torsion-free subvariety of class at most
$c$ of ${\mathfrak{N}}_{c}$ (that is, its free groups of arbitrary
rank are torsion-free).

Let $K$ be a field of characteristic zero. We identify the prime
subfield of $K$ with ${\mathbb{Q}}$. By \lq \lq Lie algebra \rq
\rq (resp. \lq \lq Lie ring \rq \rq) we mean a Lie algebra over
$K$ (resp. over ${\mathbb Z}$). For $n \in {\mathbb N}$, with $n
\geq 2$, let $L_{n}$ denote a free Lie algebra of rank $n$ freely
generated by the set $\{\ell_{1}, \ldots, \ell_{n}\}$. Let
$\mathfrak{B}$ be a variety of Lie algebras and let
${\mathfrak{B}}(L_{n})$ be the fully invariant ideal of $L_{n}$
corresponding to ${\mathfrak{B}}$. Write $L_{n}(\mathfrak{V}) =
L_{n}/{\mathfrak{B}}(L_{n})$: thus $L_{n}({\mathfrak{B}})$ is a
relatively free Lie algebra of rank $n$ in ${\mathfrak{B}}$ and it
has a free generating set $\{\omega_{1}, \ldots, \omega_{n}\}$,
where $\omega_{i} = \ell_{i} + {\mathfrak{B}}(L_{n})$, $i = 1,
\ldots, n$. As for groups the fully invariant ideals of $L_{n}$
are precisely  the verbal ideals of $L_{n}$. The same property
holds for fully invariant ideals and verbal ideals of relatively
free Lie algebras. (For further information concerning relatively
free Lie algebras and varieties of Lie algebras see
\cite[Corollary 2.5, Chapter 14]{as}.) For $c \in {\mathbb N}$,
with $c \geq 2$, we write ${\mathfrak{M}}_{c}$ for the variety of
all Lie algebras defined by the identity $[\ell_{1}, \ldots,
\ell_{c+1}] = 0$: the variety of all Lie algebras which are
nilpotent of class at most $c$. Any variety of Lie algebras is
assumed to be non-trivial. We write ${\rm var}({\mathfrak{X}})$
for the variety generated by a set or a class ${\mathfrak{X}}$ of
Lie algebras. (We use similar definition as for groups \cite[page
18]{neumann}.)

Let $\widehat{L}_{n}$ be the completion of $L_{n}$ with respect to
the lower central series of $L_{n}$. (Recall that
$\widehat{L}_{n}$ is identified with the complete (unrestricted)
direct sum $\widehat{\oplus}_{m \geq 1}L^{m}_{n}$, and it has a
natural Lie algebra structure. Furthermore, $L_{n}$ is naturally
contained in $\widehat{L}_{n}$.) At this point, we state the
Baker-Campbell-Hausdorff formula (or briefly, BCH) (see \cite[page
178]{jacobson}, \cite[Chapter 5, Theorem 5.19]{mks}, \cite[Chapter
8]{bahturin}). It is
$$
\begin{array}{ccl}
X \circ Y
 & = & X + Y + \frac{1}{2}[X, Y] + \frac{1}{12}[X, Y, Y] -
 \frac{1}{12}[X, Y, X] + \cdots,
\end{array}
$$
where each term on the right-hand side is a rational multiple of a
Lie commutator $[Z_{1}, \ldots, Z_{m}]$, $m \in {\mathbb{N}}$, and
each $Z_{i}$ is $X$ or $Y$ and only finitely many terms of each
length occur. Note that the right-hand side of the aforementioned
formula is an infinite sum. The formula states that $X \circ Y$
belongs to the completion of the free Lie algebra freely generated
by the set $\{X, Y\}$ with respect to its lower central series.

 The BCH formula defines an
associative operation in $\widehat{L}_{n}$. For $x, y \in
\widehat{L}_{n}$ the operation $x \circ y$ is defined by
$$
x \circ y = x + y + \frac{1}{2} [x,y] + \frac{1}{12} [x, y, x] -
\frac{1}{12} [x, y, y] + \cdots
$$
(see, for example, \cite[Chapter 5, page 369]{mks}). (We remark
that if $L$ is a Lie algebra and the right-hand side of the BCH
formula has a meaning for all $x, y$ in $L$, then $L$ becomes a
group with respect to $\circ$.) It is easily verified that, for
$x, y \in \widehat{L}_{n}$,
$$
(x, y) = [x, y] + \frac{1}{2} [x, y, x] + \frac{1}{2} [x, y, y] +
\cdots.
$$
Let $L$ be a relatively free Lie algebra of rank $n$, with $n \geq
2$, freely generated by the set $\{t_{1}, \ldots, t_{n}\}$, and
let $\widehat{L}$ be the completion of $L$ with respect to the
lower central series. Write $\pi_{n}^{*}$ for the natural
epimorphism from $\widehat{L}_{n}$ onto $\widehat{L}$. Then
$\pi_{n}^{*}$ preserves multiplication $\circ$ and hence is a
group homomorphism from $(\widehat{L}_{n},\circ)$ into
$(\widehat{L},\circ)$. Thus, for $i, j \in \{1, \ldots, n\}$, we
have
$$
\pi_{n}^{*}(\ell_{i} \circ \ell_{j}) = t_{i} \circ t_{j} = t_{i} +
t_{j} + \frac{1}{2} [t_{i}, t_{j}] + \frac{1}{12} [t_{i}, t_{j},
t_{i}] - \frac{1}{12} [t_{i}, t_{j}, t_{j}] + \cdots.
$$
Notice that if $L$ is nilpotent, then $\widehat{L} = L$.

In section 2, for a torsion-free finitely generated nilpotent
group $G$, we naturally associate four finitely generated
nilpotent Lie algebras, namely, ${\cal L}_{K}(G)$, ${\rm
grad}^{(\ell)}({\cal L}_{K}(G))$, ${\rm grad}^{(g)}(\exp {\cal L}_{K}(G))$ and $L_{K}(G)$. One of our main
aims, in this paper, is to prove the following theorem.

\begin{theorem}
\begin{enumerate}
 \item Let ${\mathfrak{T}}_{c}$ be a
torsion-free variety of nilpotent groups of class at most $c$. For
a positive integer $n$, with $n \geq 2$, let
$F_{n}({\mathfrak{T}}_{c})$ be the relatively free group of rank
$n$ in ${\mathfrak{T}}_{c}$. Then the Lie algebra ${\cal
L}_{K}(F_{n}({\mathfrak{T}}_{c}))$ is relatively free in some
variety of nilpotent Lie algebras, and ${\cal
L}_{K}(F_{n}({\mathfrak{T}}_{c})) \cong
L_{K}(F_{n}({\mathfrak{T}}_{c}))$ as Lie algebras in a natural
way. Moreover, $F_{n}({\mathfrak{T}}_{c})$ is a Magnus
nilpotent group.
\item Let $L$ be a relatively free nilpotent Lie
algebra over $\mathbb{Q}$ of finite rank with a free generating
set $\cal X$. Give on $L$ the structure of a group $R$ by means of
the Baker-Campbell-Hausdorff formula. Let $H$ be the subgroup of
$R$ generated by the set $\cal X$. Then $H$ is relatively free in
some variety of nilpotent groups; freely generated by the set
$\cal X$, $H$ is Magnus and $L \cong {\cal L}_{{\mathbb{Q}}}(H)
\cong L_{{\mathbb{Q}}}(H)$ as Lie algebras in a natural way.
\end{enumerate}
\end{theorem}

In \cite{kp}, Kofinas and Papistas develop a method, by making use of Theorem A
(II), in order to study the automorphism group of a relatively
free nilpotent Lie algebra over ${\mathbb{Q}}$ of finite rank.

\begin{corollary}\label{1.1}
Let $G$ be a torsion-free finitely generated nilpotent group and
$K$ a field of characteristic zero. Then
\begin{enumerate}
\item ${\rm grad}^{(\ell)}({\cal L}_{K}(G)) \cong {\rm
grad}^{(g)}(\exp {\cal L}_{K}(G)) \cong L_{K}(G)$ as Lie algebras.

\item If $G$ is relatively free, then
$$
{\cal L}_{K}(G) \cong {\rm grad}^{(\ell)}({\cal L}_{K}(G))
$$
as Lie algebras.
\end{enumerate}
\end{corollary}

Now let $(X_1,d_1)$ and $(X_2,d_2)$ be metric spaces. A map $f: X_1\rightarrow X_2$
is called $(\lambda,\varepsilon)$-\emph{quasi-isometry} if there exist constants
$\lambda\ge 1$, $\varepsilon\ge 0$ and $C\ge 0$ such that
\begin{enumerate}
 \item $$\frac {1}{\lambda}d_1(x,y)-\varepsilon \le d_2(f(x),f(y))\le \lambda d_1(x,y)+\varepsilon$$
for all $x,y\in X_1$;
\item every point of $X_2$ lies in the $C$-neighbourhood of the image of $f$.
\end{enumerate}
Notice that the above map $f$ need not be continuous. Every
finitely generated group $G$ with generating set $S$ can be turned
into a metric space with the \emph{word metric} in $G$. If $F(S)$
is the free group with generating set $S$ and $\phi:
F(S)\rightarrow G$ is the natural projection, then the word metric
in $G$ is the metric obtained by defining $d_S(g_1,g_2)$ to be the
shortest word in the pre-image of $g_1^{-1}g_2$ under $\phi$. For
more on quasi-isometries the reader could consult \cite{bh}. The
metric space $(G,d_S)$ does not depend on the choice of $S$. In
fact if $S'$ is a different generating set for $G$ then $(G, d_S)$
and $(G, d_{S'})$ are quasi-isometric.

\begin{corollary}\label{1.2}
Let $G_{1}$ and $G_{2}$ be torsion-free finitely generated
nilpotent groups which are quasi-isometric. If $G_1$ and $G_{2}$
are relatively free of finite rank, then they are isomorphic.
\end{corollary}
Obviously, Corollary \ref{1.2} implies that the simply connected
nilpotent Lie groups $L_1$ and $L_2$ given by the Mal'cev completion
of $G_1$ and $G_2$ respectively, are isomorphic. It is a well know
conjecture whether the same result is true if we drop the
relative freeness assumption (see \cite{luck,fm}).

Throughout this paper, we write ${\mathfrak{L}}$ for a residually
torsion-free nilpotent variety of groups (that is, a variety with
its free groups of arbitrary rank to be residually torsion-free
nilpotent groups). Recall that a group $G$ is called
\emph{residually torsion-free nilpotent} if for any $g \in G
\setminus \{1\}$ there exists a normal subgroup $N_{g}$ such that
$g \notin N_{g}$ and $G/N_{g}$ is a torsion-free nilpotent group.
For a positive integer $n$, with $n \geq 2$, we write $G_{n} =
F_{n}({\mathfrak{L}})$. In section 3.3, for $G_{n}$, we naturally
associate two finitely generated Lie algebras, ${\cal
L}_{K}(G_{n})$ and $L_{K}(G_{n})$. The following result extends
Theorem A to residually torsion-free nilpotent groups.

\begin{theorem}
 \begin{enumerate}
  \item
 Let ${\mathfrak{L}}$ be a
residually torsion-free nilpotent variety of groups. For a
positive integer $n$, with $n \geq 2$, let $G_{n}$ be a relatively
free group of rank $n$ in ${\mathfrak{L}}$. Then the Lie
algebra ${\cal L}_{K}(G_{n})$ is relatively free, and ${\cal
L}_{K}(G_{n}) \cong L_{K}(G_{n})$ as Lie algebras in a natural
way.
\item Let $L$ be a relatively free Lie algebra over
${\mathbb{Q}}$ of rank $n$, with $n \geq 2$, freely generated by
the set $\cal X$. Let $\widehat{L}$ be the completion of $L$ with
respect to the lower central series. Give on $\widehat{L}$ the
structure of a group $\widehat{R}$ via the
Baker-Campbell-Hausdorff formula, and consider $L$ as a Lie subalgebra of
$\widehat{L}$. Let $H$ be the subgroup of $\widehat{R}$ generated
by the set $\cal X$. Then $H$ is a relatively free residually
torsion-free nilpotent group of rank $n$ freely generated by the
set $\cal X$. A Lie algebra $\Lambda_{{\mathbb{Q}}}(H)$ over
${\mathbb{Q}}$, associated with $H$, is constructed such that $L
\cong \Lambda_{{\mathbb{Q}}}(H)$ as Lie algebras in a natural way.
Furthermore, $\Lambda_{{\mathbb{Q}}}(H)$ is a homomorphic image of
${\cal L}_{{\mathbb{Q}}}(H)$.
\end{enumerate}
\end{theorem}

Notice that there is some overlap of Theorem B with some results
in \cite{bahturin}. Namely, the Lie algebra $L_{K}(G_{n})$ is
proved to be relatively free \cite[Theorem 10, page
278]{bahturin}, and $H$ to be a relatively free residually
torsion-free nilpotent group of rank $n$ \cite[Theorem 8, pages
276-277. See, also, Comments on pages 296-297]{bahturin}.

The paper is organized as follows: In section 2, for any
torsion-free finitely generated nilpotent group $G$, four finitely
generated nilpotent Lie algebras are naturally defined, and for a
finitely generated nilpotent Lie algebra, we naturally associate a
torsion-free finitely generated nilpotent group by means of the
Baker-Campbell-Hausdorff formula. Some auxiliary lemmas are proved
in section 3. Moreover, relatively free groups and relatively Lie
algebras are studied. In section 4, we prove Theorems A and B, and
Corollaries 1.1 and 1.2. An example of a finitely generated Magnus
nilpotent group $G$, not relatively free, such that ${\cal
L}_{\mathbb{Q}}(G) \ncong L_{\mathbb{Q}}(G)$ is given in section
5.

\section{Nilpotent groups and Lie algebras}

\renewcommand{\thetheorem}{\arabic{theorem}}

Let $H$ be a nilpotent group and denote by $\tau(H)$ the set of
all elements of finite order in $H$. Then $\tau(H)$ is a subgroup
of $H$, it is characteristic in $H$, and $H/\tau(H)$ is
torsion-free. For a group $N$ and a positive integer $i$, let
$\pi_{i}$ be the natural mapping from $N$ onto $N/\gamma_{i}(N)$.
Since $N/\gamma_{i}(N)$ is nilpotent, $\tau(N/\gamma_{i}(N))$ is a
group. Write $\tau_{i}(N)$ for the complete inverse image of
$\tau(N/\gamma_{i}(N))$ in $N$ via $\pi_{i}$ i.e. $\tau_{i}(N) =
\{g \in N: g^{n} \in \gamma_{i}(N)~{\rm for~some~integer}~n\}$. We
call $\tau_{i}(N)$ the \emph{isolator} of $\gamma_{i}(N)$ in $N$.
Note that $\tau_{i}(N)/\gamma_{i}(N) = \tau(N/\gamma_{i}(N))$ for
all $i$.

Let ${\mathfrak{L}}$ be a residually torsion-free nilpotent
variety of groups. For a positive integer $n$, with $n \geq 2$, we
write $G_{n} = F_{n}({\mathfrak{L}})$. The condition of being
residually torsion-free nilpotent is equivalent to $\cap_{i \geq
1} \tau_{i}(G_{n}) = \{1\}$. Since each $\tau_{i}(G_{n})$ is a
fully invariant subgroup of $G_{n}$, it is easily verified that
there are no repetitions of terms of the series
$\{\tau_{i}(G_{n})\}_{i \geq 1}$. Notice that $(\tau_{i}(G_{n}),
\tau_{j}(G_{n})) \leq \tau_{i+j}(G)$ for all $i, j$. For each
positive integer $i$, we write $L_{i}(G_{n})$ for the quotient
group $\tau_{i}(G_{n})/\tau_{i+1}(G_{n})$. Form the (restricted)
direct sum of abelian groups $L(G_{n}) = \oplus_{i \geq
1}L_{i}(G_{n})$ and give on it a structure of a Lie ring by
defining a Lie multiplication $[a \tau_{i+1}(G_{n}), b
\tau_{j+1}(G_{n})] = (a, b) \tau_{i+j+1}(G_{n})$, where $a
\tau_{i+1}(G_{n})$ and $b \tau_{j+1}(G_{n})$ are the images of the
elements $a \in \tau_{i}(G_{n})$ and $b \in \tau_{j}(G_{n})$ in
the quotient groups $L_{i}(G_{n})$ and $L_{j}(G_{n})$,
respectively, and $(a, b) \tau_{i+j+1}(G_{n})$ is the image of the
group commutator $(a, b)$ in the quotient group $L_{i+j}(G_{n})$.
Multiplication is then extended to $L(G_{n})$ by linearity. Form
the tensor product of $K$ with $L(G_{n})$ over $\mathbb{Z}$ and
write $L_{K}(G_{n}) = K \otimes_{\mathbb{Z}} L(G_{n})$. Then
$L_{K}(G_{n})$ has the structure of a Lie algebra with
$\lambda(\lambda^{\prime} \otimes a) = \lambda \lambda^{\prime}
\otimes a$ and $[\lambda \otimes a, \lambda^{\prime} \otimes
a^{\prime}] = \lambda \lambda^{\prime} \otimes [a, a^{\prime}]$
for all $\lambda, \lambda^{\prime} \in K$ and $a, a^{\prime} \in
L(G_{n})$. Since each $L_{i}(G_{n})$ is a free ${\mathbb
Z}$-module with a basis, say $X_{i}$, every element of $K
\otimes_{\mathbb{Z}} L_{i}(G_{n})$ may be written uniquely as a
$K$-linear combination of elements $1 \otimes x$ with $x \in
X_{i}$. We write $L_{i,K}(G_{n})$ for the vector space over $K$
spanned by any $\mathbb{Z}$-basis of $L_{i}(G_{n})$. Thus we may
regard $L_{i}(G_{n})$ as a subset of $L_{i,K}(G_{n})$ and so, we
regard $L(G_{n})$ as a subset of $L_{K}(G_{n})$. Furthermore
$L_{K}(G_{n}) = \oplus_{i \geq 1} L_{i,K}(G_{n})$.

For the rest of this section, $G$ denotes a torsion-free finitely
generated nilpotent group of class $c$. The series
\begin{eqnarray}
G = \tau_{1}(G) \supset \tau_{2}(G) \supset \cdots \supset
\tau_{c}(G) \supset \tau_{c+1}(G) = \{1\}
\end{eqnarray}
is a characteristic central series of $G$ with
$\tau_{i}(G)/\tau_{i+1}(G)$ torsion-free for all $i$, $1 \leq i
\leq c$ (see \cite[page 49]{segal}). Form the direct sum of
abelian groups $L(G) = \oplus_{t=1}^{c} \tau_{t}(G)/\tau_{t+1}(G)$
and let $L_{K}(G) = \oplus_{t=1}^{c} (K \otimes
\tau_{t}(G)/\tau_{t+1}(G))$. As before, we give on $L_{K}(G)$ the
structure of a Lie algebra. Let $n_{i}$ denote the rank of the
free abelian group $L_{i}(G)$. For $i = 1, \ldots, c$, let $f(i) =
n_{1} + \cdots + n_{i-1}$, with $n_{0} = 0$ and $n_{1} = n$. Let
$X_{i} = \{a_{f(i)+1}, \ldots, a_{f(i+1)}\}$ be a subset of
$\tau_{i}(G)$ such that the set $\{a_{f(i)+1}\tau_{i+1}(G),
\ldots, a_{f(i+1)}\tau_{i+1}(G)\}$ is a ${\mathbb Z}$-basis of
$L_{i}(G)$. We refine the series $(1)$ of $G$ to obtain a central
series
\begin{eqnarray}
G = {\cal G}_{1} \supset \cdots \supset {\cal G}_{f(i)+1} \supset
\cdots \supset {\cal G}_{f(i+1)} \supset \cdots \supset {\cal
G}_{f(c+1)} \supset {\cal G}_{f(c+1)+1} = \{1\},
\end{eqnarray}
with $i = 1, \ldots, c$, such that $a_{f(i)+j}$ is a
representative in $G$ of a generating element of ${\cal
G}_{f(i)+j}$ modulo ${\cal G}_{f(i)+j+1}$, with $j = 1, \ldots,
n_{i}$. (The length $f(c+1)$ is an invariant for $G$. It is called
the \emph{Hirsch number} of $G$; denoted by ${\cal H}(G)$.)
Following \cite{jennings}, we call the aforementioned central series of
$G$, a ${\cal G}$-\emph{series} of $G$. (In \cite{jennings} it is called
${\cal F}$-series.) Thus every element $g$ of $G$ may be written
uniquely in the form
\begin{eqnarray}
g = a^{\beta_{1}}_{1} \cdots a^{\beta_{n_{1}}}_{n_{1}} \cdots
a^{\beta_{f(i)+1}}_{f(i)+1} \cdots a^{\beta_{f(i+1)}}_{f(i+1)}
\cdots a^{\beta_{f(c)+1}}_{f(c)+1} \cdots
a^{\beta_{f(c+1)}}_{f(c+1)},
\end{eqnarray}
where $\beta_{j} \in {\mathbb{Z}}$. In what follows we assume that
the aforementioned series (2) and the elements $a_{1}, \ldots,
a_{f(c+1)}$ in $(3)$ have been selected. The set $\{a_{1}, \ldots,
a_{f(c+1)}\}$ is called a \emph{canonical basis} (or,
\emph{Mal'cev basis}) of $G$.

Let $\Gamma = KG$ be the group algebra of $G$ over $K$, and let
$\Delta$ be the augmentation ideal of $\Gamma$. It has been proved
in \cite[Theorem 4.3]{jennings} that $\cap_{i \geq 1} \Delta^{i} =
\{0\}$. Take the set $\{\Delta^{i}\}_{i \geq 1}$ as a fundamental
system of neighbourhoods of the element $0$ in $\Gamma$; then a
sequence $b_{1}, b_{2}, \ldots, b_{n}, \ldots$ of elements of
$\Gamma$ converges to $b \in \Gamma$ if, for every $i$, there
exists an integer $n(i)$ such that $n > n(i)$ implies that $b_{n}
- b \in \Delta^{i}$. So $\{\Delta^i\}$ induce a topology on
$\Gamma$ and let $\Gamma^{*}$ be the completion of $\Gamma$ in
this topology, and $\Delta^{*}$ be the completion of $\Delta$. We
may consider $\Gamma^{*}$ to be the algebra of all \emph{formal
power series} $a^{*}$ of the form $a^{*} = \alpha_{0} + \sum
\alpha_{k} d_{k}$, where $\alpha_{0}, \alpha_{k} \in K$, $k = 1,
2, \ldots$, and $d_{k} \in \Delta^{k}$, while $\Delta^{*}$
consists of all elements $a^{*}$ with $\alpha_{0} = 0$. We
identify $\Gamma$ with its isomorphic image in $\Gamma^{*}$.
Define $\exp : \Delta^{*} \longrightarrow 1 + \Delta^{*}$ and
$\log : 1 + \Delta^{*} \longrightarrow \Delta^{*}$ as in
\cite{jennings}.

We associate with $\Delta^{*}$ the Lie algebra $\Lambda^{*} =
(\Delta^{*})_{L}$ in the usual way by defining the binary
operation of commutation in $\Lambda^{*}$ by means of $[x^{*},
y^{*}] = x^{*}y^{*} - y^{*}x^{*}$ for all $x^{*}, y^{*} \in
\Delta^{*}$. For a positive integer $c$, with $c \geq 3$, we
define the left normed Lie product $[y^{*}_{1}, \ldots,
y^{*}_{c-1}, y^{*}_{c}] = [[y^{*}_{1}, \ldots, y^{*}_{c-1}],
y^{*}_{c}]$. It is proved in \cite{jennings} that $\cap_{i \geq
1}(\Delta^{*})^{i} = \{0\}$. Furthermore it is easily verified
that if $\Theta^{*}$ is an ideal of $\Delta^{*}$, and $M^{*} =
(\Theta^{*})_{L}$, then $M^{*}$ is an ideal of the Lie algebra
$\Lambda^{*}$. Thus we have $\Lambda^{*} = (\Delta^{*})_{L}
\supset (\Delta^{*})^2_{L} \supset \cdots$. For a positive integer
$i$, let $\gamma_{i}(\Lambda^{*})$ denote the $i$-th term of the
lower central series of $\Lambda^{*}$. Then, for all $i$,
$\gamma_{i}(\Lambda^{*}) \subseteq (\Delta^{*})^i_{L}$. Therefore
$\cap_{i \geq 1} \gamma_{i}(\Lambda^{*}) = \{0\}$ and we obtain
$\Lambda^{*} = \gamma_{1}(\Lambda^{*}) \supset
\gamma_{2}(\Lambda^{*}) \supset \cdots \supset
\gamma_{i}(\Lambda^{*}) \supset \cdots$. If $v_{k} \in
\gamma_{k}(\Lambda^{*})$ for $k = 1, 2, \ldots$, then the infinite
series $v_{1} + v_{2} + \cdots + v_{k} + \cdots$ is an element of
$\Lambda^{*}$ (see \cite[Lemma 6.1]{jennings}).

The BCH formula reveals the intimate connection between the group
$1 + \Delta^{*}$ and the Lie algebra $\Lambda^{*}$. For $X = 1 +
x^{*}$ and $Y = 1 + y^{*} \in 1 + \Delta^{*}$,
$$
\log XY = \log X + \log Y + \frac{1}{2} [\log X, \log Y] +
\frac{1}{12} [\log X, \log Y, \log Y] - \cdots
$$
where each term on the right-hand side is a rational multiple of a
Lie commutator $[Z_{1}, \ldots, Z_{m}]$, $m \in {\mathbb N}$,
and each $Z_{i}$ is $\log X$ or $\log Y$ and only finitely many
terms of each length occur. Note that the right-hand side of the
aforementioned formula is an infinite sum and therefore
convergent. Since $\log$ and $\exp$ are mutually inverse
bijections, we define an operation on $\Lambda^{*} =
(\Delta^{*})_{L}$ as follows: Let $u, v \in \Lambda^{*}$. Then
there exists unique $X$ and $Y$ in $1 + \Delta^{*}$ such that
$\log X = u ~{\rm and}~\log Y = v$. Define
\[
\begin{array}{rcl}
u \circ_{G} v & = & \log X \circ_{G}  \log Y \\
 & = & \log XY  \\
  & = & u + v + \frac{1}{2} [u, v] + \frac{1}{12} [u, v, v] -
 \cdots.
\end{array}
\]
Notice that $(\Lambda^{*}, \circ_{G})$ is a group. We write
$\circ$ instead of $\circ_{G}$ if it is clear in the context. Let
${\cal L}_{K}(G)$ be the vector space over $K$ spanned by all
$\log g$ with $g \in G$. Then ${\cal L}_{K}(G)$ is a nilpotent Lie
subalgebra of class $c$ of $\Lambda^{*}$, and the set $\{\log
a_{1}, \ldots, \log a_{f(c+1)}\}$ is a $K$-basis of ${\cal
L}_{K}(G)$ (see \cite[Theorem 7.3]{jennings}). Thus $\dim {\cal
L}_{K}(G) = {\cal H}(G)$. Notice that ${\cal L}_{K}(G) = K
\otimes_{\mathbb{Q}} {\cal L}_{\mathbb{Q}}(G)$. Observe that $\exp
{\cal L}_{K}(G)$ is a subgroup of $1 + \Delta^{*}$ and $(\exp
u)(\exp v) = \exp (u \circ v)$ for all $u, v \in {\cal L}_{K}(G)$.
Furthermore $({\cal L}_{K}(G), \circ)$ is a subgroup of
$(\Lambda^{*}, \circ)$. It is easily verified that ${\cal
L}_{K}(G)$ is isomorphic as group to $\exp {\cal L}_{K}(G)$ by a
mapping sending $u$ to $\exp u$ for all $u \in {\cal L}_{K}(G)$.
Form the direct sum of the vector spaces over $K$
$$
{\rm grad}^{(\ell)}({\cal L}_{K}(G)) = \oplus_{t=1}^{c}
\gamma_{t}({\cal L}_{K}(G))/\gamma_{t+1}({\cal L}_{K}(G)).
$$
The Lie multiplication in ${\cal L}_{K}(G)$ induces a Lie algebra
structure on ${\rm grad}^{(\ell)}({\cal L}_{K}(G))$. Namely, for
all $i, j$,
$$
[u+\gamma_{i+1}({\cal L}_{K}(G)), v+\gamma_{j+1}({\cal L}_{K}(G))]
= [u, v]+\gamma_{i+j+1}({\cal L}_{K}(G)),
$$
where $u \in \gamma_{i}({\cal L}_{K}(G))$, $v \in \gamma_{j}({\cal
L}_{K}(G))$ and $[u, v] \in \gamma_{i+j}({\cal L}_{K}(G))$.
Multiplication is then extended to ${\rm grad}^{(\ell)}({\cal
L}_{K}(G))$ by linearity.

Let $G$ be a torsion-free finitely generated nilpotent group. A
\emph{Mal'cev completion} of $G$ is a torsion-free radicable
nilpotent group $R$ containing $G$ and such that for all $a \in
R$, there exists $r \in {\mathbb N}$ such that $a^{r} \in G$. (A
group $G$ is said to be \emph{radicable} if for every $x$ in $G$
and an arbitrary natural number $m$, there exists $y$ in $G$ such
that $y^{m} = x$. If $G$ is a torsion-free nilpotent group, $x, y
\in G$, and $x^{m} = y^{m}$ for some $m \in {\mathbb N}$, then $x
= y$ (see \cite[page 247]{kurosh}.) If $G$ is a subgroup of a
torsion-free radicable nilpotent group $S$, then $S$ contains a
Mal'cev completion of $G$. Any two Mal'cev completions of $G$ are
isomorphic by an isomorphism which fixes $G$ pointwise (that is,
$R$ is unique up to isomorphism). We note that $\exp{\cal
L}_{\mathbb{Q}}(G)$ is a Mal'cev completion of $G$ (see
\cite{baumslag,segal}).

Let $L$ be a finitely generated nilpotent Lie algebra over
${\mathbb{Q}}$ and let $L^{\prime}$ denote the derived algebra of
$L$. Let $n = \dim L/L^{\prime}$, and let $h_{1}, \ldots, h_{n}$
be elements of $L$ such that the set $\{h_{1} + L^{\prime},
\ldots, h_{n} + L^{\prime}\}$ is a ${\mathbb{Q}}$-basis of
$L/L^{\prime}$. We assert that $L$ is generated by the set
$\{h_{1}, \ldots, h_{n}\}$ as Lie algebra. Indeed, let $\{y_{1},
\ldots, y_{m}\}$ be a generating set of $L$ and $m$ be the
smallest number of generators. It is easily shown that $m = n$.
Thus $\{y_{1} + L^{\prime}, \ldots, y_{n} + L^{\prime}\}$ is a
${\mathbb{Q}}$-basis of $L/L^{\prime}$. Let $A$ be the Lie
subalgebra of $L$ generated by the set $\{h_{1}, \ldots, h_{n}\}$.
To show that $A = L$ it is enough to prove that $L \subseteq A$.
For $j = 1, \ldots, n$,
$$
y_{j} = \sum^{n}_{i=1} \alpha_{ij} h_{i} + v_{j},
$$
where $\alpha_{ij} \in {\mathbb{Q}}$, $v_{j} \in L^{\prime}$.
Since $L$ is nilpotent, we can finally express each Lie commutator
$[y_{i_{1}}, \ldots, y_{i_{\kappa}}]$, $\kappa \in {\mathbb{N}}$,
as a ${\mathbb{Q}}$-linear combination of Lie commutators
$[h_{j_{1}}, \ldots, h_{j_{\ell}}]$, $\ell \in {\mathbb{Q}}$.
Therefore $L \subseteq A$ and so, $L = A$. Hence $L$ is generated
by the set $\{h_{1}, \ldots, h_{n}\}$. We give on $L$, by means of
BCH formula, the structure of a group, denoted by $R$. It is well
known that $R$ is a torsion-free radicable nilpotent group. Let
$H$ be the subgroup of $R$ generated by the set $\{h_{1}, \ldots,
h_{n}\}$. Then $R$ is a Mal'cev completion of $H$, and $H^{\prime}
= \tau_{2}(H)$ (see \cite[Proof of Theorem B, page 457]{bg}).
Since both $\exp {\cal L}_{\mathbb{Q}}(H)$ and $R$ are Mal'cev
completions of $H$, we obtain $R \cong \exp{\cal
L}_{\mathbb{Q}}(H)$ as groups (see, for example, \cite[Chapter 6,
Corollary 4]{segal}) and so, by the Mal'cev correspondence, $L
\cong {\cal L}_{\mathbb{Q}}(H)$ as Lie algebras by an isomorphism
sending $h_{i}$ to $\log h_{i}$, with $i = 1, \ldots, n$. We state
the above observations as lemma.

\begin{lemma}\label{2.1}
Let $L$ be a finitely generated nilpotent Lie algebra over
${\mathbb{Q}}$, and let $h_{1}, \ldots, h_{n}$ be elements of $L$
such that the set $\{h_{1} + L^{\prime}, \ldots, h_{n} +
L^{\prime}\}$ is a ${\mathbb{Q}}$-basis of $L/L^{\prime}$. Then
$L$ is generated by the set $\{h_{1}, \ldots, h_{n}\}$. Consider
$L$ as a group by means of the Baker-Campbell-Hausdorff formula,
denoted by $R$. Let $H$ be the subgroup of $R$ generated by the
set $\{h_{1}, \ldots, h_{n}\}$. Then $R$ is a Mal'cev completion
of $H$, $H$ is torsion-free finitely generated nilpotent group of
class $c$, $H^{\prime} = \tau_{2}(H)$, and $L \cong {\cal
L}_{\mathbb{Q}}(H)$ as Lie algebras. \qed
\end{lemma}

\section{Groups and Lie algebras}
\subsection{Some auxiliary Lemmas}

Throughout this section
we shall give some auxiliary results helping us to prove our main
results. Let $L$ be a Lie algebra. For a positive integer $i$, let
$\gamma_{i}(L)$ be the $i$-th term of the lower central series of
$L$. Write $L^{\prime} = \gamma_{2}(L)$. The following result is
well-known (and it is easily proved).

\begin{lemma}\label{3.1}
 Let $G$ be a torsion-free finitely
generated nilpotent group. Then \linebreak
$\gamma_{i}(G)\tau_{i+1}(G)/\tau_{i+1}(G)$ has finite index in
$\tau_{i}(G)/\tau_{i+1}(G)$. \qed
\end{lemma}

For a proof of the following result, we refer the reader to \cite[Chapter VIII, Lemma 9.4, page 330]{hb}.

\begin{lemma}\label{3.2}
Suppose that $G$ is a group and that
$L$ is a Lie ring. Suppose that for $i \geq 1$, $\sigma_{i}$ is a
homomorphism of $\gamma_{i}(G)/\gamma_{i+1}(G)$ onto an additive
subgroup $L_{i}$ of $L$ such that $L = L_{1} + L_{2} + \cdots$.
Suppose further that if $x \in \gamma_{i}(G)$, $y
\in\gamma_{j}(G)$,
$$
[\sigma_{i}(x \gamma_{i+1}(G)), \sigma_{j}(y \gamma_{j+1}(G))] =
\sigma_{i+j}((x, y) \gamma_{i+j+1}(G)),
$$
where $(x, y) = x^{-1}y^{-1}xy$.

i) If $G$ is generated by a set $X$, $L$ is the Lie ring
generated by the set
$$
Y = \{\sigma_{1}(x \gamma_{2}(G)): x \in X\}.
$$

ii) For $r = 1, 2, \ldots$, $\gamma_{r}(L) = L_{r} + L_{r+1}
+ \cdots$, and $\gamma_{r}(L)/\gamma_{r+1}(L)$ is isomorphic to
$L_{r}/(L_{r} \cap \gamma_{r+1}(L))$. \qed
\end{lemma}

Let $G$ be a torsion-free finitely generated nilpotent group of
class $c$. Let $a_{1}, \ldots, a_{n}$ be elements of $G$ such that
$\{a_{1} \tau_{2}(G), \ldots,$ $a_{n} \tau_{2}(G)\}$ is a
$\mathbb{Z}$-basis of $G/\tau_{2}(G)$. Since \linebreak
$(G/G^{\prime})/(\tau_{2}(G)/G^{\prime})$ is naturally isomorphic
to $G/\tau_{2}(G)$, and $G/\tau_{2}(G)$ is a free
$\mathbb{Z}$-module of rank $n$, there are elements $v_{1},
\ldots, v_{n} \in \tau_{2}(G)$ such that the set $\{a_{1}v_{1}
G^{\prime}, \ldots, a_{n}v_{n} G^{\prime}\}$ is a part of a
generating set of $G/G^{\prime}$. Let $y_{n+1}, \ldots, y_{m}$ be
elements of $\tau_{2}(G)$ subject to $\{a_{1}v_{1} G^{\prime},
\ldots, a_{n}v_{n} G^{\prime}, y_{n+1} G^{\prime}, \ldots, y_{m}
G^{\prime}\}$ is a generating set of $G/G^{\prime}$. Set $Y =
\{a_{1}v_{1}, \ldots, a_{n}v_{n},$ $y_{n+1}, \ldots, y_{m}\}$.
Since $G/G^{\prime}$ is finitely generated, we obtain
$\gamma_{i}(G)/\gamma_{i+1}(G)$ is finitely generated for all $i$.
In fact, the set of all group commutators of the form $(u_{1},
\ldots, u_{i})$, with $u_{1}, \ldots, u_{i} \in Y$, modulo
$\gamma_{i+1}(G)$, generates $\gamma_{i}(G)/\gamma_{i+1}(G)$. It
is easily verified that $(z_{1}, \ldots, z_{i-1}) \in
\tau_{i}(G)$, with $i \geq 2$, if some $z_{j} \in \{y_{n+1},
\ldots, y_{m}\}$ with $j = 1, \ldots, i-1$. Thus the set $Z_{i} =
\{(z_{1}, \ldots, z_{i}) \tau_{i+1}(G): z_{1}, \ldots, z_{i} \in
\{a_{1}, \ldots, a_{n}\}\}$ generates $\gamma_{i}(G)
\tau_{i+1}(G)/\tau_{i+1}(G)$. Write \linebreak $L^{(S)}(G) =
\oplus_{1 \leq i \leq c}L^{(S)}_{i}(G)$, with $L^{(S)}_{i}(G) =
\gamma_{i}(G)\tau_{i+1}(G)/\tau_{i+1}(G)$. It is easily checked
that $L^{(S)}(G)$ is a Lie subring of $L(G)$. For each positive
integer $i$, let $\pi_{i}$ be the natural mapping from
$\gamma_{i}(G)/\gamma_{i+1}(G)$ onto $L^{(S)}_{i}(G)$. Since, for
$x \in \gamma_{i}(G)$ and $y \in \gamma_{j}(G)$,
$$
[\pi_{i}(x \gamma_{i+1}(G)), \pi_{j}(y \gamma_{j+1}(G))] =
\pi_{i+j}((x, y) \gamma_{i+j+1}(G)),
$$
we obtain from Lemma \ref{3.2} that $L^{(S)}(G)$ is the Lie ring
generated by the set $\{a_{1} \tau_{2}(G), \ldots,$ $a_{n}
\tau_{2}(G)\}$. Since $L^{(S)}_{i}(G)$ is free abelian, we obtain
$K \otimes L^{(S)}_{i}(G)$ is a finite-dimensional vector space
over $K$. Since the set of all elements of the form $1 \otimes z$,
with $z \in Z_{i}$, spans $K \otimes L^{(S)}_{i}(G)$, we get there
exists a subset $X_{i}$ of $Z_{i}$ such that the set of all
elements of the form $1 \otimes x$, with $x \in X_{i}$, is a
$K$-basis of $K \otimes L^{(S)}_{i}(G)$. Since $L^{(S)}_{i}(G)$
has finite index in $\tau_{i}(G)/\tau_{i+1}(G)$ (Lemma \ref{3.1}),
and $\tau_{i}(G)/\tau_{i+1}(G)$ is free abelian, we obtain $K
\otimes L^{(S)}_{i}(G) = L_{i,K}(G)$. Hence it is easy to verify
the following result.

\begin{lemma}\label{3.3}
 Let $G$ be a torsion-free finitely
generated nilpotent group, and let $a_{1}, \ldots, a_{n}$ be
elements of $G$ such that $\{a_{1} \tau_{2}(G), \ldots,$ $a_{n}
\tau_{2}(G)\}$ is a $\mathbb{Z}$-basis of $G/\tau_{2}(G)$. Then
$\{a_{1} \tau_{2}(G), \ldots,$ \linebreak $a_{n} \tau_{2}(G)\}$ is
a generating set of $L_{K}(G)$, and $\{(a_{1} \tau_{2}(G)) +
L_{K}(G)^{\prime}, \ldots$ $, (a_{n} \tau_{2}(G)) +
L_{K}(G)^{\prime}\}$ is a $K$-basis for
$L_{K}(G)/L_{K}(G)^{\prime}$. \qed
\end{lemma}

The following result gives a generating set of ${\cal L}_{K}(G)$
as Lie algebra with respect to a given ${\mathbb Z}$-basis of
$G/\tau_{2}(G)$.

\begin{lemma}\label{3.4}
Let $G$ be a torsion-free finitely generated nilpotent group, and
let $a_{1}, \ldots, a_{n}$ be elements of $G$ such that $\{a_{1}
\tau_{2}(G), \ldots, a_{n} \tau_{2}(G)\}$ is a ${\mathbb Z}$-basis of $G/\tau_{2}(G)$. Then $\{\log a_{1}, \ldots,$
\linebreak $\log a_{n}\}$ generates ${\cal L}_{K}(G)$ as Lie
algebra.
\end{lemma}

\pf Let $a_{1}, \ldots, a_{n}$ be elements of $G$ such that
$\{a_{1} \tau_{2}(G), \ldots, a_{n} \tau_{2}(G)\}$ is a ${\mathbb
Z}$-basis of $G/\tau_{2}(G)$. Choose a canonical basis $\{a_{1},
\ldots, a_{n}, \ldots, a_{f(c+1)}\}$ of $G$ where $c$ is the
nilpotency class of $G$. Let ${\cal A}$ be the Lie subalgebra of
${\cal L}_{K}(G)$ generated by the set $\{\log a_{1}, \ldots, \log
a_{n}\}$. We claim that ${\cal A} = {\cal L}_{K}(G)$. It is enough
to show that ${\cal L}_{K}(G) \subseteq {\cal A}$. In fact it is
enough to show that $\log g \in {\cal A}$ for all $g \in G$. Since
${\cal A}$ is a Lie algebra and $\{\log a_{1}, \ldots, \log
a_{f(c+1)}\}$, is a $K$-basis of ${\cal L}_{K}(G)$ (as vector
space), it is enough to show that $\log g \in {\cal A}$ for all $g
\in \{a_{1}, \ldots, a_{f(c+1)}\}$. Let $y_{n+1}, \ldots y_{m}$ be
elements of $\tau_{2}(G)$ subject to $\{a_{1}v_{1} G^{\prime},
\ldots, a_{n}v_{n} G^{\prime}, y_{n+1} G^{\prime}, \ldots,$ $y_{m}
G^{\prime}\}$, with $v_{1}, \ldots, v_{n} \in \tau_{2}(G)$, is a
generating set for $G/G^{\prime}$. The set $\{(z_{1}, \ldots,
z_{i}) \tau_{i+1}(G): z_{1}, \ldots, z_{i} \in \{a_{1}, \ldots,
a_{n}\}\}$ generates $\gamma_{i}(G) \tau_{i+1}(G)/\tau_{i+1}(G)$.
Since $\gamma_{i}(G) \tau_{i+1}(G)/\tau_{i+1}(G)$ has finite index
in $\tau_{i}(G)/\tau_{i+1}(G)$ (by Lemma \ref{3.1}), we obtain
$a^{m_{ij}}_{f(i)+j} = g_{ij} u_{i+1,j}$ for some $m_{ij} \in
{\mathbb Z}$, $g_{ij}$ is a product of group commutators of the
form $(z_{1}, \ldots, z_{i})$ and $z_{1}, \ldots, z_{i} \in
\{a_{1}, \ldots, a_{n}\}$, and $u_{i+1,j} \in \tau_{i+1}(G)$ $j =
1, \ldots, n_{i}$. For $i = c$, we have
\begin{eqnarray}
a^{m_{ij}}_{f(c)+j} = g_{cj}
\end{eqnarray}
for $j = 1, \ldots, n_{c}$. By the BCH formula, we have $\log g
\in {\cal A}$ for all $g \in \gamma_{c}(G)$. Applying the BCH
formula on (4) and since $\log g_{cj} \in {\cal A}$, we obtain
$\log a_{f(c)+1}, \ldots, \log a_{f(c+1)} \in {\cal A}$, and hence
$\log g \in {\cal A}$ for all $g \in \tau_{c}(G)$. Suppose that
for all $\kappa$, with $i < \kappa \leq c$, $\log a_{f(\kappa)+1},
\ldots, \log a_{f(c+1)} \in {\cal A}$. Let $h \in \tau_{k}(G)$.
Since $G$ is nilpotent and $\cap_{j \geq 1} \tau_{j}(G) = \{1\}$,
we write $h = g_{\kappa} g_{\kappa + 1} \cdots g_{c}$ with $g_{j}
\in \tau_{j}(G) \setminus \tau_{j+1}(G)$ and $j = \kappa, \ldots,
c$. Note that
\begin{eqnarray}
h = a^{\beta_{f(\kappa)+1}}_{f(\kappa)+1} \cdots
a^{\beta_{f(\kappa+1)}}_{f(\kappa+1)} \cdots
a^{\beta_{f(c)+1}}_{f(c)+1} \cdots a^{\beta_{f(c+1)}}_{f(c+1)}.
\end{eqnarray}
Applying the BCH formula on (5), our hypothesis and since ${\cal
A}$ is a Lie algebra, we have $\log h \in {\cal A}$ for all $h \in
\tau_{\kappa}(G)$. In particular, for $\kappa = i + 1$,
\begin{eqnarray}
\log h \in {\cal A}~~{\rm for~~all}~~h \in \tau_{i+1}(G).
\end{eqnarray}
For $z_{1}, \ldots, z_{i} \in \{a_{1}, \ldots, a_{n}\}$, the BCH
formula gives
\begin{eqnarray}
\log (z_{1}, \ldots, z_{i}) = [\log z_{1}, \ldots, \log z_{i}] +
\sum_{\mu} r_{\mu} d_{\mu}
\end{eqnarray}
where each $d_{\mu}$ is a repeated Lie commutator of length at
least $i+1$ in the arguments $\log z_{1}, \ldots, \log z_{i}$ each
of which appears at least once; and the coefficients $r_{\mu} \in
{\mathbb{Q}}$ (see \cite[Corollary 2, page 102]{segal}. The
arguments given in the proof of Corollary 2 can be adopted here.)
Once more, applying the BCH formula on a product of group
commutators of length $i$, using the equation (7) and the fact that ${\cal A}$
is a Lie algebra generated by the set $\{\log a_{1}, \ldots, \log
a_{n}\}$, we obtain
\begin{eqnarray}
\log (\prod_{{\rm finite}} (z_{1}, \ldots, z_{i})^{\nu}) \in {\cal
A},
\end{eqnarray}
where $\nu \in {\mathbb Z}$. But $a^{m_{ij}}_{f(i)+j} =
g_{ij}u_{(i+1)j}$ for some $m_{ij} \in {\mathbb Z}$, $g_{ij}$
is a product of group commutators of the form $(z_{1}, \ldots,
z_{i})$, with $z_{1}, \ldots, z_{i} \in \{a_{1}, \ldots, a_{n}\}$,
and $u_{(i+1)j} \in \tau_{i+1}(G)$ $j = 1, \ldots, n_{i}$. Apply
the BCH formula on $a^{m_{ij}}_{f(i)+j} = g_{ij}u_{(i+1)j}$. By
the equations (7) and (8) and since ${\cal A}$ is a Lie algebra,
we obtain $\log a_{f(i)+1}, \ldots, \log a_{f(i+1)} \in {\cal A}$.
Thus ${\cal L}_{K}(G) \subseteq {\cal A}$ and so ${\cal L}_{K}(G)
= {\cal A}$. \qed

\vskip .120 in

From the proof of Lemma \ref{3.4}, we obtain the following result.

\begin{corollary}\label{3.5}
Let $G$ be a torsion-free finitely generated nilpotent group of
class $c$, and let $a_{1}, \ldots, a_{n}$ be elements of $G$ such
that $\{a_{1} \tau_{2}(G), \ldots, a_{n} \tau_{2}(G)\}$ is a
${\mathbb Z}$-basis of $G/\tau_{2}(G)$. Let $\{a_{1}, a_{2},
\ldots, a_{f(c+1)}\}$ be a canonical basis of $G$. Then for $j
\geq f(t)+1$, with $t \geq 2$, $\log a_{j} \in \gamma_{t}({\cal
L}_{K}(G))$. \qed
\end{corollary}

\begin{lemma}\label{3.6}
Let $G$ be a torsion-free finitely
generated nilpotent group of class $c$. For a positive integer
$t$, with $2 \leq t \leq c$, let $\tau_{t}(G)$ be the isolator of
$\gamma_{t}(G)$ in $G$. Let $\pi_{t}$ be the natural mapping from
$G$ onto $G/\tau_{t}(G)$. Then there exists a Lie algebra
epimorphism $\xi_{\pi_{t}}$ from ${\cal L}_{K}(G)$ onto ${\cal
L}_{K}(G/\tau_{t}(G))$ such that, for all $g \in G$,
$\xi_{\pi_{t}}(\log g) = \log \pi_{t}(g)$. Furthermore ${\rm
ker}\xi_{\pi_{t}} = {\cal L}_{K}(\tau_{t}(G))$.
\end{lemma}

\pf For $H \in \{G, G/\tau_{t}(G)\}$, we write
$\Gamma_{H} = KH$ for the group algebra of $H$ over $K$, and
$\Delta_{H}$ for the augmentation ideal of $\Gamma_{H}$. Notice
that $H$ is a torsion-free finitely generated nilpotent group. Let
$\Gamma^{*}_{H}$ be the algebra of all formal power series $a^{*}$
of the form $a^{*} = \alpha_{0} + \sum \alpha_{k} d_{k}$, where
$\alpha_{0}, \alpha_{k} \in K$, $k = 1, 2, \ldots,$ and $d_{k} \in
\Delta^{k}_{H}$. In addition, we write $\Delta^{*}_{H}$ for the
subalgebra of $\Gamma^{*}_{H}$ consisting of all elements $a^{*}$
with $\alpha_{0} = 0$. We associate with $\Delta^{*}_{H}$ the Lie
algebra $\Lambda^{*}_{H} = (\Delta^{*}_{H})_{L}$. Let $h \in H$.
The natural mapping $\pi_{t}$ induces an algebra epimorphism
$\pi^{*}_{t,K}$ from $\Delta^{*}_{G}$ onto
$\Delta^{*}_{G/\tau_{t}(G)}$ such that $\pi^{*}_{t,K}(\log g) =
\log \pi_{t}(g)$ for all $g \in G$. It is easily checked that
$\pi^{*}_{t,K}$ induces a Lie algebra epimorphism from
$\Lambda^{*}_{G}$ onto $\Lambda^{*}_{G/\tau_{t}(G)}$. Hence we
obtain a Lie algebra epimorphism $\xi_{\pi_{t}}$ from ${\cal
L}_{K}(G)$ onto ${\cal L}_{K}(G/\tau_{t}(G))$ such that
$\xi_{\pi_{t}}(\log g) = \log \pi_{t}(g)$ for all $g \in G$. But
$$
\begin{array}{ccl}
\dim{\rm ker}\xi_{\pi_{t}} & = & \dim{\cal L}_{K}(G) - \dim{\cal L}_{K}(G/\tau_{t}(G)) \\
 & = & {\cal H}(G) - {\cal H}(G/\tau_{t}(G)) \\
 & = & {\cal H}(\tau_{t}(G)) \\
 & = & \dim{\cal L}_{K}(\tau_{t}(G)).
\end{array}
$$
Let $v \in {\cal L}_{K}(\tau_{t}(G))$. Then $v = k_{1} \log v_{1}
+ \cdots + k_{s} \log v_{s}$ with $v_{1}, \ldots, v_{s} \in
\tau_{t}(G)$ and $k_{1}, \ldots, k_{s} \in K$. Since
$\xi_{\pi_{t}}(v) = k_{1} \log \pi_{t}(v_{1}) + \cdots + k_{s}
\log \pi_{t}(v_{s}) = 0$, we obtain $v \in {\rm
ker}\xi_{\pi_{t}}$. Therefore ${\cal L}_{K}(\tau_{t}(G)) \subseteq
{\rm ker} \xi_{\pi_{t}}$ and so, ${\rm ker} \xi_{\pi_{t}} = {\cal
L}_{K}(\tau_{t}(G))$. \qed

\vskip .120 in

The following result gives us information about the terms of the
lower central series of ${\cal L}_{K}(G)$.

\vskip .120 in

\begin{lemma}\label{3.7}
Let $G$ be a torsion-free  finitely generated nilpotent group of
class $c$. Then $\gamma_{t}({\cal L}_{K}(G)) = {\cal
L}_{K}(\tau_{t}(G))$ for all $t$, with $t = 1, \ldots, c$.
\end{lemma}

\pf Since ${\cal L}_{K}(\tau_{t}(G))$ has a $K$-basis the set
$\{\log a_{f(t)+1}, \ldots, \log a_{f(c+1)}\}$, and by Corollary
\ref{3.5}, it is enough to show that $\gamma_{t}({\cal L}_{K}(G))
\subseteq {\cal L}_{K}(\tau_{t}(G))$. By Lemma \ref{3.6}, ${\rm
ker}\xi_{\pi_{t}} = {\cal L}_{K}(\tau_{t}(G))$ and so, it is
enough to show that $\gamma_{t}({\cal L}_{K}(G)) \subseteq {\rm
ker} \xi_{\pi_{t}}$. Let $v \in \gamma_{t}({\cal L}_{K}(G))$. By
Lemma \ref{3.4}, $v$ is written as a $K$-linear combination of Lie
commutators of the form $[\log a_{i_{1}}, \ldots, \log
a_{i_{\kappa}}]$ with $t \leq \kappa \leq c$. The BCH formula
gives
\begin{eqnarray}
[\log a_{i_{1}}, \ldots, \log a_{i_{\kappa}}] = \log (a_{i_{1}},
\ldots, a_{i_{\kappa}}) + \sum_{i} s_{i} \log v_{i},
\end{eqnarray}
where each $v_{i}$ is a left normed group commutator of length at
least $\kappa+1$ in the arguments $a_{i_{1}}, \ldots,
a_{i_{\kappa}}$, each of which appears at least once; and the
coefficients belong to ${\mathbb Q}$. Applying
$\xi_{\pi_{t}}$ to the equation (9), we get $[\log a_{i_{1}},
\ldots, \log a_{i_{\kappa}}] \in {\rm ker}\xi_{\pi_{t}} = {\cal
L}_{K}(\tau_{t}(G))$. Therefore $\gamma_{t}({\cal L}_{K}(G))
\subseteq {\rm ker}\xi_{\pi_{t}}$ and so, $\gamma_{t}({\cal
L}_{K}(G)) = {\cal L}_{K}(\tau_{t}(G))$. \qed

\vskip .120 in

The following result is probably well-known.

\begin{lemma}\label{3.8}
For any torsion-free finitely generated nilpotent group $G$ of
class $c$, ${\rm grad}^{(\ell)}({\cal L}_{K}(G))$ is isomorphic to
$L_{K}(G)$ as Lie algebras under the mapping $\phi$ sending $\log
a_{f(i)+j} + \gamma_{i+1}({\cal L}_{K}(G))$ to $a_{f(i)+j}\tau_{i+1}(G)$
for all $i = 1, \ldots, c$, $j = 1, \ldots, n_{i}$.
\end{lemma}

\pf By the proof of Lemma \ref{3.7}, the set $\{\log a_{f(t)+1},
\ldots, \log a_{f(c+1)}\}$ is a $K$-basis of \linebreak
$\gamma_{t}({\cal L}_{K}(G)) = {\cal L}_{K}(\tau_{t}(G))$. Notice
that
$$
\gamma_{t}({\cal L}_{K}(G)) = {\rm span} \{\log a_{f(t)+1},
\ldots, \log a_{f(t+1)}\} \oplus \gamma_{t+1}({\cal L}_{K}(G)).
$$
Since
$$
\gamma_{t}({\cal L}_{K}(G))/\gamma_{t+1}({\cal L}_{K}(G)) \cong K
\otimes \tau_{t}(G)/\tau_{t+1}(G)
$$
as vector spaces for $t = 1, \ldots, c$, we obtain the mapping
$\phi$ from ${\rm grad}^{(\ell)}({\cal L}_{K}(G))$ to $L_{K}(G)$
sending $\log a_{f(i)+j} + \gamma_{i+1}({\cal L}_{K}(G))$ to
$a_{f(i)+j}\tau_{i+1}(G)$ for all $i = 1, \ldots, c$, $j = 1,
\ldots, n_{i}$, is a $K$-linear isomorphism. By the BCH formula,
we obtain
\begin{flushleft}
$\phi([\log a_{f(i)+\mu}+\gamma_{i+1}({\cal L}_{K}(G)), \log
a_{f(j)+\nu}+\gamma_{j+1}({\cal L}_{K}(G))]) = \phi([\log
a_{f(i)+\mu}, \log a_{f(j)+\nu}]+\gamma_{i+j+1}({\cal L}_{K}(G)))
= \phi(\log (a_{f(i)+\mu}, a_{f(j)+\nu}) + \gamma_{i+j+1}({\cal
L}_{K}(G)).$
\end{flushleft}
Write
$$
(a_{f(i)+\mu}, a_{f(j)+\nu}) = a^{m_{1}}_{f(i+j)+1} \cdots
a^{m_{n_{i+j}}}_{f(i+j+1)} v,
$$
where $m_{1}, \ldots, m_{n_{i+j}} \in {\mathbb{Z}}$, $v \in
\tau_{i+j+1}(G)$. Then, by BCH formula,
\begin{flushleft}
$\phi([\log a_{f(i)+\mu}+\gamma_{i+1}({\cal L}_{K}(G)), \log
a_{f(j)+\nu}+\gamma_{j+1}({\cal L}_{K}(G))]) = \phi(m_{1} \log
a_{f(i+j)+1} + \cdots + m_{n_{i+j}} \log a_{f(i+j+1)} +
\gamma_{i+j+1}({\cal L}_{K}(G))) =$ $m_{1} a_{f(i+j)+1}
\tau_{i+j+1}(G) + \cdots + m_{n_{i+j}} a_{f(i+j+1)}
\tau_{i+j+1}(G) = a^{m_{1}}_{f(i+j)+1} \cdots
a^{m_{n_{i+j}}}_{f(i+j+1)} \tau_{i+j+1}(G)$ $= (a_{f(i)+\mu},
a_{f(j)+\nu}) \tau_{i+j+1}(G)= [\phi(\log
a_{f(i)+\mu}+\gamma_{i+1}({\cal L}_{K}(G))), \phi(\log
a_{f(j)+\nu}+\gamma_{j+1}({\cal L}_{K}(G)))])$
\end{flushleft}
for all $i, j \in \{1, \ldots, c\}$, $\mu \in \{1, \ldots,
n_{i}\}$ and $\nu \in \{1, \ldots, n_{j}\}$. Thus $\phi$ is a Lie
algebra isomorphism. \qed

\subsection{Relatively free groups}

 Let
${\mathfrak{N}}_{c}$ be the variety of nilpotent groups of class
at most $c$, and let ${\mathfrak{T}}_{c}$ be a torsion-free
subvariety of class at most $c$ of ${\mathfrak{N}}_{c}$. For the
rest of this section, for positive integers $n$ and $c$, with $n
\geq 2$, we write $F_{n,c} = F_{n}({\mathfrak{N}}_{c})$ and $G =
F_{n}({\mathfrak{T}}_{c})$. The groups $F_{n,c}$ and $G$ are
freely generated by the set $\{x_{1}, \ldots, x_{n}\}$, $x_{i} =
f_{i} {\mathfrak{N}}_{c}(F_{n})$, $i = 1, \ldots, n$, and the set
$\{y_{1}, \ldots, y_{n}\}$, $y_{i} = f_{i}
{\mathfrak{T}}_{c}(F_{n})$, $i = 1, \ldots, n$, respectively.
Since ${\mathfrak{T}}_{c}$ is a subvariety of
${\mathfrak{N}}_{c}$, the natural map $\pi$ from $F_{n,c}$ to $G$
is surjective. Since ${\rm ker} \pi =
{\mathfrak{T}}_{c}(F_{n})/{\mathfrak{N}}_{c}(F_{n})$, we obtain
${\rm ker} \pi$ is a fully invariant subgroup of $F_{n,c}$. Hence
${\rm ker} \pi$ is a verbal subgroup of $F_{n,c}$. Let $\mu$ and
$\nu$ be the Hirsch numbers of $F_{n,c}$ and $G$, respectively.
Since $F_{n,c}$ is finitely generated nilpotent, ${\rm ker}\pi$ is
a (torsion-free) finitely generated nilpotent group of class at
most $c$. Since $F_{n,c}/{\rm ker}\pi \cong G$ and $G$ is
torsion-free, it follows from a result of Hirsch (see
\cite[Theorem 2.312]{hirsch}, also, \cite[Theorem 2.3]{jennings})
that ${\rm ker}\pi$ has Hirsch number $\mu - \nu$. Since
$L_{K}(F_{n,c})$ is a free nilpotent Lie algebra of class $c$
freely generated by the set $\{x_{1} F^{\prime}_{n,c}, \ldots,
x_{n} F^{\prime}_{n,c}\}$ (see, for example, \cite{shmelkin}),
${\cal L}_{K}(F_{n,c})$ is a nilpotent Lie algebra of class $c$
and $\{\log x_{1}, \ldots, \log x_{n}\}$ generates ${\cal
L}_{K}(F_{n,c})$ (by Lemma \ref{3.4}), we obtain the mapping
$\chi$ from the set $\{x_{1} F^{\prime}_{n,c}, \ldots, x_{n}
F^{\prime}_{n,c}\}$ into ${\cal L}_{K}(F_{n,c})$ sending $x_{i}
F^{\prime}_{n,c}$ to $\log x_{i}$, $i = 1, \ldots, n$, extends
uniquely to a Lie algebra epimorphism. Since $\dim L_{K}(F_{n,c})
= \dim {\cal L}_{K}(F_{n,c}) = \mu$, we obtain $\chi$ is
one-to-one and so, $\chi$ is a Lie algebra isomorphism. We
summarize the above observations as follows.

\begin{lemma}\label{3.9}
For positive integers $n$ and $c$, with $n \geq 2$, let $F_{n,c}$
be the free nilpotent group of rank $n$ and class $c$; freely
generated by the set $\{x_{1}, \ldots, x_{n}\}$. Then ${\cal
L}_{K}(F_{n,c})$ is a free nilpotent Lie algebra of rank $n$ and
class $c$; freely generated by the set $\{\log x_{1}, \ldots, \log
x_{n}\}$. \qed
\end{lemma}

The next result gives us a way of how a group homomorphism of
$F_{n,c}$ onto $G$ can define a Lie algebra homomorphism from
${\cal L}_{K}(F_{n,c})$ onto ${\cal L}_{K}(G)$. The proof of the
following result is similar to the proof given in Lemma \ref{3.6}.

\begin{lemma}\label{3.10}
Let $M$ be a torsion-free finitely generated nilpotent group of
class $c$ such that $M/\tau_{2}(M)$ is a free abelian group of
rank $n$, with $n \geq 2$. Let $\varphi$ be any group homomorphism
from $F_{n,c}$ into $M$, and let $\psi_{\varphi}$ be the mapping
from the set $\{\log x_{1}, \ldots, \log x_{n}\}$ into ${\cal
L}_{K}(M)$ defined by $\psi_{\varphi}(\log x_{i}) = \log
\varphi(x_{i})$ for $i = 1, \ldots, n$. Then $\psi_{\varphi}$
extends uniquely to a Lie algebra homomorphism from ${\cal
L}_{K}(F_{n,c})$ into ${\cal L}_{K}(M)$ and, for all $u \in
F_{n,c}$, $\psi_{\varphi}(\log u) = \log \varphi(u)$. \qed
\end{lemma}

\begin{lemma}\label{3.11}
For positive integers $n$ and $c$, with
$n \geq 2$, let $G = F_{n}({\mathfrak{T}}_{c})$, with $n \geq 2$,
and let $\pi$ be the natural mapping from $F_{n,c}$ onto $G$. Then
${\rm ker}\pi \subseteq F^{\prime}_{n,c}$, and ${\cal L}_{K}({\rm
ker}\pi) \subseteq {\cal L}_{K}(F^{\prime}_{n,c}) = {\cal
L}_{K}(F_{n,c})^{\prime}$.
\end{lemma}

\pf First we shall show that ${\rm ker}\pi \subseteq
F^{\prime}_{n,c}$. To get a contradiction, we assume that ${\rm
ker}\pi \nsubseteq F^{\prime}_{n,c}$. Thus $(F_{n,c}/{\rm
ker}\pi)^{\prime} = F^{\prime}_{n,c} {\rm ker}\pi/{\rm ker}\pi$.
Since $G/G^{\prime}$ is free abelian of rank $n$, and
$G/G^{\prime} \cong F_{n,c}/F^{\prime}_{n,c} {\rm ker}\pi$, we get
$F_{n,c}/F^{\prime}_{n,c}{\rm ker}\pi$ is free abelian of rank
$n$. Since $F_{n,c}/F^{\prime}_{n,c}$ is free abelian of rank $n$,
we obtain $F^{\prime}_{n,c} {\rm ker}\pi/F^{\prime}_{n,c}$ is not
a trivial free abelian subgroup of $F_{n,c}/F^{\prime}_{n,c}$
which is a contradiction. Therefore ${\rm ker}\pi \subseteq
F^{\prime}_{n,c}$. Thus ${\cal L}({\rm ker}\pi) \subseteq {\cal
L}_{K}(F^{\prime}_{n,c})$. By Lemma \ref{3.7}, we obtain ${\cal
L}_{K}(F^{\prime}_{n,c}) = {\cal L}_{K}(F_{n,c})^{\prime}$. Hence
${\cal L}_{K}({\rm ker}\pi) \subseteq {\cal
L}_{K}(F_{n,c})^{\prime}$. \qed

For a torsion-free finitely generated nilpotent group $G$ of class
$c$, we write ${\cal L}_{t}(G)$ for the vector subspace of ${\cal
L}_{K}(G)$ spanned by all Lie commutators of the form $[\log
a_{i_{1}}, \ldots, \log a_{i_{t}}]$ with $i_{1}, \ldots, i_{t} \in
\{1, \ldots, n\}$.

\begin{lemma}\label{3.12} For positive integers $n$ and $c$, with
$n \geq 2$, let $G = F_{n}({\mathfrak{T}}_{c})$ and let $\pi$ be
the natural mapping from $F_{n,c}$ onto $G$. Then
\[
{\cal L}_{K}({\rm ker}\pi) = \oplus^{c}_{t = 2} ({\cal L}_{K}({\rm
ker}\pi) \cap {\cal L}_{t}(F_{n,c})).
\]
\end{lemma}

\pf By Lemma \ref{3.11}, and since ${\cal
L}_{K}(F_{n,c})=\oplus_{t=1}^c{\cal L}_{K}(F_{n,c})$, we have
$$
{\cal L}_{K}({\rm ker}\pi) \subseteq {\cal
L}_{K}(F_{n,c})^{\prime} = \oplus^{c}_{t = 2} {\cal
L}_{t}(F_{n,c}).
$$
Let $w \in {\cal L}_{K}({\rm ker}\pi)$. For a fixed $i$, with $i =
1, \ldots, n$, express $w$ as a sum
$$
w = w_{0} + w_{1} + \cdots + w_{s}
$$
where the number of times that $\log x_{i}$ occurs in $w_{j}$ is
$j$, $j = 0, \ldots, s$. Pick distinct non-zero integer numbers
$\alpha_{0}, \alpha_{1}, \ldots, \alpha_{s}$. Then
$$
w(\log x_{1}, \ldots, \alpha_{j} \log x_{i}, \ldots, \log x_{n}) =
\sum^{s}_{\kappa = 0} \alpha^{\kappa}_{j} w_{\kappa}(\log x_{1},
\ldots, \log x_{n}).
$$
Now the determinant
$$
 \left |
 \begin{array}{ccccc}
 1 & \alpha_{0} & \alpha^{2}_{0} & \cdots & \alpha^{s}_{0} \\
 1 & \alpha_{1} & \alpha^{2}_{1} & \cdots & \alpha^{s}_{1} \\
 \vdots & \vdots & \vdots & & \vdots \\
 1 & \alpha_{s} & \alpha^{2}_{s} & \cdots & \alpha^{s}_{s}
 \end{array}
 \right |,
$$
is a Vandermonte determinant with value $\prod_{i < j} (\alpha_{i}
- \alpha_{j})$ so is non-zero. Consequently, each $w_{j}$ is a
${\mathbb Q}$-linear combination of the elements $w(\log
x_{1}, \ldots, \alpha_{j} \log x_{i}, \ldots, \log x_{n})$. For
the fixed $i$ and $j = 0, \ldots, s$, let $\varphi_{i,j}$ be the
mapping from the set $\{x_{1}, \ldots, x_{n}\}$ into $F_{n,c}$
defined by $\varphi_{i,j}(x_{\kappa}) = x^{\alpha_{j}}_{i}$ if $i
= \kappa$ and $\varphi_{i,j}(x_{\kappa}) = x_{\kappa}$ if $i \neq
\kappa$. Since $F_{n,c}$ is relatively free on the set $\{x_{1},
\ldots, x_{n}\}$, $\varphi_{i,j}$ extends uniquely to a group
endomorphism of $F_{n,c}$. By Lemma \ref{3.10} (for $M = F_{n,c}$),
$\psi_{\varphi_{i,j}}$ is a Lie algebra endomorphism of ${\cal
L}_{K}(F_{n,c})$. Write
$$
w = \lambda_{1} \log u_{1} + \cdots + \lambda_{s} \log u_{s},
$$
where $\lambda_{1}, \ldots, \lambda_{s} \in K$ and $u_{1}, \ldots,
u_{s} \in {\rm ker} \pi$. Applying $\psi_{\varphi_{i,j}}$ on $w$,
and since ${\rm ker} \pi$ is a fully invariant subgroup of
$F_{n,c}$, we obtain $\psi_{\varphi_{i,j}}(w) \in {\cal
L}_{K}({\rm ker}\pi)$. Hence
$$
w(\log x_{1}, \ldots, \alpha_{j} \log x_{i}, \ldots, \log x_{n})
\in {\cal L}_{K}({\rm ker}\pi)
$$
for $j = 0, \ldots, s$. Therefore, for $\kappa = 0, \ldots, s$, we
have
$$
w_{\kappa}(\log x_{1}, \ldots, \log x_{n}) \in {\cal L}_{K}({\rm
ker}\pi).
$$
Consequently, each $w_{i}$ is a ${\mathbb Q}$-linear
combination of the elements
$$
w(\log x_{1}, \ldots, \alpha_{j} \log x_{i}, \ldots, \log x_{n}),
$$
for $j = 0, \ldots, s$. If we repeat the process on each $w_{i}$
using different $\log x_{\kappa}$ (with $\kappa \neq i$), then
eventually we obtain each homogeneous component of $w$ belongs to
${\cal L}_{K}({\rm ker}\pi)$. Therefore if $w \in {\cal
L}_{K}({\rm ker}\pi)$ then the homogeneous components of $w$
belong to ${\cal L}_{K}({\rm ker}\pi)$ and so, we obtain the
required result. \qed

\begin{proposition}\label{3.13}
For all positive integers $n$ and $c$, with $n \geq 2$, ${\cal
L}_{K}(F_{n,c}/{\rm ker}\pi)$ is isomorphic to ${\cal
L}_{K}(F_{n,c})/{\cal L}_{K}({\rm ker} \pi)$ as a Lie algebra.
Furthermore, ${\cal L}_{\mathbb{Q}}({\rm ker}\pi)$ is a fully
invariant ideal of ${\cal L}_{\mathbb{Q}}(F_{n,c})$.
\end{proposition}

\pf By Lemma \ref{3.10} (for $M = G$), the mapping $\psi_{\pi}$
from the set $\{\log x_{1}, \ldots, \log x_{n}\}$ into ${\cal
L}_{K}(G)$ defined by $\psi_{\pi}(\log x_{i}) = \log \pi(x_{i})$
for $i = 1, \ldots, n$ extends uniquely to a Lie algebra
homomorphism from ${\cal L}_{K}(F_{n,c})$ into ${\cal L}_{K}(G)$
and, for all $u \in F_{n,c}$, $\psi_{\pi}(\log u) = \log \pi(u)$.
Since $\pi$ is onto and by the definition of ${\cal L}_{K}(G)$, we
obtain $\psi_{\pi}$ is onto. First we shall show that ${\rm ker}
\psi_{\pi} = {\cal L}_{K}({\rm ker} \pi)$. Let $w \in {\cal
L}_{K}({\rm ker} \pi)$. Then
$$
w = \kappa_{1} \log v_{1} + \cdots + \kappa_{s} \log v_{s},
$$
where $v_{1}, \ldots, v_{s} \in {\rm ker} \pi$ and $\kappa_{1},
\ldots, \kappa_{s} \in K$. Applying $\psi_{\pi}$, we have
$$
\psi_{\pi}(w) = \kappa_{1} \log \pi(v_{1}) + \cdots + \kappa_{s}
\log \pi(v_{s}) = 0.
$$
Therefore $w \in {\rm ker} \psi_{\pi}$. Hence
\begin{eqnarray}
{\cal L}_{K}({\rm ker} \pi) \subseteq {\rm ker} \psi_{\pi}.
\end{eqnarray}
But $\dim {\cal L}_{K}({\rm ker}\pi) = {\cal H}({\rm ker}\pi) =
\mu - \nu$. In addition, $\dim ({\rm ker}\psi_{\pi}) = \dim {\cal
L}_{K}(F_{n,c}) - \dim {\cal L}_{K}(G) = {\cal H}(F_{n,c}) - {\cal
H}(G) = \mu - \nu$. By the equation (10), we obtain ${\rm
ker}\psi_{\pi} = {\cal L}_{K}({\rm ker} \pi)$.

Next we show that ${\cal L}_{\mathbb{Q}}({\rm ker}\pi)$ is a fully
invariant ideal of ${\cal L}_{\mathbb{Q}}(F_{n,c})$. Let $\xi$ be
a Lie algebra homomorphism of ${\cal L}_{\mathbb{Q}}(F_{n,c})$,
and let $\xi(\log x_{i}) = u_{i}$ for $i = 1, \ldots, n$. Since
$\exp {\cal L}_{\mathbb{Q}}(F_{n,c})$ is a Mal'cev completion of
$F_{n,c}$ there are positive integers $m_{i}$ such that $\exp
(m_{i} u_{i}) \in F_{n,c}$ for $i = 1, \ldots, n$. Let $\varphi$
be the endomorphism of $F_{n,c}$ satisfying the conditions
$\varphi(x_{i}) = \exp (m u_{i})$, with $m$ the least common
multiple of $m_{1}, \ldots, m_{n}$ for $i = 1, \ldots, n$. (Notice
that $\exp(mu_{i}) = 1$ if and only if $u_{i} = 0$.) Since ${\rm
ker}\pi$ is a fully invariant subgroup of $F_{n,c}$, we have
$\varphi(v) \in {\rm ker}\pi$ for all $v \in {\rm ker}\pi$. Let $w
\in {\cal L}_{\mathbb Q}({\rm ker}\pi)$ and, by Lemma \ref{3.12},
we assume that $w$ is a homogeneous element. Write
$$
w = w(\log x_{1}, \ldots, \log x_{n}) = \lambda_{1} \log v_{1} +
\cdots + \lambda_{s} \log v_{s},
$$
where $\lambda_{1}, \ldots, \lambda_{s} \in {\mathbb Q}$ and
$v_{1}, \ldots, v_{s} \in {\rm ker}\pi$. By Lemma \ref{3.10} (for $M =
F_{n,c}$), $\psi_{\varphi}(\log u) = \log \varphi(u)$ for all $u
\in F_{n,c}$. Applying $\psi_{\varphi}$ on $w$, we obtain
$\psi_{\varphi}(w) \in {\cal L}_{\mathbb{Q}}({\rm ker}\pi)$. Since
$\psi_{\varphi}$ is a Lie algebra homomorphism of ${\cal
L}_{\mathbb{Q}}(F_{n,c})$, we obtain $\lambda w(u_{1}, \ldots,
u_{n}) \in {\cal L}_{\mathbb{Q}}({\rm ker}\pi)$ for some non-zero
integer $\lambda$. Hence $\xi(w) \in {\cal L}_{\mathbb{Q}}({\rm
ker}\pi)$ for all homogeneous elements $w \in {\cal
L}_{\mathbb{Q}}({\rm ker}\pi)$. Since for any element $w$ in
${\cal L}_{\mathbb{Q}}({\rm ker}\pi)$ the homogeneous components
of $w$ belong to ${\cal L}_{\mathbb{Q}}({\rm ker}\pi)$ (by Lemma
\ref{3.12}) and $\xi$ is a Lie algebra homomorphism, we obtain ${\cal
L}_{\mathbb{Q}}({\rm ker}\pi)$ is a fully invariant ideal of
${\cal L}_{\mathbb{Q}}(F_{n,c})$. \qed

\begin{proposition}\label{3.14}
For any relatively free of finite
rank $n$ torsion-free nilpotent group $G$, its Lie algebra ${\cal
L}_{K}(G)$ is relatively free in some variety of nilpotent Lie
algebras.
\end{proposition}

\pf Since ${\cal L}_{\mathbb{Q}}(F_{n,c})$ is a free nilpotent Lie
algebra (by Lemma \ref{3.9}) and any fully invariant ideal of
${\cal L}_{\mathbb{Q}}(F_{n,c})$ is verbal (The arguments given in
the proof of Theorem 13.31 in \cite{neumann} are still valid for
relatively free Lie algebras.), it follows from Proposition
\ref{3.13} that ${\cal L}_{\mathbb{Q}}(F_{n,c}/{\rm ker}\pi)$ is
relatively free. Since $G \cong F_{n,c}/{\rm ker}\pi$, we obtain
${\cal L}_{\mathbb{Q}}(G)$ is relatively free. Recall that ${\cal
L}_{K}(G) = K \otimes_{\mathbb{Q}} {\cal L}_{\mathbb{Q}}(G)$.
Since relatively freeness is preserved by a field extension, we
obtain the required result (see \cite[Chapter 14]{as}). \qed

\begin{theorem}\label{3.15}
For any relatively free, torsion-free nilpotent group $G$ of
finite rank $n$, ${\cal L}_{K}(G)$ is isomorphic to $L_{K}(G)$ as
a Lie algebra via an isomorphism $\eta$ sending $\log y_{i}$ to
$y_{i} G^{\prime}$, $i = 1, \ldots, n$.
\end{theorem}

\pf Write $G = F_{n}({\mathfrak{T}}_{c})$ for some torsion-free
subvariety ${\mathfrak{T}}_{c}$ of class at most $c$ of
${\mathfrak{N}}_{c}$. It is freely generated by the set $\{y_{1},
\ldots, y_{n}\}$, where $y_{i} = f_{i}{\mathfrak{T}}_{c}(F_{n})$.
Then the factor  group $G/G^{\prime}$ is either a free abelian
group freely generated by the set $\{y_{1} G^{\prime}, \ldots,
y_{n} G^{\prime}\}$ or it is the abelian group of exponent $m > 0$
every element of whose basis $\{y_{1} G^{\prime}, \ldots, y_{n}
G^{\prime}\}$ has order $m$ (see \cite[page 11]{neumann}). If
$G/G^{\prime}$ has finite exponent $m$, then, since $G$ is
nilpotent of class $c$, we obtain $G$ has finite exponent dividing
$m^{c}$ (see, for example, \cite[page 13]{segal}) which is a
contradiction. Therefore $G/G^{\prime}$ is a free abelian group of
rank $n$. Hence $\tau_{2}(G) = G^{\prime}$.

By Lemma \ref{3.3}, $L_{K}(G)$ is generated as Lie algebra by the
set $\{y_{1} G^{\prime}, \ldots,$ $y_{n} G^{\prime}\}$. In
addition, $L_{K}(G)$ is a nilpotent Lie algebra of class $c$. Then
$L_{K}(G) \in {\mathfrak{M}}_{c}$. Since ${\cal L}_{K}(F_{n,c})$
is free in ${\mathfrak{M}}_{c}$ freely generated by the set
$\{\log x_{1}, \ldots, \log x_{n}\}$ (by Lemma \ref{3.9}), the
mapping $\rho$ from $\{\log x_{1}, \ldots, \log x_{n}\}$ into
$L_{K}(G)$ defined by $\rho(\log x_{i}) = y_{i} G^{\prime}$, $i =
1, \ldots, n$, extends uniquely to a Lie algebra epimorphism from
${\cal L}_{K}(F_{n,c})$ onto $L_{K}(G)$. We claim that $\rho({\cal
L}_{K}({\rm ker}\pi)) = \{0\}$. Recall that ${\cal L}_{K}({\rm
ker}\pi) = {\rm ker} \psi_{\pi}$ (see the proof of Proposition
\ref{3.13}), where $\psi_{\pi}$ is the Lie algebra epimorphism
from ${\cal L}_{K}(F_{n,c})$ onto ${\cal L}_{K}(G)$ such that
$\psi_{\pi}(\log x_{i}) = \log \pi(x_{i}) = \log y_{i}$, $i = 1,
\ldots, n$. Thus, by Lemma \ref{3.12},
\begin{eqnarray}
{\rm ker} \psi_{\pi} = \oplus^{c}_{t = 2} ({\rm ker} \psi_{\pi}
\cap {\cal L}_{t}(F_{n,c})).
\end{eqnarray}
To prove that $\rho({\rm ker}\psi_{\pi}) = \{0\}$ it is enough to
show that $\rho(v) = 0$ for all $v \in {\rm ker}\psi_{\pi}$. By
the equation (11), it is enough to prove that $\rho(v) = 0$ for
all homogeneous elements in ${\rm ker}\psi_{\pi}$. For any
homogeneous Lie commutator $w$ in ${\cal L}_{K}(F_{n,c})$, we
write $\widetilde{w}$ for the Lie commutator in $L_{K}(G)$
obtained from $w$ by replacing the Lie multiplication in ${\cal
L}_{K}(F_{n,c})$ by the Lie multiplication in $L_{K}(G)$, and, in
addition, by replacing each $\log x_{i}$ by $y_{i} G^{\prime}$
with $i = 1, \ldots, n$.  Let $v \in {\rm ker} \psi_{\pi} \cap
{\cal L}_{t}(F_{n,c})$ for some $t \geq 2$. Write
$$
v = \sum \beta_{(i_{1}, \ldots, i_{t})} [ \log x_{i_{1}}, \ldots,
\log x_{i_{t}}],
$$
where $\beta_{(i_{1}, \ldots, i_{t})} \in K$. Then
\begin{eqnarray}
0 = \psi_{\pi}(v) = \sum \beta_{(i_{1}, \ldots, i_{t})} [ \log
y_{i_{1}}, \ldots, \log y_{i_{t}}].
\end{eqnarray}
Since
$$
[\log y_{i_{1}}, \ldots, \log y_{i_{t}}] = \log (y_{i_{1}},
\ldots, y_{i_{t}}) + u,
$$
where $u \in {\cal L}_{K}(\tau_{t+1}(G))$, we obtain from the
equation (12)
\begin{eqnarray}
\sum \beta_{(i_{1}, \ldots, i_{t})}\log (y_{i_{1}}, \ldots,
y_{i_{t}}) \in {\cal L}(\tau_{t+1}(G)).
\end{eqnarray}
Choose a canonical basis $\{a_{1}, \ldots, a_{f(c+1)}\}$ of $G$
with $a_{i} = y_{i}$, $i = 1, \ldots, n$. Write $(y_{i_{1}},
\ldots, y_{i_{t}}) = a^{m_{f(t)+1}}_{f(t)+1} \cdots
a^{m_{f(t+1)}}_{f(t+1)} x$, where $x \in \tau_{t+1}(G)$. By the
BCH formula and since ${\cal L}_{K}(\tau_{t+1}(G))$ has a
$K$-basis the set $\{\log a_{f(t+1)+1}, \ldots,$ $\log
a_{f(c+1)}\}$, we get from (13)
$$
\sum \beta_{(i_{1}, \ldots, i_{t})} (\sum^{f(t+1)}_{k=f(t)+1}
m_{k} \log a_{k}) = 0.
$$
Let $\zeta$ be the $K$-linear isomorphism from ${\cal L}_{K}(G)$
onto $L_{K}(G)$ such that $\zeta(\log a_{f(i)+j}) = a_{f(i)+j}
\tau_{i+1}(G)$ for all $i = 1, \ldots, c$ and $j = 1, \ldots,
n_{i}$. Since ${\cal L}_{K}(\tau_{t+1}(G))$ has a $K$-basis the
set $\{\log a_{f(t+1)+1}, \ldots, \log a_{f(c+1)}\}$ and since
$\gamma_{t+1}(L_{K}(G)) = \oplus_{j=t+1}^{c} L_{j, K}(G)$, we have
$\zeta({\cal L}_{K}(\tau_{t+1}(G))) = \gamma_{t+1}(L_{K}(G))$.
Therefore
\begin{eqnarray}
\sum \beta_{(i_{1}, \ldots, i_{t})} (\prod^{f(t+1)}_{k=f(t)+1}
(a_{k} \tau_{t+1}(G))^{m_{k}}) = 0
.
\end{eqnarray}
Since
$$
[y_{i_{1}} G^{\prime}, \ldots, y_{i_{t}} G^{\prime}] = (y_{i_{1}},
\ldots, y_{i_{t}}) \tau_{t+1}(G) = \prod^{f(t+1)}_{k=f(t)+1}
(a_{k} \tau_{t+1}(G))^{m_{k}},
$$
we obtain from (14) $\sum \beta_{(i_{1}, \ldots, i_{t})}
[y_{i_{1}} G^{\prime}, \ldots,$ $y_{i_{t}} G^{\prime}] = 0$. Thus
$$
\rho(v) = \widetilde{v} = \sum \beta_{(i_{1}, \ldots, i_{t})}
[y_{i_{1}} G^{\prime}, \ldots, y_{i_{t}} G^{\prime}] = 0.
$$
Therefore $\rho({\cal L}_{K}({\rm ker}\pi)) = \{0\}$. Hence ${\cal
L}_{K}({\rm ker}\pi) \subseteq {\rm ker}\rho$. Since $\dim {\rm
ker}\rho = \dim {\cal L}_{K}(F_{n,c}) - \dim L_{K}(G) = \mu - \nu
= \dim {\cal L}_{K}({\rm ker}\pi)$, we obtain ${\cal L}_{K}({\rm
ker}\pi) = {\rm ker}\rho$. Therefore \linebreak ${\cal
L}_{K}(F_{n,c})/{\cal L}_{K}({\rm ker}\pi)\cong L_{K}(G)$ by a Lie
algebra isomorphism $\rho_{1}$ such that $\rho_{1}(\log
x_{i}+{\cal L}({\rm ker}\pi)) = y_{i} G^{\prime}$, $i = 1, \ldots,
n$. By the proof Proposition \ref{3.13}, ${\cal
L}_{K}(F_{n,c})/{\cal L}_{K}({\rm ker}\pi) \cong {\cal L}_{K}(G)$
by a Lie algebra isomorphism $\psi_{\pi,1}$ such that
$\psi_{\pi,1}(\log x_{i}+{\cal L}_K({\rm ker}\pi)) = \log y_{i}$,
$i = 1, \ldots, n$. Let $\eta = \rho_{1} \circ \psi^{-1}_{\pi,1}$.
Then $\eta$ is a Lie algebra isomorphism from ${\cal L}_{K}(G)$
into $L_{K}(G)$ such that $\eta(\log y_{i}) = y_{i} G^{\prime}$
for $i = 1, \ldots, n$. \qed

\vskip .120 in

Let $L$ be a relatively free nilpotent Lie algebra over
$\mathbb{Q}$ of finite rank $n$; freely generated by the set
$\{h_{1}, \ldots, h_{n}\}$. It is easily verified that the set
$\{h_{1}+L^{\prime}, \ldots, h_{n}+L^{\prime}\}$ is a
${\mathbb{Q}}$-basis of $L/L^{\prime}$. Give on $L$, by means of
the Baker-Campbell-Hausdorff formula, the structure of a group
denoted by $R$. Let $H$ be the subgroup of $R$ generated by the
set $\{h_{1}, \ldots, h_{n}\}$ and let $c$ be the nilpotency class
of $L$. Notice that the nilpotency class of $H$ is $c$ as well.
Since $H$ is generated by $n$ elements and it is nilpotent of
class $c$ (that is, $H \in {\mathfrak{N}}_{c}$), and since
$F_{n,c}$ is relatively free in ${\mathfrak{N}}_{c}$, the map
$\tau$ from $F_{n,c}$ into $H$ sending $x_{i}$ to $h_{i}$, $i = 1,
\ldots, n$, is a group epimorphism. By Lemma \ref{2.1},
$H^{\prime} = \tau_{2}(H)$. Thus, the first part of Lemma
\ref{3.11} shows that ${\rm ker} \tau \subseteq F^{\prime}_{n,c}$.
By Lemma \ref{3.10} (for $M = H$), $\tau$ induces a Lie algebra
homomorphism $\psi_{\tau}$ from ${\cal L}_{\mathbb{Q}}(F_{n,c})$
into ${\cal L}_{\mathbb{Q}}(H)$ such that $\psi_{\tau}(\log u) =
\log \tau(u)$ for all $u \in F_{n,c}$. Note that the set
$\{h_{1}H^{\prime}, \ldots, h_{n}H^{\prime}\}$ is a
${\mathbb{Z}}$-basis for $H/H^{\prime}$ (see \cite[Proof of
Theorem B, page 457]{bg}). By Lemma \ref{3.10}, $\psi_{\tau}$ is
surjective. Thus ${\cal L}_{\mathbb{Q}}(F_{n,c})/{\rm ker}
\psi_{\tau} \cong {\cal L}_{\mathbb{Q}}(H)$. By applying similar
arguments as in the proof of Proposition \ref{3.13}, we obtain
${\cal L}_{\mathbb{Q}}({\rm ker} \tau) = {\rm ker}\psi_{\tau}$ and
so, ${\cal L}_{\mathbb{Q}}(F_{n,c})/{\cal L}_{\mathbb{Q}}({\rm
ker}\tau) \cong {\cal L}_{\mathbb{Q}}(H)$. By Lemma \ref{2.1}, $L
\cong {\cal L}_{\mathbb{Q}}(H)$ and so, ${\cal L}_{\mathbb{Q}}(H)$
is relatively free. Thus ${\cal L}_{\mathbb{Q}}(F_{n,c})/{\cal
L}_{\mathbb{Q}}({\rm ker} \tau)$ is relatively free. It is easy to
verify that ${\cal L}_{\mathbb{Q}}({\rm ker} \tau)$ is a fully
invariant ideal of ${\cal L}_{\mathbb{Q}}(F_{n,c})$, and ${\rm
ker} \tau$ is a fully invariant subgroup of $F_{n,c}$. So, we
obtain the following result.

\begin{proposition}\label{3.16}
Let $\tau$ and $H$ be as above. \\
(i) ${\cal L}_{\mathbb{Q}}(F_{n,c})/{\cal L}_{\mathbb{Q}}({\rm ker} \tau) \cong {\cal L}_{\mathbb{Q}}(H)$. \\
(ii) ${\cal L}_{\mathbb{Q}}({\rm ker}\tau)$ \emph{is a fully
invariant ideal of}
${\cal L}_{\mathbb{Q}}(F_{n,c})$. \\
(iii) ${\rm ker} \tau$ is a fully invariant subgroup of $F_{n,c}$.
\qed
\end{proposition}

\subsection{Relatively free Lie algebras}

 An inverse
system $(G_{i}, \varphi_{ij})$ of sets indexed by a directed
nonempty set $I$ consists of a family $\{ G_{i} |~ i \in I\}$ of
sets and a family $\{ \varphi_{ij}: G_{j} \rightarrow G_{i} |~ i,
j \in I, i \leq j\}$ of maps such that $\varphi_{ii}$ is the
identity map ${\rm Id}_{G_{i}}$ for each $i$ and $\varphi_{ij}
\varphi_{jk} = \varphi_{ik}$ whenever $i \leq j \leq k$. We shall
call a family $\{ \psi_{i}: Y \rightarrow G_{i} |~ i \in I\}$ of
maps \emph{compatible} to the inverse system $(G_i,\phi_{ij})$ if $\varphi_{ij} \psi_{j} = \psi_{i}$
whenever $i \leq j$. An \emph{inverse limit} $(\widetilde{G},
\varphi_{i})$ of an inverse system $(G_{i}, \varphi_{ij})$ of sets
is a set $\widetilde{G}$ together with a compatible family $\{
\varphi_{i} : \widetilde{G} \rightarrow G_{i}\}$ of maps with the
following universal property: whenever $(\psi_{i}: Y \rightarrow
G_{i})$ is a compatible family of maps from a set $Y$, there is a
unique map $\psi: Y \rightarrow \widetilde{G}$ such that
$\varphi_{i} \psi = \psi_{i}$ for each $i$. Write $C = \prod_{i
\in I} G_{i}$ for the cartesian product of all sets $G_{i}$, and
for each $i$ write $\pi_{i}$ for the projection map from $C$ to
$G_{i}$. Define
$$
\widetilde{G} = \{c \in C: \varphi_{ij} \pi_{j}(c) = \pi_{i}(c)
~~{\rm for~~all}~~i, j ~{\rm with}~i \leq j\}
$$
and $\varphi_{i} = \pi_{i}|_{\widetilde{G}}$ for each $i$. Then
$(\widetilde{G}, \varphi_{i})$ is an inverse limit of $(G_{i},
\varphi_{ij})$, denoted by $\varprojlim G_i$. Suppose each $G_{i}$
is a Lie algebra and each $\varphi_{ij}$ is a Lie algebra
homomorphism. Then $C$ is a Lie algebra, and it is easy to see
that $\widetilde{G}$ is a Lie subalgebra of $C$. The coordinate
projections $\pi_{i}$ are obviously homomorphisms in this case.

Let $x_{i} = f_{i} {\mathfrak{L}}(F_{n})$, with $i = 1, \ldots,
n$. That is, the set $\{x_{1}, \ldots, x_{n}\}$ is a basis (i.e. a
free generating set) for $G_{n}$. Since $G_{n}$ is relatively free
and $\tau_{c+1}(G_{n})$ is a fully invariant subgroup of $G_{n}$,
we obtain $\tau_{c+1}(G_{n})$ is verbal. Write $G_{n} =
F_{n}/{\mathfrak{L}}(F_{n})$ and let $N_{c+1}$ be the complete
inverse image in $F_{n}$ of $\tau_{c+1}(G_{n})$. Then
$N_{c+1}/{\mathfrak{L}}(F_{n}) = \tau_{c+1}(G_{n})$. It is easily
verified that $N_{c+1}$ is a fully invariant subgroup of $F_{n}$
and so, $N_{c+1}$ is verbal. That is, there exists a set
$\Omega_{c+1}$ of words such that $N_{c+1} = \Omega_{c+1}(F_{n})$.
But then $N_{c+1}/{\mathfrak{L}}(F_{n}) =
\Omega_{c+1}(F_{n}/{\mathfrak{L}}(F_{n}))$ (see \cite{neumann}, the
proof of 13.31 Theorem, page 10]). Thus $G_{n}/\tau_{c+1}(G_{n})
\cong F_{n}/N_{c+1}$ and so, $G_{n}/\tau_{c+1}(G_{n})$ is a
relatively free torsion-free nilpotent group of rank $n$ and nilpotency class
at most $c$. We write $G_{n,c} = G_{n}/\tau_{c+1}(G_{n})$ and
claim that $G_{n,c}$ has class exactly $c$. Indeed, to get a
contradiction we assume that $\gamma_{d}(G_{n}) \subseteq
\tau_{c+1}(G_{n})$ for some $d \leq c$. Let $x \in
\tau_{d}(G_{n})$. Thus there exists a positive integer $m$ such
that $x^{m} \in \gamma_{d}(G_{n}) \subseteq \tau_{c+1}(G_{n})$.
Since $\tau_{c+1}(G_{n}) \subseteq \tau_{d+1}(G_{n})$ and
$\tau_{d}(G_{n})/\tau_{d+1}(G_{n})$ is torsion-free, we obtain $x
\in \tau_{d+1}(G_{n})$ and so, $\tau_{d}(G_{n}) =
\tau_{d+1}(G_{n})$ which contradicts to residually torsion-free
nilpotency. For positive integers $i$ and $c$, with $1 \leq i \leq
n$, let $x_{i,c} = x_{i} \tau_{c+1}(G_{n})$. Thus $G_{n,c}$ is
freely generated by the set $\{x_{i,c}: i = 1, \ldots, n\}$. As
shown in Proposition \ref{3.14}, and the proof of Theorem \ref{3.15}, ${\cal
L}_{K}(G_{n,c})$ is a relatively free nilpotent Lie algebra of
rank $n$; freely generated by the set $\{\log x_{1,c}, \ldots,
\log x_{n,c}\}$. Since $G_{n,c}$ has class $c$, we obtain ${\cal
L}_{K}(G_{n,c})$  has class $c$ (see \cite[Theorem 7.3]{jennings}).

The following result help us to construct the inverse limit of
${\cal L}_{K}(G_{n,c})$ with $c \geq 1$ and $n$ a fixed positive
integer with $n \geq 2$. The proof is based on some ideas of the
proof given in Lemma \ref{3.6}. (See, also, Lemma \ref{3.10}.)

\begin{lemma}\label{pap}
For positive integers $c, d$, with $c
\leq d$, we write $\pi_{c,d}$ for the natural epimorphism from
$G_{n,d}$ onto $G_{n,c}$ sending $g \tau_{d+1}(G_{n})$ to $g
\tau_{c+1}(G_{n})$ for all $g \in G_{n}$. Then there exists a Lie
algebra epimorphism $\xi_{\pi_{c,d}}$ from ${\cal L}_{K}(G_{n,d})$
onto ${\cal L}_{K}(G_{n,c})$ such that $\xi_{\pi_{c,d}}(\log (g
\tau_{d+1}(G_{n}))) = \log (g \tau_{c+1}(G_{n}))$ for all $g \in
G_{n}$. In particular, $\xi_{\pi_{c,d}}(\log x_{i,d}) = \log
x_{i,c}$ for $i = 1, \ldots, n$, and $\xi_{\pi_{c,d}}(v(\log
x_{1,d}, \ldots, \log x_{n,d})) =
 v(\log x_{1,c}, \ldots, \log x_{n,c})$ for all $v(\log x_{1,d},
\ldots, \log x_{n,d}) \in {\cal L}_{K}(G_{n,d})$. \qed
\end{lemma}

Form the following inverse limit
$$
\widetilde{{\cal L}(G_{n})} = \varprojlim  ({\cal L}_{K}(G_{n,c}),
\xi_{\pi_{c,d}}).
$$
Throughout this paper we abbreviate $\widetilde{{\cal L}}_{n} =
\widetilde{{\cal L}(G_{n})}$, with $n \geq 2$. A typical element
of $\widetilde{{\cal L}}_{n}$ has the form $(v_{1}, v_{2},
\ldots)$ where $v_{c} \in {\cal L}_{K}(G_{n,c})$, with $c \geq 1$,
and $\xi_{\pi_{c,d}}(v_{d}) = v_{c}$ for $c \leq d$. Write, for $i
= 1, \ldots, n$,
$$
t_{i} = (\log (x_{i}\tau_{2}(G_{n})), \log (x_{i}
\tau_{3}(G_{n})), \ldots, \log (x_{i}\tau_{c+1}(G_{n})), \ldots)
\in \widetilde{{\cal L}}_{n}.
$$
Let ${\cal L}_{n} = {\cal L}(G_{n})$ be the Lie subalgebra of
$\widetilde{{\cal L}}_{n}$ generated by the set $\{t_{1}, \ldots,
t_{n}\}$. Recall that for a positive integer $c$, $L_{c}(G_{n}) =
\tau_{c}(G_{n})/\tau_{c+1}(G_{n})$. Since $L_{c}(G_{n})$ is a free
abelian group of finite rank $\geq 1$, we may tensor it by $K$ to
obtain a vector space $L_{c,K}(G_{n})$ over $K$. In particular, if
$X_{c}$ is any ${\mathbb{Z}}$-basis of $L_{c}(G_{n})$ then every
element of $L_{c,K}(G_{n})$ may be uniquely written as a
$K$-linear combination of elements $1 \otimes x$, with $x \in
X_{c}$. Let $L_{K}(G_{n}) = \oplus_{c \geq 1} L_{c,K}(G_{n})$. For
a positive integer $c$, and a Lie algebra $L$, we write
$\gamma_{c}(L)$ for the $c$-th term of the lower central series of
$L$. Notice that $\gamma_{d}(L_{K}(G_{n})) = \oplus_{i \geq
d}L_{i,K}(G_{n})$. \qed

\begin{theorem}\label{3.18}
Let ${\cal L}_{n}$ be the Lie
subalgebra of $\widetilde{{\cal L}}_{n}$ generated by the set
$\{t_{1}, \ldots, t_{n}\}$. Then ${\cal L}_{n}$ is a relatively
free Lie algebra freely generated by the set $\{t_{1}, \ldots,
t_{n}\}$.
\end{theorem}

\pf For a positive integer $n$, with $n \geq 2$, let $A_{n}$ be
the free associative algebra over $K$; freely generated by the set
$\{\ell_{1}, \ldots, \ell_{n}\}$. Give on $A_{n}$ the structure of
a Lie algebra by defining $[u, v] = uv-vu$ for all $u, v \in
A_{n}$. Let $L_{n}$ be the Lie subalgebra of $A_{n}$ generated by
the set $\{\ell_{1}, \ldots, \ell_{n}\}$. It is well-known that
$L_{n}$ is a free Lie algebra; freely generated by the set
$\{\ell_{1}, \ldots, \ell_{n}\}$ (see, for example,
\cite{jacobson}). Consider the natural epimorphism
$$
\sigma_{n}: L_{n} \longrightarrow {\cal L}_{n}, \mbox{\ with\ }
\ell_i\mapsto t_i, \ i=1,\ldots,n.
$$
Then $L_{n}/{\rm ker}\sigma_{n} \cong {\cal L}_{n}$. Notice that
$v(t_{1}, \ldots, t_{n}) = 0$ if and only if $v(\ell_{1}, \ldots,
\ell_{n}) \in {\rm ker} \sigma_{n}$. To prove that ${\cal L}_{n}$
is relatively free, it is enough to show that ${\rm ker}
\sigma_{n}$ is a verbal ideal (see \cite[Chapter 14, page
275]{as}). As for groups \cite[Theorem 12.34, page 5]{neumann} verbal
ideals turn out to be precisely the fully invariant ideals of
$L_{n}$ (that is, those ideals which are invariant under all
endomorphisms of $L_{n}$). Thus it is enough to show that ${\rm
ker} \sigma_{n}$ is fully invariant i.e. if $v(\ell_{1}, \ldots,
\ell_{n}) \in {\rm ker} \sigma_{n}$, then $v(v_{1}, \ldots, v_{n})
\in {\rm ker} \sigma_{n}$ for all $v_{1}, \ldots, v_{n} \in L_{n}$
or, equivalently, if $v(t_{1}, \ldots, t_{n}) = 0$, then $v(u_{1},
\ldots, u_{n}) = 0$ for all $u_{1}, \ldots, u_{n} \in {\cal
L}_{n}$. Each $t_{i}$ may be considered as a function on
${\mathbb{N}}$ such that $t_{i}(c) = \log x_{i,c} \in {\cal
L}_{K}(G_{n,c})$ for all $c \in {\mathbb{N}}$. The element
$v(t_{1}, \ldots, t_{n})$ is considered as a function on
${\mathbb{N}}$. Thus $v(t_{1}(c), \ldots, t_{n}(c)) = 0$ for all
$c \in {\mathbb{N}}$. Hence $v(t_{1}(c), \ldots, t_{n}(c)) = 0$ in
each ${\cal L}_{K}(G_{n,c})$. For any positive integer $c$, ${\cal
L}_{K}(G_{n,c})$ is a relatively free nilpotent Lie algebra;
freely generated by the set $\{\log x_{1,c}, \ldots, \log
x_{n,c}\}$. Since $t_{i}(c) = \log x_{i,c}$ for $i = 1, \ldots,
n$, we obtain $\{t_{1}(c), \ldots, t_{n}(c)\}$ is a free
generating set for ${\cal L}_{K}(G_{n,c})$. Therefore $v(w_{1,c},
\ldots, w_{n,c}) = 0$ for all $w_{1,c}, \ldots, w_{n,c} \in$
${\cal L}(G_{n,c})$. By the construction of $\widetilde{{\cal
L}}_{n}$, we have $\widetilde{{\cal L}}_{n} \in {\rm var}(\{{\cal
L}_{K}(G_{n,c}): c \in {\mathbb{N}}\})$. But $v(\ell_{1}, \ldots,
\ell_{n})$ is an identity in the cartesian product $\prod_{c \geq
1} {\cal L}_{K}(G_{n,c})$. (The arguments given in the proof of
15.1 of \cite[page 15]{neumann} can be adopted here without any
changes.) Hence $v(\ell_{1}, \ldots, \ell_{n})$ is an identity in
$\widetilde{{\cal L}}_{n}$. Since ${\cal L}_{n} \subseteq
\widetilde{{\cal L}}_{n}$, we obtain $v(\ell_{1}, \ldots,
\ell_{n})$ is an identity in ${\cal L}_{n}$. Therefore, ${\cal
L}_{n}$ is a relatively free Lie algebra; freely generated by the
set $\{t_{1}, \ldots, t_{n}\}$. \qed

\subsection{Verbal ideals}

 Let ${\cal F}_{\infty}$ be the
free Lie algebra on a countably infinite set $\{\omega_{1},
\omega_{2}, \ldots, \}$. For each $n \geq 1$, the free Lie algebra
${\cal F}_{n}$ freely generated by $\omega_{1}, \ldots,
\omega_{n}$ will then be embedded in ${\cal F}_{\infty}$, in a
natural way. This free Lie algebra ${\cal F}_{\infty}$ is
introduced for the special purpose to provide \lq words\rq: an
element $w \in {\cal F}_{\infty}$ is called a \emph{word} in the
variables $\omega_{1}, \omega_{2}, \ldots$. Each word $w$ involves
only finitely many variables, in the sense that $w \in {\cal
F}_{n}$ for some $n \in {\mathbb{N}}$. A set of words $V$ is
\emph{closed} if it is, as subset of ${\cal F}_{\infty}$, a fully
invariant ideal of ${\cal F}_{\infty}$. In particular, if $v \in
V$ is a word in $n$ variables, and $(u_{1}, \ldots, u_{n}) \in
{\cal F}_{\infty}^{n}$ an $n$-tuple of words, then $v(u_{1},
\ldots, u_{n}) \in V$. The \emph{closure} of an arbitrary set
${\mathfrak{w}}$ of words is defined as the intersection of all
closed sets containing ${\mathfrak{w}}$. (Since the set of all
words is a closed set containing ${\mathfrak{w}}$, the
aforementioned definition makes sense.) Let $Y$ be a set of words
involving $\omega_{1}, \ldots, \omega_{n}$ only such that $Y$ is a
fully invariant ideal of ${\cal F}_{n}$ (hence it is a verbal
ideal). If $V$ denotes the closure of $Y$ then $V \cap {\cal
F}_{n} = Y$ (see, for example, \cite[page 7]{neumann}. The
arguments given in page 7 can be adopted here without any
changes]). That is, the closure of a fully invariant ideal of
${\cal F}_{n}$ intersected with ${\cal F}_{n}$ leads back to the
original ideal. Let $L_{\infty}$ be a free Lie algebra on a
countably infinite set $\{\ell_{1}, \ell_{2}, \ldots\}$. For a
positive integer $n$, with $n \geq 2$, let $L_{n}$ denote the free
Lie algebra freely generated by the set $\{\ell_{1}, \ldots,
\ell_{n}\}$. Let $\xi$ be the natural mapping of ${\cal
F}_{\infty}$ onto $L_{\infty}$; it is given by $\xi(\omega_{i}) =
\ell_{i}$ for $i \geq 1$. It is an isomorphism between the two
free Lie algebras which, restricted to ${\cal F}_{n} \subset {\cal
F}_{\infty}$ maps ${\cal F}_{n}$ onto $L_{n} \subset L_{\infty}$.
So we have the following diagram:
$$\begin{diagram}
\node{{\cal F}_{\infty}}
\arrow{e,t}{} \arrow{s,t}{\xi} \node{{\cal F}_{n}}\arrow{s,r,..}{}\\
\node{L_{\infty}} \arrow{e,t}{} \node{L_n} \arrow{e,t}{\sigma_n}
\node{{\cal L}_{n}}
\end{diagram}$$
Let $V$ be a fully invariant ideal of ${\cal F}_{\infty}$. The
words in $V$ that involve only $\omega_{1}, \ldots, \omega_{n}$
are given by $V \cap {\cal F}_{n}$. Of course, $V \cap {\cal
F}_{n}$ is a fully invariant ideal of ${\cal F}_{n}$. For a
positive integer $n$, with $n \geq 2$, let $V_{n}$ be the fully
invariant ideal of ${\cal F}_{n}$ such that $\xi(V_{n}) = {\rm
ker} \sigma_{n}$. Let $Y_{n}$ be the closure of $V_{n}$. Then
$Y_{n} \cap {\cal F}_{n} = V_{n}$. Since $\xi$ is bijective, we
obtain $\xi(Y_{n}) \cap L_{n} = \xi(V_{n}) = {\rm ker}
\sigma_{n}$. Recall that there is a $1-1$ correspondence between
the varieties ${\mathfrak{V}}$ of Lie algebras and the fully
invariant ideals $V$ of ${\cal F}_{\infty}$ (or the closed sets of
words $V$ or the verbal ideals $V$). Let ${\mathfrak{V}}_{n}$ (or
simply ${\mathfrak{V}}$) be the variety of Lie algebras
corresponding to $Y_{n}$. Write ${\mathfrak{V}}(L_{n}) =
\xi(Y_{n}) \cap L_{n}$. That is, the verbal ideal of $L_{n}$
corresponding to ${\mathfrak{V}}$. Notice that
${\mathfrak{V}}(L_{n}) = {\rm ker} \sigma_{n}$. It is well known
that $L_{n} = \oplus_{m \geq 1} L^{m}_{n}$, where $L^{m}_{n}$
denotes the vector subspace of $L_{n}$ spanned by all Lie
commutators $[\ell_{i_{1}}, \ldots, \ell_{i_{m}}]$ with $i_{1},
\ldots, i_{m} \in \{1, \ldots, n\}$. Since $K$ is an infinite
field, we obtain ${\mathfrak{V}}$ is (multi-)homogeneous. That is,
given an identity $v \equiv 0$ in ${\mathfrak{V}}$ such that $v =
\sum_{\alpha} v_{\alpha}$ is the (multi-)homogeneous decomposition
of $v$, then for all $\alpha$, $v_{\alpha} \equiv 0$ is also an
identity in ${\mathfrak{V}}$. The proof of the following result is
elementary.

\begin{lemma}\label{3.19}
For a positive integer $n$, with $n \geq 2$, let $L_{n}$ be a free
Lie algebra of rank $n$. Let $I$ be a proper fully invariant ideal
of $L_{n}$. Then $I \subseteq L^{\prime}_{n}$. \qed
\end{lemma}

Since ${\rm ker}\sigma_{n} \subseteq L^{\prime}_{n}$ (by the proof
of Theorem \ref{3.18} and Lemma \ref{3.19}),
$$
{\rm ker}\sigma_{n} = {\mathfrak{V}}(L_{n}) = \oplus_{m \geq 2}
({\rm ker}\sigma_{n} \cap L^{m}_{n}).
$$
Since ${\mathfrak{V}}(L_{n}) = \oplus_{m \geq
2}({\mathfrak{V}}(L_{n}) \cap L^{m}_{n})$, it is easily verified
that ${\cal L}_{n} = \oplus_{m \geq 1} {\cal L}^{m}_{n}$, where
${\cal L}^{m}_{n} = (L^{m}_{n} +
{\mathfrak{V}}(L_{n}))/{\mathfrak{V}}(L_{n})$. Furthermore,
$\gamma_{c}({\cal L}_{n}) = \oplus_{m \geq c}{\cal L}^{m}_{n}$ for
all $c$. For positive integers $i$ and $c$, with $i \leq c$, let
$n_{i}$ denote the rank of the free abelian group $L_{i}(G_{n,c})$
$(= \gamma_{i}(G_{n,c})/\gamma_{i+1}(G_{n,c}))$. Furthermore, let
$f(i) = n_{1} + \cdots + n_{i-1}$, with $n_{0} = 0$ and $n_{1} =
n$. Let $B_{i,c} = \{a_{f(i)+1,c}, \ldots, a_{f(i+1,c)}\}$ be a
subset of $\gamma_{i}(G_{n,c})$ such that the set
$\{a_{f(i)+1,c}\gamma_{i+1}(G_{n,c}), \ldots,
a_{f(i+1),c}\gamma_{i+1}(G_{n,c})\}$ is a ${\mathbb Z}$-basis of
$L_{i}(G_{n,c})$. Thus $B_{c} = \cup B_{i,c} = \{a_{1,c}, \ldots,
a_{f(c+1),c}\}$ is a canonical basis of $G_{n,c}$. For any
positive integer $c$, we choose a canonical basis $B_{c}$ of
$G_{n,c}$ subject to $a_{\kappa,c} = x_{\kappa,c}$, $\kappa = 1,
\ldots, n$. For a positive integer $c$, with $c \geq 2$, let
$\pi_{c-1,c}$ be the natural epimorphism from $G_{n,c}$ onto
$G_{n,c-1}$ sending $g \tau_{c+1}(G_{n})$ to $g \tau_{c}(G_{n})$
for all $g \in G_{n}$. Then ${\rm ker} \pi_{c-1,c} =
\tau_{c}(G_{n})/\tau_{c+1}(G_{n})$. By Lemma \ref{pap}, there
exists a Lie algebra epimorphism $\xi_{\pi_{c-1,c}}$ from ${\cal
L}_{K}(G_{n,c})$ onto ${\cal L}_{K}(G_{n,c-1})$ such that
$\xi_{\pi_{c-1,c}}(\log(g \tau_{c+1}(G_{n}))) = \log(g
\tau_{c}(G_{n}))$ for all $g \in G_{n}$. In particular,
$\xi_{\pi_{c-1,c}}(\log x_{i,c}) = \log x_{i,c-1}$ for $i = 1,
\ldots, n$.

\begin{proposition}\label{3.20}
For a positive integer $n$, with
$n \geq 2$, let $H_{n}$ be a relatively free group of rank $n$.
Then $\tau_{c}(H_{n}) = \gamma_{c}(H_{n}) \tau_{c+1}(H_{n})$ for
all $c$.
\end{proposition}

\pf Let $H_{n}$ be a relatively free group of rank $n$,
with $n \geq 2$ freely generated by the set ${\mathfrak{h}}$. The
commutator group $H_{n}/H^{\prime}_{n}$ is an abelian relatively
free group freely generated by ${\mathfrak{h}}$ taken modulo
$H^{\prime}_{n}$. Suppose first that $H_{n}/H^{\prime}_{n}$ is a
free abelian group of exponent $m$, with $m > 0$. Then
$\tau(H_{n}/H^{\prime}_{n}) = H_{n}/H^{\prime}_{n}$. Since
$\tau_{2}(H_{n})/\gamma_{2}(H_{n}) = \tau(H_{n}/H^{\prime}_{n})$,
we obtain $\tau_{2}(H_{n}) = H_{n}$. We claim that $H_{n} =
\tau_{c}(H_{n})$ for all $c$. Since $H_{n} = \tau_{2}(H_{n})$, we
have $x \in \tau_{2}(H_{n})$ for all $x \in H_{n}$. Since
$H_{n}/H^{\prime}_{n}$ has exponent $m$, we get $x^{m} \in
H^{\prime}_{n}$. Then $(x^{m}, y) \in \gamma_{3}(H_{n})$ for all
$y \in H_{n}$. Using repeatedly the commutator identity $(ab, c) =
(a, c) (a, c, b)(b, c)$, we obtain $(x, y)^{m} \in
\gamma_{3}(H_{n})$ for all $y \in H_{n}$. So, $(x, y) \in
\tau_{3}(H_{n})$ for all $y \in H_{n}$. It is clearly enough that
$(x, y) \in \tau_{3}(H_{n})$ for all $x, y \in H_{n}$. Since
$\tau_{3}(H_{n})$ is fully invariant, we have $\gamma_{2}(H_{n})
\subseteq \tau_{3}(H_{n})$. Let $x \in \tau_{2}(H_{n})$. Then
$x^{m} \in \gamma_{2}(H_{n})$ and so $x^{m} \in \tau_{3}(H_{n})$.
Since $\tau_{2}(H_{n})/\tau_{3}(H_{n})$ is torsion-free, we get $x
\in \tau_{3}(H_{n})$. Therefore $\tau_{2}(H_{n}) =
\tau_{3}(H_{n})$. Continuing this process we obtain $H_{n} =
\tau_{c}(H_{n})$ for all $c$. So, we get the required result.

Thus we may assume that $H_{n}/H^{\prime}_{n}$ is a free abelian
group with a basis (i.e. a free generating set) ${\mathfrak{h}}$
taken modulo $H^{\prime}_{n}$. Assume that there are no
repetitions of terms of the series $\{\tau_{c}(H_{n})\}_{c \geq
1}$. Fix a positive integer $c$. It is clearly enough that we may
assume that $c \geq 2$. It is enough to show that $\tau_{c}(H_{n})
\subseteq \gamma_{c}(H_{n}) \tau_{c+1}(H_{n})$. Write $H_{n,c} =
H_{n}/\tau_{c+1}(H_{n})$. Notice that $H_{n,c}$ is a relatively
free nilpotent torsion-free group of rank $n$ and class $c$. (If
the class of $H_{n,c}$ is strictly less than $c$, then it be can
easily shown that there are repetitions of the series
$\{\tau_{c}(H_{n})\}_{c \geq 1}$.) By Theorem A (I), $H_{n,c}$ is
Magnus and so, $\tau_{\kappa}(H_{n,c}) = \gamma_{\kappa}(H_{n,c})$
for all $\kappa$, $\kappa \in \{1, \ldots, c\}$. Write $E_{n,c} =
H_{n}/\gamma_{c}(H_{n}) \tau_{c+1}(H_{n})$. Since
$\gamma_{c}(H_{n,c}) = \gamma_{c}(H_{n})
\tau_{c+1}(H_{n})/\tau_{c+1}(H_{n})$, we have
$H_{n,c}/\gamma_{c}(H_{n,c}) \cong E_{n,c}$. Thus $E_{n,c}$ is
torsion-free. Let $\rho_{c-1,c}$ denote the natural epimorphism
from $E_{n,c}$ onto $H_{n,c-1}$ sending $h (\gamma_{c}(H_{n})
\tau_{c+1}(H_{n}))$ to $h \tau_{c}(H_{n})$ for all $h \in H_{n}$.
It is easily verified that ${\rm ker} \rho_{c-1,c} =
\tau_{c}(H_{n})/\gamma_{c}(H_{n}) \tau_{c+1}(H_{n})$. Furthermore,
it is easy to see that ${\rm ker} \rho_{c-1,c}$ is an homomorphic
image of $\tau_{c}(H_{n})/\gamma_{c}(H_{n})$. Since
$H_{n}/\gamma_{c}(H_{n})$ is a finitely generated nilpotent group,
we have $\tau(H_{n}/\gamma_{c}(H_{n}))$ is a finite group and so,
$\tau_{c}(H_{n})/\gamma_{c}(H_{n})$ is a finite group. Therefore
${\rm ker} \rho_{c-1,c}$ is finite. Since $E_{n,c}$ is
torsion-free, we obtain ${\rm ker} \rho_{c-1,c}$ is trivial and so
$\tau_{c}(H_{n}) \subseteq \gamma_{c}(H_{n}) \tau_{c+1}(H_{n})$.
Hence $\tau_{c}(H_{n}) = \gamma_{c}(H_{n}) \tau_{c+1}(H_{n})$ for
all $c$. Finally we assume that there are repetitions of terms of
the series $\{\tau_{c}(H_{n})\}_{c \geq 1}$. Since
$H_{n}/H^{\prime}_{n}$ is free abelian, we have $\tau_{2}(H_{n}) =
H^{\prime}_{n}$. Let $\kappa$ be the smallest positive integer
such that $\tau_{\kappa}(H_{n}) = \tau_{\kappa + 1}(H_{n})$ and
$H_{n} \supset \tau_{2}(H_{n}) \supset \cdots \supset \tau_{\kappa
- 1}(H_{n}) \supset \tau_{\kappa}(H_{n})$. Since
$\tau_{\kappa}(H_{n}) = \tau_{\kappa + 1}(H_{n})$, we get
$\tau_{j}(H_{n}) = \tau_{\kappa}(H_{n})$ for all $j \geq \kappa +
1$. Hence our claim holds for all $j \geq \kappa$. If $\kappa = 2$
then $\tau_{c}(H_{n}) = \gamma_{c}(H_{n}) \tau_{c+1}(H_{n})$ for
all $c$. Thus we concentrate on $j < \kappa$ and $\kappa \geq 3$.
Write $H_{n,c} = H_{n}/\tau_{c+1}(H_{n})$ for $c = 1, \ldots,
\kappa - 1$, and fix $c$. Using similar arguments as before we
obtain the required result. Hence, in any case, we have
$\tau_{c}(H_{n}) = \gamma_{c}(H_{n}) \tau_{c+1}(H_{n})$. \qed

\begin{remark}\label{3.21}
{\upshape The proof given in the first part of the
proof of Proposition \ref{3.20} is independent up to the rank of
$H_{n}$.}
\end{remark}

\noindent We deduce the following result.

\begin{corollary}\label{3.22}
Let ${\mathfrak{L}}$ be a residually
torsion-free nilpotent variety of groups. For positive integers
$n$ and $c$, with $n \geq 2$, let $G_{n} = F_{n}({\mathfrak{L}})$
be the relatively free group of rank $n$ in ${\mathfrak{L}}$ and
$G_{n,c} = G_{n}/\tau_{c+1}(G_{n})$. Then, for all $c$,
$\gamma_{c}(G_{n,c}) = \tau_{c}(G_{n})/\tau_{c+1}(G_{n})$.
\end{corollary}

\pf Since $\gamma_{c}(G_{n,c}) = \gamma_{c}(G_{n})
\tau_{c+1}(G_{n})/\tau_{c+1}(G_{n})$, we obtain from Proposition
\ref{3.20} the required result. \qed

For a positive integer $c$, let $\psi_{c}$ be the mapping from
$\{\ell_{1}, \ldots, \ell_{n}\}$ into ${\cal L}_{K}(G_{n,c})$ such
that $\psi_{c}(\ell_{i}) = t_{i}(c) = \log x_{i,c}$, $i = 1,
\ldots, n$. Since $L_{n}$ is free on $\{\ell_{1}, \ldots,
\ell_{n}\}$, $\psi_{c}$ extends to a Lie algebra homomorphism from
$L_{n}$ into ${\cal L}_{K}(G_{n,c})$. Since ${\cal
L}_{K}(G_{n,c})$ is generated by the set $\{\log x_{1,c}, \ldots,
\log x_{n,c}\}$, we have $\psi_{c}$ is surjective. Let
$v(\ell_{1}, \ldots, \ell_{n}) \in {\rm ker} \sigma_{n}$. Then
$v(t_{1}, \ldots, t_{n}) = 0$ in ${\cal L}_{n}$ and so, $v(\log
x_{1,c}, \ldots, \log x_{n,c}) = 0$ for all $c \in {\mathbb{N}}$
(see the proof of Theorem \ref{3.18}). Hence $\psi_{c}(v(\ell_{1},
\ldots, \ell_{n})) = 0$ i.e. $v(\ell_{1}, \ldots, \ell_{n}) \in
{\rm ker} \psi_{c}$. Therefore ${\rm ker} \sigma_{n} \subseteq
{\rm ker} \psi_{c}$ and so, $\psi_{c}$ induces a Lie algebra
epimorphism $\widetilde{\psi}_{c}$, say, from ${\cal L}_{n}$ onto
${\cal L}_{K}(G_{n,c})$ (for all $c$) sending $t_{i}$ to $\log
x_{i,c}$ for $i = 1, \ldots, n$. So
$$\begin{diagram}
\node{L_n} \arrow{e,t}{\sigma_n} \arrow{se,b}{\psi_c} \node{{\cal L}_n} \arrow{s,r,..}{\widetilde{\psi}_c}\\
\node{} \node{{\cal L}_{K}(G_{n,c})}
  \end{diagram}$$

\begin{lemma}\label{3.23}
For a positive integer $c$,
$\widetilde{\psi}_{c}({\cal L}_{n}^{c}) = {\cal
L}_{K}(\gamma_{c}(G_{n,c}))$. In particular, for $c \geq 2$,
$\widetilde{\psi}_{c}(\gamma_{c}({\cal L}_{n})) = {\rm ker}
\xi_{\pi_{c-1,c}}$. Furthermore, ${\rm ker} \sigma_{n} = \cap_{c
\geq 1} {\rm ker} \psi_{c}$.
\end{lemma}

\pf Recall that, for a positive integer $c$, ${\cal
L}^{c}_{n}$ is the vector $K$-space spanned by all Lie commutators
$[t_{i_{1}}, \ldots, t_{i_{c}}]$ with $i_{1}, \ldots, i_{c} \in
\{1, \ldots, n\}$. Observe that
\begin{center}
$[t_{i_{1}}, \ldots, t_{i_{c}}] = (\underbrace{0,
\ldots,0}_{c-1},[t_{i_{1}}(c), \ldots, t_{i_{c}}(c)],
[t_{i_{1}}(c+1), \ldots, t_{i_{c}}(c+1)], \ldots) \in
\widetilde{{\cal L}}_{n}.$
\end{center}
Applying $\widetilde{\psi}_{c}$ on $[t_{i_{1}}, \ldots,
t_{i_{c}}]$, we obtain $\widetilde{\psi}_{c}([t_{i_{1}}, \ldots,
t_{i_{c}}]) = [t_{i_{1}}(c), \ldots, t_{i_{c}}(c)]$. Thus
$\widetilde{\psi}_{c}({\cal L}_{n}^{c}) \subseteq \gamma_{c}({\cal
L}(G_{n,c}))$. By Lemma \ref{3.7} and since $G_{n,c}$ is Magnus,
we have $\widetilde{\psi}_{c}({\cal L}_{n}^{c}) \subseteq {\cal
L}(\gamma_{c}(G_{n,c}))$. Let $c = 1$. Since ${\rm ker} \sigma_{n}
\subseteq L^{\prime}_{n}$, we obtain ${\cal L}_{n}/{\cal
L}^{\prime}_{n} \cong L_{n}/L^{\prime}_{n}$ as vector spaces.
Since ${\cal L}_{n}^{\prime} = \oplus_{c \geq 2} {\cal
L}_{n}^{c}$, we get $\dim {\cal L}_{n}^{1} = n$ and so $\{t_{1},
\ldots, t_{n}\}$ is a $K$-basis for ${\cal L}_{n}^{1}$. Since
$\dim {\cal L}(G_{n,1}) = {\cal H}(G_{n,1}) = n$, we obtain
$\{t_{1}(c), \ldots, t_{n}(c)\}$ is a $K$-basis for ${\cal
L}_{K}(G_{n,1})$. It is clearly enough that
$\widetilde{\psi}_{1}({\cal L}_{n}^{1}) = {\cal L}_{K}(G_{n,1})$.
Thus we may assume that $c \geq 2$. By Lemma \ref{3.6}, ${\cal
L}_{K}(\gamma_{c}(G_{n,c})) = {\rm ker} \xi_{\pi_{c-1,c}}$ and so,
$\widetilde{\psi}_{c}({\cal L}_{n}^{c}) \subseteq {\rm ker}
\xi_{\pi_{c-1,c}}$. To prove that $\gamma_{c}({\cal L}(G_{n,c}))
\subseteq \widetilde{\psi}_{c}({\cal L}_{n}^{c})$ it is enough to
show that $[t_{i_{1}}(c), \ldots, t_{i_{c}}(c)] \in
\widetilde{\psi}_{c}({\cal L}_{n}^{c})$ for all $i_{1}, \ldots,
i_{c} \in \{1, \ldots, n\}$. Since $\gamma_{c}({\cal
L}_{K}(G_{n,c})) = {\cal L}_{K}(\gamma_{c}(G_{n,c})) = {\rm ker}
\xi_{\pi_{c-1,c}}$, we obtain $[t_{i_{1}}(c-1), \ldots,
t_{i_{c}}(c-1)] = 0$. Consider the element
\begin{center}
$v = (\underbrace{0, \ldots,0}_{c-1},[t_{i_{1}}(c), \ldots,
t_{i_{c}}(c)], [t_{i_{1}}(c+1), \ldots, t_{i_{c}}(c+1)], \ldots)
\in \widetilde{{\cal L}}_{n}.$
\end{center}
It is easily seen that $[t_{i_{1}}, \ldots, t_{i_{c}}] = v$, and
$\widetilde{\psi}_{c}(v) = [t_{i_{1}}(c), \ldots, t_{i_{c}}(c)]$.
Therefore $\widetilde{\psi}_{c}({\cal L}_{n}^{c}) =
\gamma_{c}({\cal L}_{K}(G_{n,c})) = {\rm ker} \xi_{\pi_{c-1,c}}$.
Since $\gamma_{c}({\cal L}_{n}) = \oplus_{m \geq c} {\cal
L}_{n}^{m}$, it is easily seen that
$\widetilde{\psi}_{c}(\gamma_{c}({\cal L}_{n})) = {\rm ker}
\xi_{\pi_{c-1,c}}$ for all $c \geq 2$.

Since ${\rm ker}\sigma_{n} \subseteq {\rm ker} \psi_{c}$ for all
$c$, we have ${\rm ker} \sigma_{n} \subseteq \cap_{c \geq 1} {\rm
ker} \psi_{c}$. Let $u = u(\ell_{1}, \ldots, \ell_{n}) \in \cap_{c
\geq 1} {\rm ker} \psi_{c}$. Since $u \in {\rm ker} \psi_{c}$ for
all $c$, we have $\psi_{c}(u(\ell_{1}, \ldots, \ell_{n}))$ $=
u(t_{1}(c), \ldots, t_{n}(c)) = 0$ for all $c$. Since ${\cal
L}_{K}(G_{n,c})$ is relatively free on $\{t_{1}(c), \ldots,
t_{n}(c)\}$, we have $u(\ell_{1}, \ldots, \ell_{n})$ is an
identity for ${\cal L}_{K}(G_{n,c})$ for all $c$. This means that
$u(\ell_{1}, \ldots, \ell_{n})$ is an identity of the cartesian
product $\prod_{c \geq 1} {\cal L}_{K}(G_{n,c})$. Thus,
$u(\ell_{1}, \ldots, \ell_{n})$ is an identity of
$\widetilde{{\cal L}}_{n}$ and so, it is an identity of ${\cal
L}_{n}$. That is, $u(\ell_{1}, \ldots,\ell_{n}) \in
{\mathfrak{V}}(L_{n}) = {\rm ker} \sigma_{n}$. Therefore ${\rm
ker} \sigma_{n} = \cap_{c \geq 1} {\rm ker} \psi_{c}$. \qed

In the next few lines, we write $I_{c} = {\rm ker}
\widetilde{\psi}_{c}$. We claim that $I_{c}$ is a fully invariant
ideal of ${\cal L}_{n}$. It is enough to show that
$\vartheta(I_{c}) \subseteq I_{c}$ for any $\vartheta$
endomorphism of ${\cal L}_{n}$. Let $u = u(t_{1}, \ldots, t_{n})
\in I_{c}$. Then $\widetilde{\psi}_{c}(u(t_{1}, \ldots, t_{n})) =
u(\log x_{1,c}, \ldots, \log x_{n,c}) = 0$ in ${\cal
L}_{K}(G_{n,c})$. Since ${\cal L}_{K}(G_{n,c})$ is relatively free
on $\{\log x_{1,c}, \ldots, \log x_{n,c}\}$, we obtain $u(u_{1},
\ldots, u_{n}) = 0$ for all $u_{1}, \ldots, u_{n} \in {\cal
L}_{K}(G_{n,c})$. Let $\vartheta$ be an endomorphism of ${\cal
L}_{n}$. Then $\vartheta(t_{i}) = v_{i}$ for $i = 1, \ldots, n$.
Since $\widetilde{\psi}_{c}$ is a Lie algebra homomorphism from
${\cal L}_{n}$ into ${\cal L}(G_{n,c})$,
$\widetilde{\psi}_{c}(v_{i}) = w_{i}$ for $i = 1, \ldots, n$. But
$u(w_{1}, \ldots, w_{n}) = 0$ (in ${\cal L}_{K}(G_{n,c})$) and so,
$\widetilde{\psi}_{c}(\vartheta(u(t_{1}, \ldots, t_{n}))) =
\widetilde{\psi}_{c}(u(v_{1}, \ldots, v_{n})) = u(w_{1}, \ldots,
w_{n}) = 0$. Therefore $\vartheta(u(t_{1}, \ldots, t_{n})) \in
I_{c}$. Hence $I_{c}$ is fully invariant (and so, $I_{c}$ is
verbal). Since ${\cal L}_{n}^{\prime} = \oplus_{m \geq 2} {\cal
L}_{n}^{m}$ and $I_{c}$ is a proper fully invariant ideal of
${\cal L}_{n}$, it is easy to see that $I_{c} \subseteq {\cal
L}^{\prime}_{n}$.

Next we claim that $I_{c} = \gamma_{c+1}({\cal L}_{n})$. Since
${\cal L}_{n}/I_{c} \cong {\cal L}_{K}(G_{n,c})$ and ${\cal
L}_{K}(G_{n,c})$ has class $c$, it is enough to show that $I_{c}
\subseteq \gamma_{c+1}({\cal L}_{n})$. To get a contradiction, let
$w \in I_{c}$ but not in $\gamma_{c+1}({\cal L}_{n})$. Since
$I_{c}$ is a proper fully invariant ideal and since $K$ is an
infinite field, we may assume that $w \in I_{c} \cap {\cal
L}_{n}^{d}$ for some $d$, with $2 \leq d \leq c$. We write
$$
w = w(t_{1}, \ldots, t_{n}) = \sum \alpha_{(i_{1}, \ldots,
i_{d})}[t_{i_{1}}, \ldots, t_{i_{d}}],
$$
$\alpha_{(i_{1}, \ldots, i_{d})} \in K$. Since
$\widetilde{\psi}_{c}(w) = 0$, we obtain
$$
w(t_{1}(c), \ldots, t_{n}(c)) = \sum \alpha_{(i_{1}, \ldots,
i_{d})}[t_{i_{1}}(c), \ldots, t_{i_{d}}(c)] = 0
$$
in ${\cal L}_{K}(G_{n,c})$. Making use of $\xi_{\pi_{\mu,c}}$
(Lemma \ref{pap}), with $\mu \leq c$, we have $w(t_{1}(\mu),
\ldots, t_{n}(\mu)) = 0$ for all $\mu \leq c$. Notice that
$$
\widetilde{\psi}_{\kappa}(w) = w(t_{1}(\kappa), \ldots,
t_{n}(\kappa)) \in \gamma_{d}({\cal L}(G_{n,\kappa}))
$$
for all $\kappa \geq c$. Since $\xi_{\pi_{c,c+1}}(w(t_{1}(c+1)),
\ldots, t_{n}(c+1)) = w(t_{1}(c), \ldots, t_{n}(c)) = 0$, we
obtain $w(t_{1}(c+1), \ldots, t_{n}(c+1)) \in {\rm ker}
\xi_{\pi_{c,c+1}}$. By Lemma \ref{3.6} and (Lemma \ref{3.7}),
$w(t_{1}(c+1), \ldots,$ \linebreak  $t_{n}(c+1)) \in
\gamma_{c+1}({\cal L}(G_{n,c+1}))$. For a positive integer $m$,
with $m = 1, \ldots, c$, let ${\cal L}_{m}(G_{n,c})$ denote the
vector subspace of ${\cal L}_{K}(G_{n,c})$ spanned by all Lie
commutators of the form \linebreak $[\log x_{i_{1},c}, \ldots,
\log x_{i_{m},c}]$. By Theorem \ref{3.15} and since
$L_{K}(G_{n,c})$ is graded,
$$
{\cal L}_{K}(G_{n,c}) = \oplus_{m=1}^{c} {\cal L}_{m}(G_{n,c}).
$$
Since $w(t_{1}(c+1), \ldots, t_{n}(c+1))$ has length $d$ and $d
\leq c$, we obtain $w(t_{1}(c+1), \ldots,$ \linebreak $t_{n}(c+1))
= 0$ (in ${\cal L}_{K}(G_{n,c+1})$). Thus
$\widetilde{\psi}_{c+1}(w) = 0$. Continue this process, we finally
obtain $w = 0$ which is a contradiction. Therefore $I_{c} =
\gamma_{c+1}({\cal L}_{n})$. Furthermore, ${\cal L}_{n}^{c} \cong
\gamma_{c}({\cal L}_{n})/\gamma_{c+1}({\cal L}_{n}) \cong
\gamma_{c}({\cal L}_{K}(G_{n,c}) = {\cal
L}_{K}(\gamma_{c}(G_{n,c}))$ (by Lemma \ref{3.7}). Thus we obtain
the following result.

\begin{proposition}\label{3.24}
For a positive integer $c$, ${\cal L}_{n}/\gamma_{c+1}({\cal
L}_{n}) \cong {\cal L}_{K}(G_{n,c})$ as Lie algebras under the
isomorphism $\varphi_{c+1}$ sending $t_{i} + \gamma_{c+1}({\cal
L}_{n})$ to $\log x_{i,c}$, $i = 1, \ldots, n$. In particular,
$\varphi_{c+1}(u(t_{1}, \ldots, t_{n}) + \gamma_{c+1}({\cal
L}_{n})) = u(\log x_{1,c}, \ldots, \log x_{n,c})$ for all
$u(t_{1}, \ldots, t_{n}) \in {\cal L}_{n}$ and $c \geq 1$.
Furthermore, for all positive integers $c$, ${\cal L}_{n}^{c}
\cong {\cal L}_{K}(\gamma_{c}(G_{n,c}))$ as vector spaces via the
linear isomorphism $\varphi_{c+1} \eta_{c}$, where $\eta_{c}$ is
the natural linear isomorphism from ${\cal L}^{c}_{n}$ onto
$\gamma_{c}({\cal L}_{n})/\gamma_{c+1}({\cal L}_{n})$ sending $u$
to $u + \gamma_{c+1}({\cal L}_{n})$. \qed
\end{proposition}

Since $G_{n} = F_{n}({\mathfrak{L}}) =
F_{n}/{\mathfrak{L}}(F_{n})$ is relatively free; freely generated
by the set $\{x_{1}, \ldots, x_{n}\}$, we obtain
$G_{n}/G^{\prime}_{n} \cong
F_{n}/{\mathfrak{L}}(F_{n})F^{\prime}_{n}$ is an abelian
relatively free group freely generated by the set
$\{x_{1}G^{\prime}_{n}, \ldots,$ $x_{n}G^{\prime}_{n}\}$. Thus
$G_{n}/G^{\prime}_{n}$ is either free abelian with basis
$\{x_{1}G^{\prime}_{n}, \ldots, x_{n}G^{\prime}_{n}\}$ and so,
$F_{n}/{\mathfrak{L}}(F_{n})F^{\prime}_{n}$ has exponent zero,
that is, ${\mathfrak{L}}(F_{n}) \subseteq F^{\prime}_{n}$ and
$G_{n}$ satisfies commutator laws, or it is the abelian group of
exponent $m$ every element of whose basis $\{x_{1}G^{\prime}_{n},
\ldots, x_{n}G^{\prime}_{n}\}$ has order $m$ for some $m > 0$. If
$F_{n}/{\mathfrak{L}}(F_{n})F^{\prime}_{n}$ has exponent $m$ then
$x^{m} \in {\mathfrak{L}}$. Since $G_{n}$ is torsion-free, we
obtain $G_{n}/G^{\prime}_{n}$ is a free abelian of rank $n$ with
basis $\{x_{1}G^{\prime}_{n}, \ldots, x_{n}G^{\prime}_{n}\}$. Thus
$\tau_{2}(G_{n}) = G_n'$. Our next result is about relative
freeness of $L_{K}(G_{n})$.

\begin{theorem}\label{3.27}
The Lie algebra ${\cal L}_{n}$ is isomorphic to $L_{K}(G_{n})$ via
a Lie algebra isomorphism sending $t_{i}$ to
$x_{i}G^{\prime}_{n}$, $i = 1, \ldots, n$. In particular,
$L_{K}(G_{n})$ is relatively free Lie algebra; freely generated by
the set $\{x_{1}G^{\prime}_{n}, \ldots, x_{n}G^{\prime}_{n}\}$.
\end{theorem}

\pf Let $\sigma_{n}$ be the natural mapping from $L_{n}$
onto ${\cal L}_{n}$ sending $\ell_{i}$ to $t_{i}$, $i = 1, \ldots,
n$. By the proof of Theorem \ref{3.18}, we obtain ${\rm ker} \sigma_{n}$
is a fully invariant ideal of $L_{n}$. Let $v \in {\rm
ker}\sigma_{n} \cap L^{m}_{n}$ for some $m \geq 2$, and write $v$
as a linear combination of Lie commutators $[\ell_{i_{1}}, \ldots,
\ell_{i_{m}}]$ i.e. $v = \sum \alpha_{(i_{1}, \ldots, i_{m})}
[\ell_{i_{1}}, \ldots, \ell_{i_{m}}]$. Thus $\sigma_{n}(v) = \sum
\alpha_{(i_{1}, \ldots, i_{m})} [t_{i_{1}}, \ldots, t_{i_{m}}] =
0$. Therefore, for all $c$, $\sum \alpha_{(i_{1}, \ldots, i_{m})}
[t_{i_{1}}(c), \ldots, t_{i_{m}}(c)] = 0$ (as in proof of Theorem
\ref{3.18}). That is,
$$
\sum \alpha_{(i_{1}, \ldots, i_{m})} [\log x_{i_{1},c}, \ldots,
\log x_{i_{m},c}] = 0 \eqno (15)
$$
in ${\cal L}_{K}(G_{n,c})$. By the BCH formula, we obtain
$$
[\log x_{i_{1},c}, \ldots, \log x_{i_{m},c}] = \log (x_{i_{1},c},
\ldots, x_{i_{m},c}) + u,
$$
where $u \in {\cal L}_{K}(\gamma_{m+1}(G_{n,c}))$, and so, we have
from (15)
$$
\sum \alpha_{(i_{1}, \ldots, i_{m})} \log (x_{i_{1},c}, \ldots,
x_{i_{m},c}) \in {\cal L}_{K}(\gamma_{m+1}(G_{n,c})).
$$
Using similar arguments as in the proof of Theorem \ref{3.15}, we
obtain $\sum \alpha_{(i_{1}, \ldots, i_{m})}
[x_{i_{1},c}G^{\prime}_{n,c}, \ldots,$
$x_{i_{m},c}G^{\prime}_{n,c}] = 0$ in $L_{K}(G_{n,c})$ for all
$c$. But ${\cal L}_{K}(G_{n,c})$ is naturally isomorphic to
$L_{K}(G_{n,c})$. By Theorem A (I), ${\cal L}_{K}(G_{n,c})$ is a
relatively free nilpotent Lie algebra and so is $L_{K}(G_{n,c})$.
So, $v$ is an identity in each $L_{K}(G_{n,c})$. Hence, $v$ is an
identity in the cartesian product $L_{K}(G_{n,1}) \times
L_{K}(G_{n,2}) \times \cdots$ therefore, $v$ is an identity in
$\varprojlim L_{K}(G_{n,c})$. Since $L_{K}(G_{n,c}) \cong
L_{K}(G_{n})/\gamma_{c+1}(L_{K}(G_{n}))$, as Lie algebras, in a
natural way, it is easily verified that the \linebreak
$\varprojlim L_{K}(G_{n,c})$ is isomorphic, as a Lie algebra, to
$\varprojlim L_{K}(G_{n})/\gamma_{c+1}(L_{K}(G_{n}))$. Thus $v$ is
an identity in $\varprojlim
L_{K}(G_{n})/\gamma_{c+1}(L_{K}(G_{n}))$. Since $L_{K}(G_{n})$ is
embedded into $\varprojlim
L_{K}(G_{n})/\gamma_{c+1}(L_{K}(G_{n}))$ in a natural way, we have
$v$ is an identity in $L_{K}(G_{n})$. Hence $L_{K}(G_{n}) \in
{\mathfrak{V}}$. By Proposition \ref{3.20} and Lemma \ref{3.2},
$L_{K}(G_{n})$ is generated by the set $\{x_{1}G^{\prime}_{n},
\ldots, x_{n}G^{\prime}_{n}\}$. Let $\varphi$ be the natural
epimorphism from ${\cal L}_{n}$ into $L_{K}(G_{n})$ sending
$t_{i}$ to $x_{i} G^{\prime}_{n}$, $i = 1, \ldots, n$. We claim
that $\varphi$ is one-to-one. By Proposition \ref{3.24}, ${\cal
L}_{n}/\gamma_{c+1}({\cal L}_{n})$ is isomorphic as a Lie algebra
to ${\cal L}_{K}(G_{n,c})$ via an isomorphism $\varphi_{c+1}$ for
all $c$. Thus
$$
\begin{array}{ccl}
{\cal L}_{n}/\gamma_{c+1}({\cal L}_{n}) &
\stackrel{\varphi_{c+1}}{\cong} &
{\cal L}_{K}(G_{n,c}) \\

& \cong & L_{K}(G_{n,c}) \\

 & \cong &
 L_{K}(G_{n})/\gamma_{c+1}(L_{K}(G_{n}))
\end{array}
 $$
for all $c$. Since ${\cal L}_{n}$ is residually nilpotent, i.e.
$\cap_{c \geq 1} \gamma_{c+1}({\cal L}_{n}) = \{0\}$, and for all
$c$, ${\cal L}_{n}/\gamma_{c+1}({\cal L}_{n}) \cong
 L_{K}(G_{n})/\gamma_{c+1}(L_{K}(G_{n}))$,
 we obtain ${\cal L}_{n}$ is isomorphic as a Lie algebra
 to $L_{K}(G_{n})$. By Theorem \ref{3.18}, we obtain the Lie
 algebra $L_{K}(G_{n})$ is relatively free. It is
 easily verified that the set $\{x_{1}G^{\prime}_{n}, \ldots,
 x_{n}G^{\prime}_{n}\}$ freely generates
 $L_{K}(G_{n})$. \qed

\section{Proofs of Theorems A and B}

\subsection{Proof of Theorem A} (I) It follows from
Proposition \ref{3.14} and Theorem \ref{3.15} that ${\cal
L}_{K}(F_{n}({\mathfrak{T}}_{c}))$ is relatively free in some
variety of nilpotent Lie algebras, and ${\cal
L}_{K}(F_{n}({\mathfrak{T}}_{c})) \cong
L_{K}(F_{n}({\mathfrak{T}}_{c}))$ as Lie algebras in a natural
way. Write $G = F_{n}({\mathfrak{T}}_{c})$, with $n \geq 2$,
freely generated by the set $\{y_{1}, \ldots, y_{n}\}$. From the
proof of Theorem 3, $G/G^{\prime}$ is free abelian group of rank
$n$, and so, $\tau_{2}(G) = G^{\prime}$. Recall that $L^{(S)}(G) =
\oplus_{1 \leq i \leq c} L^{(S)}_{i}(G)$, where $L^{(S)}_{i}(G) =
\gamma_{i}(G) \tau_{i+1}(G)/\tau_{i+1}(G)$, $i = 1, \ldots, c$.
The additive group of $L^{(S)}(G)$ is free abelian,
$L_{\mathbb{Q}}(G) = {\mathbb{Q}}
\otimes_{{\mathbb{Z}}}L^{(S)}(G)$ and $L^{(S)}(G)$ is a subset of
$L_{\mathbb{Q}}(G)$. By Lemma 3.2, $L^{(S)}(G)$ is generated as a
Lie ring by the set $\{y_{1}G^{\prime}, \ldots,
y_{n}G^{\prime}\}$. We give on $L_{\mathbb{Q}}(G)$ the structure
of a group by means of BCH formula, denoted by $R$. (Notice that
$R = L_{\mathbb{Q}}(G)$ as sets.) By Lemma 3.3,
$\{y_{1}G^{\prime}+L_{\mathbb{Q}}(G)^{\prime}, \ldots,
y_{n}G^{\prime}+L_{\mathbb{Q}}(G)^{\prime}\}$ is a
$\mathbb{Q}$-basis for
$L_{\mathbb{Q}}(G)/L_{\mathbb{Q}}(G)^{\prime}$. Let $H$ be the
subgroup of $R$ generated by the set $\{y_{1}G^{\prime}, \ldots,
y_{n}G^{\prime}\}$. (Notice that the identity element of $H$ is
the zero element in $L_{\mathbb{Q}}(G)$.) By Lemma 2.1, $H$ is a
torsion-free finitely generated nilpotent group of class $c$ and
$\tau_{2}(H) = H^{\prime}$.

Next we shall show that there is a group isomorphism
$$
\gamma_{i}(H)/\gamma_{i+1}(H) \cong
\gamma_{i}(L^{(S)}(G))/\gamma_{i+1}(L^{(S)}(G)) \eqno (16)
$$
for all $i$, with $i = 1, \ldots, c$. (For the proof of the
aforementioned isomorphism, we use some arguments given in
\cite{shmelkin}.) First we shall show that $\gamma_{c}(H) =
\gamma_{c}(L^{(S)}(G))$. Every element of $\gamma_{c}(H)$ can be
written as a product of group commutators (in the sense of the
operation $\circ$) of length $c$ in the $y_{1}G^{\prime}, \ldots,
y_{n}G^{\prime}$. By the BCH formula, we may deduce that every
group commutator of length $c$ in the $y_{1}G^{\prime}, \ldots,
y_{n}G^{\prime}$ lies in $L_{c,{\mathbb{Q}}}(G)$ $(= {\mathbb{Q}}
\otimes_{\mathbb{Z}} \gamma_{c}(G))$. Since
$\gamma_{c+1}(L^{(S)}(G)) = \{0\}$, the multiplication of group
commutators of length $c$ in $H$ is equal to their addition in the
ring $L^{(S)}(G)$. Hence it follows that $\gamma_{c}(H) \subseteq
\gamma_{c}(L^{(S)}(G))$. Conversely, if $u = u(y_{1}G^{\prime},
\ldots, y_{n}G^{\prime})$ is a Lie commutator of length $c$ in the
$y_{1}G^{\prime}, \ldots, y_{n}G^{\prime}$ then (since
$\gamma_{c+1}(L^{(S)}(G)) = \{0\}$) it is equal to the group
commutator in the $y_{1}G^{\prime}, \ldots, y_{n}G^{\prime}$
obtained from $u$ by replacing the operation on ring
multiplication by the operation of commutation in the group $H$.
Since every element of $\gamma_{c}(L^{(S)}(G))$ is a linear
combination with integer coefficients of Lie commutators of length
$c$ in the $y_{1}G^{\prime}, \ldots, y_{n}G^{\prime}$, we have the
inverse inclusion $\gamma_{c}(L^{(S)}(G)) \subseteq
\gamma_{c}(H)$. Therefore, we get $\gamma_{c}(H) =
\gamma_{c}(L^{(S)}(G)) = \gamma_{c}(G)$. Write $G^{(c-1)} =
G/\tau_{c}(G)$. Thus $G^{(c-1)}$ is a relatively free nilpotent
torsion-free group of rank $n$ and class $c-1$. It is easy to
verify that $\tau_{i+1}(G^{(c-1)}) = \tau_{i+1}(G)/\tau_{i}(G)$
for $i = 1, \ldots, c-1$. Since $(G^{(c-1)})^{\prime} =
G^{\prime}/\tau_{c}(G)$, we get $G^{(c-1)}/(G^{(c-1)})^{\prime}$
is a free abelian group of rank $n$. Moreover $G^{(c-1)}$ is
generated by the set $\{\overline{y}_{1}, \ldots,
\overline{y}_{n}\}$, where $\overline{y}_{i} = y_{i} \tau_{c}(G)$,
$i = 1, \ldots, n$  and $\tau_{2}(G^{(c-1)}) =
G^{\prime}/\tau_{c}(G) = (G^{(c-1)})^{\prime}$. Further
$L^{(S)}(G^{(c-1)})$ is free abelian, $L_{\mathbb{Q}}(G^{(c-1)}) =
{\mathbb{Q}} \otimes_{\mathbb{Z}} L^{(S)}(G^{(c-1)})$, and
$L^{(S)}(G^{(c-1)})$ is regarded as a subset of
$L_{\mathbb{Q}}(G^{(c-1)})$. By Lemma 3.2, $L^{(S)}(G^{(c-1)})$ is
generated as a Lie ring by the set
$\{\overline{y}_{1}(G^{(c-1)})^{\prime}, \ldots,
\overline{y}_{n}(G^{(c-1)})^{\prime}\}$. It is easy to check that
$\{\overline{y}_{1}(G^{(c-1)})^{\prime}, \ldots,
\overline{y}_{n}(G^{(c-1)})^{\prime}\}$ is a $\mathbb{Z}$-basis
for $G^{(c-1)}/(G^{(c-1)})^{\prime}$ and so, by Lemma 3.3,
$\{\overline{y}_{1}(G^{(c-1)})^{\prime}+L_{\mathbb{Q}}(G^{(c-1)})^{\prime},
\ldots,
\overline{y}_{n}(G^{(c-1)})^{\prime}+L_{\mathbb{Q}}(G^{(c-1)})^{\prime}\}$
is a $\mathbb{Q}$-basis for
$L_{\mathbb{Q}}(G^{(c-1)})/L_{\mathbb{Q}}(G^{(c-1)})^{\prime}$. As
before, we give on $L_{\mathbb{Q}}(G^{(c-1)})$ the structure of a
group by means of the BCH formula, denoted by $R^{(c-1)}$. Let
$H^{(c-1)}$ be the subgroup of $R^{(c-1)}$ generated by the set
$\{\overline{y}_{1}(G^{(c-1)})^{\prime}, \ldots,
\overline{y}_{n}(G^{(c-1)})^{\prime}\}$. By Lemma 2.1, $H^{(c-1)}$
is a torsion-free finitely generated nilpotent group and
$\tau_{2}(H^{(c-1)}) = (H^{(c-1)})^{\prime}$. Using similar
arguments as before, $\gamma_{c-1}(H^{(c-1)}) =
\gamma_{c-1}(L^{(S)}(G^{(c-1)})$. In the next few lines, we write
$b_{i} = y_{i}G^{\prime}$ and $b^{(c-1)}_{i} =
\overline{y}_{i}(G^{(c-1)})^{\prime}$, $i = 1, \ldots, n$. As the
group operation $\circ$ can be expressed in terms of the Lie
algebra operation, the natural Lie algebra epimorphism
$\alpha_{L}$ from $L_{\mathbb{Q}}(G)$ onto
$L_{\mathbb{Q}}(G^{(c-1)})$ induces a group epimorphism
$\alpha_{H}$ from $H$ onto $H^{(c-1)}$ such that
$\alpha_{H}(y_{i}G^{\prime}) = \alpha_{L}(y_{i}G^{\prime}) =
\overline{y}_{i}(G^{(c-1)})^{\prime}$, $i = 1, \ldots, n$. We
claim that the kernel of $\alpha_{H}$, ${\rm ker} \alpha_{H}$, is
equal to $\gamma_{c}(H)$. Since $\gamma_{c}(H) =
\gamma_{c}(L^{(S)}(G)) = \gamma_{c}(G)$ and $\gamma_{c}(G)
\subseteq {\rm ker} \alpha_{H}$, it is enough to show that ${\rm
ker} \alpha_{H} \subseteq \gamma_{c}(G)$. Let $v \in {\rm ker}
\alpha_{H}$ and let $\mu$ be the smallest natural integer such
that $v \in \gamma_{\mu}(H) \setminus \gamma_{\mu + 1}(H)$. Write
$$
v = (\prod_{{\rm finite}} (b_{i_{1}}, \ldots,
b_{i_{\mu}})^{a_{(i_{1}, \ldots, i_{\mu})}}) v^{\prime}(b_{1},
\ldots, b_{n}),
$$
where $a_{(i_{1}, \ldots, i_{\mu})} \in {\mathbb{Z}}$ and
$v^{\prime}(b_{1}, \ldots, b_{n}) \in \gamma_{\mu + 1}(H)$. Since
$v \in {\rm ker} \alpha_{H}$, we obtain
$$
(\prod_{{\rm finite}} (b^{(c-1)}_{i_{1}}, \ldots,
b^{(c-1)}_{i_{\mu}})^{a_{(i_{1}, \ldots, i_{\mu})}})
v^{\prime}(b^{(c-1)}_{1}, \ldots, b^{(c-1)}_{n}) = 0.
$$
By the BCH formula,
$$
\sum a_{(i_{1}, \ldots, i_{\mu})} [b^{(c-1)}_{i_{1}}, \ldots,
b^{(c-1)}_{i_{\mu}}] \in \gamma_{\mu +
1}(L_{\mathbb{Q}}(G^{(c-1)})).
$$
Since $L_{\mathbb{Q}}(G^{(c-1)}) = \oplus_{1 \leq i \leq c-1}
({\mathbb{Q}} \otimes L^{(S)}_{i}(G^{(c-1)}))$ and, for $t = 1,
\ldots, c-1$, $\gamma_{t}(L_{\mathbb{Q}}(G^{(c-1)}) = \oplus_{i
\geq t} ({\mathbb{Q}} \otimes L^{(S)}_{i}(G^{(c-1)}))$, we have
$$
\sum a_{(i_{1}, \ldots, i_{\mu})} [b^{(c-1)}_{i_{1}}, \ldots,
b^{(c-1)}_{i_{\mu}}] = 0
$$
and so,
$$
\sum a_{(i_{1}, \ldots, i_{\mu})} [b_{i_{1}}, \ldots, b_{i_{\mu}}]
\in {\rm ker} \alpha_{L},
$$
where ${\rm ker} \alpha_{L}$ denotes the kernel of $\alpha_{L}$.
Suppose that $\mu \neq c$. Since ${\rm ker} \alpha_{L} =
{\mathbb{Q}} \otimes_{\mathbb{Z}} \gamma_{c}(G)$ and
$L_{\mathbb{Q}}(G)$ is graded, we obtain
$$
\sum a_{(i_{1}, \ldots, i_{\mu})} [b_{i_{1}}, \ldots, b_{i_{\mu}}]
= 0~~ ({\rm in}~~L_{\mathbb{Q}}(G)).
$$
Thus
$$
\prod_{\rm finite} (y_{i_{1}}, \ldots, y_{i_{\mu}})^{a_{(i_{1},
\ldots, i_{\mu})}} \in \gamma_{\mu + 1}(G).
$$
Recall from the proof of Theorem \ref{3.15} that $\eta$ is the Lie
algebra isomorphism from ${\cal L}_{\mathbb{Q}}(G)$ into
$L_{{\mathbb{Q}}}(G)$ satisfying the conditions $\eta(\log y_{i})
= y_{i}G^{\prime}$ for $i = 1, \ldots, n$. Let $\xi_{\eta}$ be the
mapping from the group $\exp {\cal L}_{\mathbb{Q}}(G)$ to the
group $R$ defined by $\xi_{\eta}(\exp u) = \eta(u)$ for all $u \in
{\cal L}_{\mathbb{Q}}(G)$. Since $\eta$ is a Lie algebra
homomorphism and the group operation $\circ$ is expressed in terms
of Lie commutators, we have $\xi_{\eta}((\exp u)(\exp v)) =
\eta(u) \circ \eta(v)$ for all $u, v \in {\cal
L}_{\mathbb{Q}}(G)$. It is easily verified that $\xi_{\eta}$ is a
group isomorphism. But $\xi_{\eta}(y_{i}) = \eta(\log y_{i}) =
b_{i}$ for $i = 1, \ldots, n$. Therefore $\xi_{\eta}(G) = H$ and
so, for $t = 1, \ldots, c$, $\xi_{\eta}(\gamma_{t+1}(G)) =
\gamma_{t+1}(H)$. Hence
$$
\prod_{\rm finite} (b_{i_{1}}, \ldots, b_{i_{\mu}})^{a_{(i_{1},
\ldots, i_{\mu})}} \in \gamma_{\mu + 1}(H),
$$
which is a contradiction. Therefore, $\mu = c$ and $v \in
\gamma_{c}(H) = \gamma_{c}(G)$. Hence ${\rm ker} \alpha_{H}
\subseteq \gamma_{c}(G)$ and ${\rm ker} \alpha_{H} = \gamma_{c}(G)
= \gamma_{c}(H)$. Thus
$$
\gamma_{c-1}(H^{(c-1)}) = \gamma_{c-1}(H)/\gamma_{c}(H) \cong
\gamma_{c-1}(L^{(S)}(G))/\gamma_{c}(L^{(S)}(G)).
$$
Eventually we see that the isomorphism (16) holds for every $i
\leq c$. Since
$$
\gamma_{i}(L^{(S)}(G))/\gamma_{i+1}(L^{(S)}(G)) \cong
\gamma_{i}(G)\tau_{i+1}(G)/\tau_{i+1}(G)
$$
for $i = 1, \ldots, c$, we obtain from (16) that $H$ is Magnus.
Since $G \cong H$ by means of $\xi_{\eta}$, we have the required
result.

(II) Since $H \cong F_{n,c}/{\rm ker}\tau$ and ${\rm ker}\tau$ is
fully invariant (by Proposition \ref{3.16} (iii)), we obtain $H$
is relatively free of finite rank. Furthermore, since $H$ is
torsion-free nilpotent, we obtain from (I) that $H$ is Magnus. By
Lemma \ref{3.11}, we obtain ${\cal L}_{\mathbb{Q}}({\rm ker} \tau)
\subseteq {\cal L}_{\mathbb{Q}}(F_{n,c})^{\prime}$. By Lemma
\ref{2.1} and (I), $L \cong {\cal L}_{\mathbb{Q}}(H) \cong
L_{{\mathbb{Q}}}(H)$ as Lie algebras. \qed

\vskip .120 in

\begin{remark}\label{4.1}
{\upshape Let $L$ be a relatively free nilpotent Lie algebra over
${\mathbb{Q}}$ of finite rank $n$, with $n \geq 2$. Let $\{h_{1},
\ldots, h_{n}\}$ be a free generating set of $L$. Then
$\{h_{1}+L^{\prime}, \ldots, h_{n}+L^{\prime}\}$ is a
${\mathbb{Q}}$-basis for $L$. Let $y_{1}, \ldots, y_{n}$ be
elements of $L$ such that the set $\{y_{1}+L^{\prime}, \ldots,
y_{n}+L^{\prime}\}$ is a ${\mathbb{Q}}$-basis of $L$. For each
$j$, with $j = 1, \ldots, n$,
$$
h_{j} = \sum_{i=1}^{n} \alpha_{ij} y_{i} + v_{i},
$$
where $v_{i} \in L^{\prime}$, $i = 1, \ldots, n$, and $\alpha_{ij}
\in {\mathbb{Q}}$. It is easily verified that $L$ is generated by
the set $\{y_{1}, \ldots, y_{n}\}$. Let $\varphi$ be the map from
$L$ into $L$ satisfying the conditions $\varphi(h_{j}) = y_{j}$,
$j = 1, \ldots, n$. Since $L$ is relatively free on $\{h_{1},
\ldots, h_{n}\}$ and $L$ is generated by the set $\{y_{1}, \ldots,
y_{n}\}$, $\varphi$ extends uniquely to a Lie algebra epimorphism
of $L$. Since $\varphi$ induces a group automorphism on
$L/L^{\prime}$ and $L$ is nilpotent, it is easily checked that
$\varphi$ is an automorphism of $L$. Thus the set $\{y_{1},
\ldots, y_{n}\}$ is a free generating set of $L$. Consider $L$ as
a group, denoted $R$, by means of the BCH formula. Let $H_{1}$ and
$H_{2}$ be the subgroups of $R$ generated by the sets $\{h_{1},
\ldots, h_{n}\}$ and $\{y_{1}, \ldots, y_{n}\}$, respectively. By
the proof of Theorem A (II) and since both $H_{1}$ and $H_{2}$
have rank $n$, we get $H_{1} \cong H_{2}$. Hence for a relatively
free nilpotent Lie algebra of finite rank $n$ over ${\mathbb{Q}}$,
we associate (via BCH formula) a unique (up to isomorphism)
relatively free Magnus nilpotent group of rank $n$. }
\end{remark}

\vskip .120 in

\emph{Proof of Corollary \ref{1.1}.}  Let $G$ be a torsion-free
finitely generated nilpotent group of class $c$, and let $K$ be a
field of characteristic zero. By Lemma \ref{3.8}, we obtain ${\rm
grad}^{(\ell)}({\cal L}_{K}(G)) \cong L_{K}(G)$ as Lie algebras in
a natural way. Write $R_{K}(G) = \exp {\cal L}_{K}(G)$, and give
on  ${\cal L}_{K}(G)$ the structure of a group by means of the BCH
formula. Then $({\cal L}_{K}(G), \circ)$ is isomorphic to
$R_{K}(G)$ by a group isomorphism sending $u$ to $\exp u$ for all
$u \in {\cal L}_{K}(G)$. Thus $\exp \gamma_{t}({\cal L}_{K}(G)) =
\gamma_{t}(R_{K}(G))$ for $t = 1, \ldots, c$. Form the direct sum
of the abelian groups
$$
{\rm grad}^{(g)}(R_{K}(G)) = \oplus_{t=1}^{c}
\gamma_{t}(R_{K}(G))/\gamma_{t+1}(R_{K}(G))
$$
and give it the structure of a Lie algebra by defining a Lie
multiplication
$$
[u\gamma_{i+1}(R_{K}(G)), v\gamma_{j+1}(R_{K}(G))] =
(u,v)\gamma_{i+j+1}(R_{K}(G))
$$
for $u \in \gamma_{i}(R_{K}(G))$, $v \in \gamma_{j}(R_{K}(G))$,
$i, j \in \{1, \ldots, c\}$. Extend this multiplication to
\linebreak ${\rm grad}^{(g)}(R_{K}(G))$ by linearity. Similarly,
we form the direct sum of the abelian groups
$$
{\rm grad}^{(g)}({\cal L}_{K}(G)) = \oplus_{t=1}^{c}
\gamma_{t}({\cal L}_{K}(G))/\gamma_{t+1}({\cal L}_{K}(G))
$$
and give it the structure of a Lie algebra. Since $({\cal
L}_{K}(G), \circ) \cong R_{K}(G)$ as groups, we have
$$
{\rm grad}^{(g)}(R_{K}(G)) \cong {\rm grad}^{(g)}({\cal L}_{K}(G))
$$
as Lie algebras in a natural way. Let $\{a_{1}, \ldots,
a_{f(c+1)}\}$ be a canonical basis for $G$. Then the set $\{\log
a_{1}, \ldots, \log a_{f(c+1)}\}$ is a $K$-basis of ${\cal
L}_{K}(G)$. By Lemma \ref{3.6} and the proof of Lemma \ref{3.7},
$\gamma_{t}({\cal L}_{K}(G))/\gamma_{t+1}({\cal L}_{t+1}(G))$ has
a $K$-basis the set $\{\log a_{f(c+1)} + \gamma_{t+1}({\cal
L}_{K}(G)), \ldots, \log a_{f(c+1)} + \gamma_{t+1}({\cal
L}_{K}(G))\}$. Using the BCH formula,
\begin{flushleft}
$(\lambda_{1} \log a_{f(t)+1} + \cdots + \lambda_{n_{t}} \log
a_{f(t+1)}) + \gamma_{t+1}({\cal L}_{K}(G)) =$
\end{flushleft}
\begin{flushright}
$(\lambda_{1} \log a_{f(t)+1} \circ \cdots \circ \lambda_{n_{t}}
\log a_{f(t+1)}) \circ \gamma_{t+1}({\cal L}_{K}(G))$
\end{flushright}
(as sets). Therefore every element of $({\cal L}_{K}(G), \circ)$
is written uniquely as $\lambda_{1} \log a_{1} \circ \cdots \circ
\lambda_{f(c+1)} \log a_{f(c+1)}$, where $\lambda_{1}, \ldots,
\lambda_{f(c+1)} \in K$. Hence
$$
{\rm grad}^{(g)}({\cal L}_{K}(G)) = {\rm grad}^{(\ell)}({\cal
L}_{K}(G))
$$
as Lie algebras. Since $({\cal L}_{K}(G)), \circ) \cong R_{K}(G)$,
we have every element of $R_{K}(G)$ is written uniquely as
$\exp(\lambda_{1} \log a_{1}) \cdots \exp(\lambda_{f(c+1)} \log
a_{f(c+1)})$ with $\lambda_{1}, \ldots, \lambda_{f(c+1)} \in K$
and so,
$$
{\rm grad}^{(g)}(R_{K}(G)) \cong {\rm grad}^{(\ell)}({\cal
L}_{K}(G))
$$
as Lie algebras in a natural way. Since ${\rm grad}^{(\ell)}({\cal
L}_{K}(G)) \cong L_{K}(G)$ (as Lie algebras), we have ${\rm
grad}^{(g)}(R_{K}(G)) \cong L_{K}(G)$ as Lie algebras in a natural
way.

(II) Suppose that $G$ is relatively free. By Theorem \ref{3.15},
we obtain ${\cal L}_{K}(G) \cong L_{K}(G)$ as Lie algebras and so,
by (I), ${\cal L}_{K}(G) \cong {\rm grad}^{(g)}(\exp {\cal
L}_{K}(G))$.

\vskip .120 in

\begin{remark}\label{4.2}
{\upshape
It is easy to verify that
$$
K \otimes_{{\mathbb{Q}}} {\rm grad}^{(g)}(\exp {\cal
L}_{{\mathbb{Q}}}(G)) \cong {\rm grad}^{(g)}(\exp {\cal
L}_{K}(G)).
$$
as Lie algebras, where $K$ is a field of characteristic zero. For $K = {\mathbb{R}}$, $\exp {\cal
L}_{\mathbb{R}}(G)$ is a real simply connected Lie group whose rational Lie
algebra is ${\cal L}_{\mathbb{Q}}(G)$. Moreover, if $\{a_{1},
\ldots, a_{f(c+1)}\}$ is a canonical basis of $G$, then every
element of $\exp {\cal L}_{\mathbb{R}}(G)$ is written uniquely as
$$\exp(\lambda_{1} \log a_{1})\cdots\exp(\lambda_{f(c+1)} \log
a_{f(c+1)})$$ with $\lambda_{1}, \ldots, \lambda_{f(c+1)} \in
{\mathbb{R}}$. Next, we recall a standard procedure of a
construction of a Lie group from a finite-dimensional nilpotent
Lie algebra $L$ over $\mathbb{Q}$. Let $b_{1}, \ldots, b_{n} \in
L$ such that the set $\{b_{1} + L^{\prime}, \ldots, b_{n} +
L^{\prime}\}$ is a $\mathbb{Q}$-basis of $L/L^{\prime}$. Give on
$L$ the structure of a group via the BCH formula. Let $H$ be the
subgroup of $(L, \circ)$ generated by the set $\{b_{1}, \ldots,
b_{n}\}$. Then $(L, \circ)$ is a Mal'cev completion of $H$, and
$({\mathbb{R}} \otimes_{{\mathbb{Q}}} L, \circ)$ is a real simply
connected Lie group containing $H$ as a discrete subgroup with
rational Lie algebra $L$. }
\end{remark}

\emph{Proof of Corollary \ref{1.2}.} Let $G_{1}, G_{2}$ be
torsion-free finitely generated nilpotent groups which are
quasi-isometric. By Remark 4.2, $\exp {\cal
L}_{\mathbb{R}}(G_{j})$, with $j = 1, 2$ is a real simply
connected Lie group whose rational Lie algebra is ${\cal
L}_{\mathbb{Q}}(G_{j})$. It follows from a result of Pansu
\cite[Theorem 3]{pansu} that
$$
{\rm grad}^{(g)}(\exp {\cal L}_{\mathbb{R}}(G_{1})) \cong {\rm
grad}^{(g)}(\exp {\cal L}_{\mathbb{R}}(G_{2}))
$$
as Lie algebras. Suppose that $G_{j}$ ($j = 1, 2$) is a relatively
free group of finite rank. By Theorem A (I), $G_{j}$ is a Magnus
group, with $j = 1, 2$. By Proposition \ref{3.14} and Theorem
\ref{3.15}, ${\cal L}_{\mathbb{R}}(G_{j})$ is relatively free and
${\cal L}_{\mathbb{R}}(G_{j}) \cong L_{\mathbb{R}}(G_{j})$ as Lie
algebras for $j = 1, 2$. By Corollary \ref{1.1}, ${\cal
L}_{\mathbb{R}}(G_{2}) \cong {\cal L}_{\mathbb{R}}(G_{1})$ as Lie
algebras. Hence both $G_{2}$ and $G_{1}$ have the same finite rank
$n$ and nilpotency class $c$, and
$\gamma_{i}(G_{2})/\gamma_{i+1}(G_{2}) \cong
\gamma_{i}(G_{1})/\gamma_{i+1}(G_{1})$, $i = 1, \ldots, c$. Let
$F_{n,c}$ be the free nilpotent group of rank $n$ and class $c$;
freely generated by the set $\{x_{1}, \ldots, x_{n}\}$. For $j =
1, 2$, let $\{y_{1j}, \ldots, y_{nj}\}$ be a free generating set
for $G_{j}$. Furthermore, we write $\pi_{j}$ for the natural group
epimorphism from $F_{n,c}$ onto $G_{j}$ such that $\pi_{j}(x_{i})
= y_{ij}$, $i = 1, \ldots, n$. Thus $G_{j} \cong F_{n,c}/{\rm ker}
\pi_{j}$, with $j = 1, 2$. By Lemma \ref{3.4} and the proof of
Theorem \ref{3.15}, the set $\{\log y_{1j}, \ldots, \log y_{nj}\}$
is a free generating set for ${\cal L}_{\mathbb{Q}}(G_{j})$. We
claim ${\cal L}_{\mathbb{Q}}(G_{2}) \cong {\cal
L}_{\mathbb{Q}}(G_{1})$ as Lie algebras. To get a contradiction we
assume that ${\cal L}_{\mathbb{Q}}(G_{2}) \ncong {\cal
L}_{\mathbb{Q}}(G_{1})$, and let $w$ be a word (for Lie algebras
over $\mathbb{Q}$) such that $w(\log y_{11}, \ldots, \log y_{n1})
= 0$ and $w(\log y_{12}, \ldots, \log y_{n2}) \neq 0$. Since
${\cal L}_{\mathbb{R}}(G_{j}) = {\mathbb{R}} \otimes_{\mathbb{Q}}
{\cal L}_{\mathbb{Q}}(G_{j})$, $j = 1, 2$, and ${\cal
L}_{\mathbb{R}}(G_{1}) \cong {\cal L}_{\mathbb{R}}(G_{2})$ (as Lie
algebras), we obtain $w$ is an identity in ${\cal
L}_{\mathbb{R}}(G_{2})$ which is a contradiction. Since both
${\cal L}_{\mathbb{Q}}(G_{1})$ and ${\cal L}_{\mathbb{Q}}(G_{2})$
are relatively free of finite rank $n$ in the same variety, we get
${\cal L}_{\mathbb{Q}}(G_{1}) \cong {\cal L}_{\mathbb{Q}}(G_{2})$.
Since, $G_1$ and $G_2$ are relatively free, in order to prove that
$G_{1}$ is isomorphic to $G_{2}$, it is enough to show that ${\rm
ker} \pi_{1} \subseteq {\rm ker}\pi_{2}$. Let $\widetilde{w} \in
{\rm ker} \pi_{1} \subseteq F^{\prime}_{n,c}$. By Lemma \ref{3.10}
(for $M = G_{1}$), the natural epimorphism $\pi_{1}$ gives rise a
Lie algebra homomorphism $\psi_{\pi_{1}}$ from ${\cal
L}_{\mathbb{Q}}(F_{n,c})$ into ${\cal L}_{\mathbb{Q}}(G_{1})$ such
that $\psi_{\pi_{1}}(\log u) = \log \pi_{1}(u)$ for all $u \in
{\rm ker} \pi_{1}$. Thus $\psi_{\pi_{1}}(\log \widetilde{w}) = 0$
and so, $\log \widetilde{w} \in {\cal L}_{\mathbb{Q}}({\rm
ker}\pi_{1})$. Thus
$$
\log \widetilde{w} = \widetilde{w}_{2} + \cdots +
\widetilde{w}_{c},
$$
where $\widetilde{w}_{t} \in {\cal L}_{\mathbb{Q}}({\rm
ker}\pi_{1}) \cap {\cal L}_{t}(F_{n,c})$ (by Lemma \ref{3.12}) for
$t = 2, \ldots, c$. By Proposition \ref{3.13}, $\widetilde{w}_{2},
\ldots, \widetilde{w}_{c}$ are identities in ${\cal
L}_{\mathbb{Q}}(G_{1})$. Since ${\cal L}_{\mathbb{Q}}(G_{1}) \cong
{\cal L}_{\mathbb{Q}}(G_{2})$, we have $\widetilde{w}_{2}, \ldots,
\widetilde{w}_{c}$ are identities for ${\cal
L}_{\mathbb{Q}}(G_{2})$ as well. Therefore $\widetilde{w}_{t} \in
{\cal L}_{\mathbb{Q}}({\rm ker}\pi_{2}) \cap {\cal
L}_{t}(F_{n,c})$ for $t = 2, \ldots, c$. Hence $ \log
\widetilde{w} \in {\cal L}_{\mathbb{Q}}({\rm ker}\pi_{2})$. Thus
$\psi_{\pi_{2}}(\log \widetilde{w}) = 0$ and so, $\widetilde{w}
\in {\rm ker}\pi_{2}$. Therefore, there exists a group epimorphism
from $G_{1}$ onto $G_{2}$. Since
$\gamma_{i}(G_{1})/\gamma_{i+1}(G_{1}) \cong
\gamma_{i}(G_{2})/\gamma_{i+1}(G_{2})$, $i = 1, \ldots, c$, we
obtain $G_{1} \cong G_{2}$. \qed

\subsection{Proof of Theorem B.} (I) It has been proved in
Theorems \ref{3.18} and \ref{3.27}.

(II) Let $L$ be a relatively free Lie algebra over ${\mathbb{Q}}$
of rank $n$, with $n \geq 2$. If $L$ is nilpotent, then our claim
follows from Theorem A (II). Thus we assume that $L$ is not
nilpotent. Let ${\mathfrak{V}}$ be a variety of Lie algebras such
that $L$ is relatively free of rank $n$. Therefore $L \cong
L_{n}/{\mathfrak{V}}(L_{n})$, where $L_{n}$ is the absolutely free
Lie algebra freely generated by the set $\{\ell_{1}, \ldots,
\ell_{n}\}$. Without loss of generality, we may write $L =
L_{n}/{\mathfrak{V}}(L_{n})$. Let $y_{i} = \ell_{i} +
{\mathfrak{V}}(L_{n})$, $i = 1, \ldots, n$. Thus the set $\{y_{1},
\ldots, y_{n}\}$ is a free generating set of $L$. For a positive
integer $m$, let $L^{m}_{n}$ be the subspace of $L_{n}$ spanned by
all Lie commutators of total degree $m$ in $\ell_{1}, \ldots,
\ell_{n}$. Since ${\mathbb{Q}}$ is infinite, $L_{n} \cap
{\mathfrak{V}}(L_{n}) = \oplus_{m \geq 1} (L^{m}_{n} \cap
{\mathfrak{V}}(L_{n}))$. Thus we may write $L$ as a sum of
homogeneous components, $L = \oplus_{m \geq 1} L^{m}$, where
$L^{m} \cong L^{m}_{n}/(L^{m}_{n} \cap {\mathfrak{V}}(L_{n}))$ and
$L^{m}$ is the subspace of $L$ spanned by all Lie commutators of
total degree $m$ in $y_{1}, \ldots, y_{n}$. Each element $u$ of
$L$ may be uniquely written in the form $u = \sum_{m \geq 1}
u_{m}$ with $u_{m} \in L^{m}$ for all $m$ and $u_{m} = 0$ for all
but finitely many $m$. Furthermore, for $c \geq 1$, $\gamma_{c}(L)
= \oplus_{m \geq c} L^{m}$. We write $L(c) = L/\gamma_{c+1}(L)$,
and $y_{i,c} = y_{i} + \gamma_{c+1}(L)$, $i = 1, \ldots, n$. Since
$L$ is relatively free, then $L(c)$ is a relatively free nilpotent
Lie algebra of rank $n$ and class $c$ with a free generating set
$\{y_{1,c}, \ldots, y_{n,c}\}$. Notice that
$\{y_{1,c}+L(c)^{\prime}, \ldots, y_{n,c}+L(c)^{\prime}\}$ is a
basis of $L(c)/L(c)^{\prime}$. Give on $L(c)$ the structure of a
group, denoted by $R_{c}$, by means of the BCH formula, and let
$Y_{c}$ be the subgroup of $R_{c}$ generated by the set
$\{y_{1,c}, \ldots, y_{n,c}\}$. By Theorem A (II), $Y_{c}$ is a
finitely generated Magnus nilpotent group of class $c$. By the BCH
formula, $\gamma_{c}(L(c))$ is spanned by all group commutators
$(y_{i_{1},c}, \ldots, y_{i_{c},c})$ with $i_{1}, \ldots, i_{c}
\in \{1, \ldots, n\}$. For positive integers $c$ and $d$, with $c
\leq d$, let $\rho_{c,d}$ be the natural Lie algebra epimorphism
from $L(d)$ onto $L(c)$ sending $y_{i,d}$ to $y_{i,c}$ for $i = 1,
\ldots, n$. As the group operation in $Y_{c}$ (for all $c$) can be
expressed in terms of the Lie algebra operations, $\rho_{c,d}$
induces a group homomorphism, say $\widetilde{\rho}_{c,d}$, from
$Y_{d}$ onto $Y_{c}$ such that $\widetilde{\rho}_{c,d}$ sends
$y_{i,d}$ to $y_{i,c}$ for $i = 1, \ldots, n$. It is clearly
enough that $Y_{d}/\gamma_{d}(Y_{d}) \cong Y_{d-1}$. Since
$\gamma_{d}(L(d))$ is spanned by all group commutators
$(y_{i_{1},d}, \ldots, y_{i_{d},d})$, with $i_{1}, \ldots, i_{d}
\in \{1, \ldots, n\}$, and ${\rm ker} \rho_{(d-1),d} =
\gamma_{d}(L(d))$, we obtain $\gamma_{d}(Y_{d}) \subseteq {\rm
ker} \widetilde{\rho}_{(d-1),d}$. Since $Y_{d}/{\rm ker}
\widetilde{\rho}_{(d-1),d} \cong Y_{d-1}$ and each $Y_{d}$ is
torsion-free finitely generated  nilpotent, we get the kernel of
$\widetilde{\rho}_{(d-1),d}$ is equal to $\gamma_{d}(Y_{d})$ for
all $d \geq 2$.

Let $\widehat{L} = \varprojlim L(c)$ be the completion of $L$ with
respect to the lower central series. ($\varprojlim L(c)$ may be
identified with the complete (unrestricted) direct sum
$\widehat{\oplus}_{m \geq 1} L^{m}$, and it has a natural Lie
algebra structure.) Moreover, $L$ is naturally contained in
$\widehat{L}$. Give on $\widehat{L}$ the structure of a group,
denoted by $\widehat{R}$, via the BCH formula. That is,
$\widehat{R} = (\widehat{L}, \circ) = \varprojlim R(c)$. Let $H$
be the subgroup of $\widehat{R}$ generated by the set $\{y_{1},
\ldots, y_{n}\}$. Notice that, for $i = 1, \ldots, n$,
$$
y_{i} = (y_{i,1}, y_{i,2}, \ldots),
$$
and, for $i, j \in \{1, \ldots, n\}$,
$$
y_{i} \circ y_{j} = (y_{i,1} \circ y_{j,1}, y_{i,2} \circ y_{j,2},
\cdots).
$$
Moreover, $H \subseteq \prod_{c \geq 1} Y_{c}$, and it is easy to
verify that $H = \varprojlim (Y_{c}, \widetilde{\rho}_{c,d})$.

Let $F_{n}$ be a free group of rank $n$, with $n \geq 2$, freely
generated by the set $\{f_{1}, \ldots, f_{n}\}$. We write
$F_{n,c}$ for the free nilpotent group of rank $n$ and class $c$;
freely generated by the set $\{x_{1}, \ldots, x_{n}\}$, with
$x_{i} = f_{i} \gamma_{c+1}(F_{n})$, $i = 1, \ldots, n$. Let
$\alpha_{c}$ be the natural epimorphism from $F_{n,c}$ onto
$Y_{c}$ sending $x_{i}$ to $y_{i,c}$, $i = 1, \ldots, n$. Then
$F_{n,c}/{\rm ker} \alpha_{c} \cong Y_{c}$ via an isomorphism
$\widetilde{\alpha}_{c}$ induced by $\alpha_{c}$, and ${\rm ker}
\alpha_{c}$ is a fully invariant subgroup of $F_{n,c}$. Let
$\pi_{c}$ be the natural epimorphism from $F_{n}$ onto $F_{n,c}$
sending $f$ to $f \gamma_{c+1}(F_{n})$ for all $f \in F_{n}$, and
let $\delta_{c} = \alpha_{c} \pi_{c}$. Thus $F_{n}/{\rm ker}
\delta_{c} \cong Y_{c}$ by an isomorphism $\widetilde{\delta}_{c}$
induced by $\delta_{c}$. The group ${\rm ker} \delta_{c}$ is a
fully invariant subgroup of $F_{n}$ since $Y_{n}$ is a relatively
free group. So we have

$$\begin{diagram}
\node{F_n} \arrow{e,t}{\pi_c} \arrow{se,b}{\delta_c} \node{F_{n,c}} \arrow{s,b}{\alpha_c}\\
\node {} \node{Y_c}  \node{Y_d}
\arrow{w,t}{\widetilde{\rho}_{c,d}}
\end{diagram}$$
Denote $h_{i,c} = f_{i} {\rm ker} \delta_{c}$ with $i = 1,
\ldots, n$. Notice that $\widetilde{\delta}_{c}(h_{i,c}) =
\delta_{c}(f_{i}) = y_{i,c}$ for $i = 1, \ldots, n$. Let $c \leq
d$, and let $w(f_1,\ldots, f_n) \in {\rm ker} \delta_{d}$. Then
$\delta_{d}(w(f_1,\ldots, f_n)) = w(y_{1,d}, \ldots, y_{n,d}) =
1_{Y_{d}}$. Since $\widetilde{\rho}_{c,d}$ is a group
homomorphism, we have $w(y_{1,c}, \ldots, y_{n,c}) = 1_{Y_{c}}$
and so, ${\rm ker} \delta_{d} \subseteq {\rm ker} \delta_{c}$. Set
$\psi_{c,d} = (\widetilde{\delta}_{c})^{-1} \widetilde{\rho}_{c,d}
\widetilde{\delta}_{d}$.

$$\begin{diagram}
\node{F_n/\mbox{ker}\delta_d} \arrow{e,t}{\widetilde{\delta}_d} \arrow{s,t}{\psi_{c,d}} \node{Y_d} \arrow{s,b}{\widetilde{\rho}_{c,d}}\\
\node{F_n/\mbox{ker}\delta_c} \arrow{e,t}{\widetilde{\delta}_c}
\node{Y_c}
\end{diagram}$$

It is clearly enough that $\psi_{c,d}(x
{\rm ker} \delta_{d}) = x {\rm ker} \delta_{c}$ for all $x \in
F_{n}$. In particular, $\psi_{c,d}(h_{i,d}) = h_{i,c}$ for $i = 1,
\ldots, n$.

For the rest of the proof, we identify $Y_{c}$ with $F_{n}/{\rm
ker} \delta_{c}$ under the group isomorphism
$\widetilde{\delta}_{c}$. In the light of this identification,
$y_{i,c} = h_{i,c}$ and $y_{i} = (h_{1,c}, h_{2,c}, \ldots)$ with
$i = 1, \ldots, n$. Moreover $H = \varprojlim (F_{n}/{\rm ker}
\delta_{c}, \psi_{c,d})$. Throughout the proof we use both $Y_{c}$
and $F_{n,c}/{\rm ker} \delta_{c}$ without making any distinction.
We claim that $H$ is a relatively free residually torsion-free
nilpotent group of rank $n$. First we shall show that $H$ is
relatively free of rank $n$. Let $\vartheta_{n}$ be the natural
homomorphism from $F_{n}$ into $H$ sending $w = w(f_{1}, \ldots,
f_{n})$ to $(w {\ker} \delta_{1}, w {\rm ker} \delta_{2},
\ldots)$. Since $H$ is generated by the set $\{y_{1}, \ldots,
y_{n}\}$, we have $F_{n}/(\cap {\rm ker} \delta_{c}) \cong H$ by
the isomorphism $\widetilde{\vartheta}_{n}$ induced by
$\vartheta_{n}$. Set $M = \cap {\rm ker} \delta_{c}$. Since $M$ is
a fully invariant subgroup of $F_{n}$, we obtain $H$ is
relatively free of rank $n$.

Let $u = u(y_{1}, \ldots, y_{n}) \in H$. Since $H$ is generated by
the set $\{y_{1}, \ldots, y_{n}\}$, $u$ is written as $u =
 y^{a_{1}}_{i_{1}} \cdots
y_{i_{\kappa}}^{a_{\kappa}}$ with $a_{1}, \ldots,
a_{\kappa} \in {\mathbb{Z}}$ and $i_{1}, \ldots, i_{\kappa}
\in \{1, \ldots, n\}$. Then
$$
u = u(y_{1}, \ldots, y_{n}) = (u(h_{1,1}, \ldots, h_{n,1}),
\ldots, u(h_{1,c}, \ldots, h_{n,c}), \ldots).
$$
Notice that $u(h_{1,c}, \ldots, h_{n,c}) \in Y_{c}$ for all $c$
and, for $c \leq d$, $\psi_{c,d}(u(h_{1,d}, \ldots, h_{n,d})) =$
$u(h_{1,c}, \ldots, h_{n,c})$. Suppose that $y_{1}^{a_{1}} \cdots
y_{n}^{a_{n}} \in H^{\prime}$ for some $a_{1}, \ldots, a_{n} \in
{\mathbb{Z}}$. Then $h_{1,1}^{a_{1}} \cdots h_{n,1}^{a_{n}} \in
Y^{\prime}_{1}$. Since $Y_{1}/Y^{\prime}_{1}$ is free abelian, we
obtain $H/H^{\prime}$ is a free abelian group of rank $n$. Suppose
that $\gamma_{s}(H) = \gamma_{s+1}(H)$ for some $s$. Then
$(y_{i_{1}}, \ldots, y_{i_{s}}) \in \gamma_{s+1}(H)$ for all
$i_{1}, \ldots, i_{s} \in \{1, \ldots, n\}$. Then $(h_{i_{1},s},
\ldots, h_{i_{s},s}) = 1_{Y_{s}}$ for all $i_{1}, \ldots, i_{s}
\in \{1, \ldots, n\}$. Since $Y_{s}$ is freely generated by the
set $\{h_{1,s}, \ldots, h_{n,s}\}$, we obtain $Y_{s}$ has
nilpotency class $s-1$ which is a contradiction. Therefore
$\gamma_{s}(H) \neq \gamma_{s+1}(H)$ for all $s$ and so, there are
no repetitions of terms of the lower central series of $H$.

Our next step is to prove that $H$ is residually torsion-free
nilpotent. For a positive integer $c$, let $\zeta_{c}$ be the
natural epimorphism from $F_{n}/M$ onto $Y_{c}$. (That is,
$\zeta_{c}(wM) = w {\rm ker} \delta_{c}$ for all $w \in F_{n}$.)
For $c \leq d$, it is easily verified that $\psi_{c,d}
\zeta_{d}(f M) = \zeta_{c}(f M)$ for all $f \in F_{n}$. Write
$\widetilde{\zeta}_{c} = \zeta_{c}
(\widetilde{\vartheta}_{n})^{-1}$. Thus
$\widetilde{\zeta}_{c}(y_{i}) = h_{i,c}$ for $i = 1, \ldots, n$.
Let $N_{c} = H/\gamma_{c+1}(H)$. The group $N_{c}$ has
nilpotency class $c$ since $\gamma_s(H)\neq \gamma_{s+1}(H)$ for
all $s$. Let $\beta_{c}$ be the natural epimorphism from $F_{n,c}$
onto $N_{c}$ such that $\beta_{c}(x_{i}) = \overline{y}_{i}$,
where $\overline{y}_{i} = y_{i}\gamma_{c+1}(H)$ $i = 1, \ldots,
n$.  Since $F_{n,c}/F^{\prime}_{n,c} \cong Y_{c}/Y^{\prime}_{c}
\cong N_{c}/N^{\prime}_{c}$, we obtain both ${\rm ker} \alpha_{c}$
and ${\rm ker} \beta_{c}$ are subgroups of $F^{\prime}_{n,c}$.
$$\begin{diagram}
\node[1]{F_{n,c}} \arrow[1]{s,t}{\beta_{c}}    \arrow[1]{e,t}{\alpha_{c}}\node[1]{Y_{c}}\\
\node[1]{N_{c}} \arrow{ne,b,..}{\xi_{c}}
  \end{diagram}$$
We claim that ${\rm ker} \beta_{c} \subseteq {\rm ker}
\alpha_{c}$. Let $v = v(x_{1}, \ldots, x_{n}) \in {\rm ker}
\beta_{c}$. Then $\beta_{c}(v) = v(\overline{y}_{1}, \ldots,
\overline{y}_{n}) = 1_{N_{c}}$ or, equivalently, $v(y_{1}, \ldots,
y_{n}) \in \gamma_{c+1}(H)$. Thus $v(h_{1,q}, \ldots, h_{n,q}) \in
\gamma_{c+1}(Y_{q})$ for all $q$. Hence $v(h_{1,c}, \ldots,
h_{n,c}) = 1_{Y_{c}}$ and so, $v = v(x_{1}, \ldots, x_{n}) \in
{\rm ker} \alpha_{c}$. Therefore ${\rm ker} \beta_{c} \subseteq
{\rm ker} \alpha_{c}$. Let $\gamma_{c}$ be the natural epimorphism
from $F_{n,c}/{\rm ker} \beta_{c}$ onto $F_{n,c}/{\rm ker}
\alpha_{c}$. Moreover, we write $\xi_{c} = \widetilde{\alpha}_{c}
\gamma_{c} \widetilde{\beta}_{c}^{-1}$, where
$\widetilde{\alpha}_{c}$ and $\widetilde{\beta}_{c}$ are the
isomorphisms induced by $\alpha_c$ and $\beta_{c}$ respectively.
$$\begin{diagram}
\node{F_{n,c}/\mbox{ker}\beta_c} \arrow{e,t}{\widetilde{\beta}_c}
\arrow{s,t}{\gamma_c} \node{N_c} \arrow{s,b}{\xi_c}\\
\node{F_{n,c}/\mbox{ker}\alpha_c} \arrow{e,t}{\widetilde{\alpha}_c} \node{Y_c}
\end{diagram}$$
It is easily
verified that $\xi_{c}(v(y_{1}, \ldots, y_{n}) \gamma_{c+1}(H)) =
v(y_{1,c}, \ldots, y_{n,c})$ for $v(y_{1}, \ldots, y_{n}) \in H$.
Let $u = u(y_{1}, \ldots, y_{n}) \in \tau_{c+1}(H)$. Then there
exists $m \in {\mathbb{N}}$ and $v \in \gamma_{c+1}(H)$ such that
$u^{m} = v$. By applying $\xi_{c}$, we have
$$
\begin{array}{cll}
1_{Y_{c}} & = & \xi_{c}(u^{m} \gamma_{c+1}(H)) \\
& = & \xi_{c}(u \gamma_{c+1}(H))^{m} \\
& = & u(y_{1,c}, \ldots, y_{n,c})^{m}.
\end{array}
$$
Since $Y_{c}$ is torsion-free, we get $u(y_{1,c}, \ldots, y_{n,c})
= 1_{Y_{c}}$ and so, $u \gamma_{c+1}(H) \in {\rm ker} \xi_{c}$.
Since $\tau(N_{c}) = \tau_{c+1}(H)/\gamma_{c+1}(H)$, we obtain
$\tau(N_{c}) \subseteq {\rm ker} \xi_{c}$. Let $w = w(y_{1},
\ldots, y_{n})$ $\in \cap_{c \geq 1} \tau_{c+1}(H)$. Then $w \in
\tau_{c+1}(H)$ for all $c$ and so, $w(y_{1,c}, \ldots, y_{n,c}) =
1_{Y_{c}}$ for all $c$. Hence $w$ is a law in $\prod_{c \geq 1}
Y_{c}$. Since $H \leq \prod_{c \geq 1} Y_{c}$, we have $w$ is a
law in $H$ and so, $w = 1_{H}$. Therefore $H$ is a residually
torsion-free nilpotent group.

Finally we show the last part of Theorem B (II). By Theorem A
(II), we have $L(c) \cong {\cal L}_{\mathbb{Q}}(Y_{c})$ for all
$c$ via an isomorphism $\lambda_{c}$, say, sending $y_{i,c}$ to
$\log h_{i,c}$ for $i = 1, \ldots, n$. As in the proof of Lemma
\ref{pap}, there exists a Lie algebra epimorphism
$\xi_{\psi_{c,d}}$, with $c \leq d$, from ${\cal
L}_{\mathbb{Q}}(Y_{d})$ onto ${\cal L}_{\mathbb{Q}}(Y_{c})$ such
that $\xi_{\psi_{c,d}}(\log (f {\rm ker} \delta_{d})) = \log (f
{\rm ker} \delta_{c})$ for all $f \in F_{n}$. Form the inverse
limit of the family $({\cal L}_{\mathbb{Q}}(Y_{c}),
\xi_{\psi_{c,d}})$, $\varprojlim {\cal L}_{\mathbb{Q}}(Y_{c})$,
and define a mapping
$$\lambda : \widehat{L} \rightarrow \varprojlim {\cal L}_{\mathbb{Q}}(Y_{c}),$$ as
follows:
$$
\lambda(u_{1}+L^{\prime}, (u_{1}+u_{2})+\gamma_{3}(L), \ldots) =
(\lambda_{1}(u_{1}+L^{\prime}),
\lambda_{2}((u_{1}+u_{2})+\gamma_{3}(L)), \ldots),
$$
where $u_{i} \in \gamma_{i}(L)$ for all $i$. It is easily verified
that $\lambda$ is a Lie algebra monomorphism. Let $v = (v_{1},
v_{2}, \ldots) \in \varprojlim {\cal L}_{\mathbb{Q}}(Y_{c})$. Thus
$\xi_{\psi_{c,d}}(\lambda_{d}(u_{d})) = v_{c}$ for some $u_{d} \in
L(d)$. We claim that $(u_{1}, u_{2}, \ldots) \in \widehat{L}$.
Since $\xi_{\psi_{c,d}} \lambda_{d} = \lambda_{c} \rho_{c,d}$
$$\begin{diagram}
\node[1]{L(d)} \arrow{s,t}{\lambda_d} \arrow[1]{e,t}{\rho_{c,d}}  \node[1]{L(c)} \arrow{s,b}{\lambda_c}\\
\node[1]{{\cal L}_{\mathbb{Q}}(Y_{d})}
\arrow{e,t}{\xi_{\psi_{c,d}}} \node[1]{{\cal
L}_{\mathbb{Q}}(Y_{c})}
  \end{diagram}$$
for all $c \leq d$ and $\lambda_{c}$ is $1-1$, we obtain
$\rho_{c,d}(u_{d}) = u_{c}$. Hence the Lie algebra monomorphism
$\lambda$ is onto and so, $\lambda$ is a Lie algebra isomorphism.

For $i$, $i = 1, \ldots, n$, let
$$
t^{\prime}_{i} = (\log h_{i,1}, \log h_{i,2}, \ldots, \log
h_{i,c}, \ldots) \in \varprojlim {\cal L}_{\mathbb{Q}}(Y_{c}),
$$
and let $\Lambda_{\mathbb{Q}}(H)$ be the Lie subalgebra of
$\varprojlim ({\cal L}_{\mathbb{Q}}(Y_{c}), \xi_{\psi_{c,d}})$
generated by the set \linebreak $\{t^{\prime}_{1}, \ldots,
t^{\prime}_{n}\}$. Since $L$ is residually nilpotent, $L$ is
embedded into $\widehat{L}$ via a Lie algebra monomorphism
$\lambda^{\prime}$, say, sending $y_{i}$ to $(y_{i} + L^{\prime},
y_{i} + \gamma_{3}(L), \ldots)$ for $i = 1, \ldots, n$. Since
$\lambda$ is a Lie algebra isomorphism, we obtain $L$ is
isomorphic to $\Lambda_{\mathbb{Q}}(H)$ via $\lambda
\lambda^{\prime}$. Let $L_{n}$ be the free Lie algebra freely
generated by the set $\{\ell_{1}, \ldots, \ell_{n}\}$. By the
proof of Theorem \ref{3.18} (for $G_{n} = H$), we have ${\cal
L}(H) \cong L_{n}/{\rm ker} \sigma_{n}$. Furthermore, by applying
the proof of Theorem \ref{3.18} for ${\cal L}_{\mathbb{Q}}(Y_{c})$
(for all $c$), we obtain $\Lambda_{\mathbb{Q}}(H) \cong L_{n}/{\rm
ker} \sigma^{\prime}_{n}$. To prove that $\Lambda_{\mathbb{Q}}(H)$
is a homomorphic image of ${\cal L}(H)$, it is enough to show that
${\rm ker} \sigma_{n} \subseteq {\rm ker} \sigma^{\prime}_{n}$.
Since $H/\tau_{c+1}(H)$ is mapped onto $Y_{c}$ via an epimorphism
induced by $\xi_{c}$, we have ${\cal
L_{\mathbb{Q}}}(H/\tau_{c+1}(H))$ is mapped onto ${\cal
L}_{\mathbb{Q}}(Y_{c})$ via a Lie algebra epimorphism sending
$\log (y_{i} \tau_{c+1}(H))$ to $\log h_{i,c}$ for all $i$. Let $v
= v(\ell_{1}, \ldots, \ell_{n}) \in {\rm ker} \sigma_{n}$. Then
$v(\log (y_{1} \tau_{c+1}(H)), \ldots, \log (y_{n} \tau_{c+1}(H)))
= 0$ for all $c$ (see the proof of Theorem \ref{3.18}). Hence
$v(\log h_{1,c}, \ldots, \log h_{n,c}) = 0$ for all $c$ and so, $v
= v(\ell_{1}, \ldots, \ell_{n}) \in {\rm ker} \sigma^{\prime}_{n}$
(by similar arguments as in the proof of Theorem \ref{3.18}).
Therefore $\Lambda_{\mathbb{Q}}(H)$ is a homomorphic image of
${\cal L}(H)$.  \qed

\section{An Example}

We shall give an example of a finitely generated Magnus nilpotent
group $G$, not relatively free, such that ${\cal
L}_{\mathbb{Q}}(G)$ is not isomorphic to $L_{\mathbb{Q}}(G)$ as
Lie algebra. We modify (and analyze) an example of a finitely
generated nilpotent Lie algebra given in \cite[page 210]{andreev}.
Let ${\cal L}_{4}$ be a free Lie algebra of rank $4$; freely
generated by the set $\{\ell_{1}, \ldots, \ell_{4}\}$. Let
$L_{4,3} = {\cal L}_{4}/\gamma_{4}({\cal L}_{4})$ and let $x_{i} =
\ell_{i} + \gamma_{4}({\cal L}_{4})$, $i = 1, \ldots, 4$. Thus
$L_{4,3}$ is a free nilpotent Lie algebra of rank $4$ and class
$3$ with a free generating set $\{x_{1}, \ldots, x_{4}\}$. Set $v
= [x_{1}, x_{2}] + [x_{3}, x_{4}, x_{3}]$, and write $L =
L_{4,3}/I$, where $I$ is the ideal of $L_{4,3}$ generated by $v$.
It is easily verified that every element of $I$ has the form
$$
a ([x_{1}, x_{2}] + [x_{3}, x_{4}, x_{3}]) + \sum_{i=1}^{4} a_{i}
[x_{1}, x_{2}, x_{i}],
$$
where $a, a_{i} \in {\mathbb{Q}}$, $i = 1, \ldots, 4$. For $i = 1,
\ldots, 4$, let $y_{i} = x_{i} + I$. Then $[y_{1}, y_{2}] = -
[y_{3}, y_{4}, y_{3}] \in \gamma_{3}(L)$ and $[y_{1}, y_{2},
y_{i}] = 0$ with $i = 1, \ldots, 4$. Since $I$ is a proper subset
of $L^{\prime}_{4,3}$, we obtain $L$ is not abelian, and the set
$\{y_{i}+L^{\prime}; i = 1, \ldots, 4\}$ is a ${\mathbb{Q}}$-basis
of $L/L^{\prime}$. Since $L$ is finitely generated nilpotent Lie
algebra and $\dim (L/L^{\prime}) = 4$, we obtain $L$ is generated
by the set $\{y_{1}, \ldots, y_{4}\}$. Suppose that $L$ is
relatively free. Since $\dim (L/L^{\prime}) = 4$ and $L$ is
nilpotent, it is easily verified that the set $\{y_{1}, \ldots,
y_{4}\}$ freely generates $L$. Then $[y_{1}, y_{2}] = [y_{3},
y_{4}, y_{3}] = 0$ and so, $[x_{1}, x_{2}] \in I$ which is not
valid. Thus $L$ is not relatively free.

Give on $L$ the structure of a group, denoted $R$, by means of the
BCH formula. That is, for $x, y \in L$,
$$
x\circ y = x + y + \frac{1}{2}[x, y] + \frac{1}{12}[x,y,y] -
\frac{1}{12} [x, y, x]. \eqno (17)
$$
Recall that $0$ is the unit of $R$ and $-x$ is the inverse of $x$
with respect to the group operation $\circ$. Let $H$ be the
subgroup of $R$ generated by the set $\{y_{i}; i = 1, \ldots,
4\}$. For $u, v \in H$,
$$
(u, v) = [u, v] + \frac{1}{2}[u,v,u] + \frac{1}{2} [u,v,v]. \eqno
(18)
$$
Thus
$$
(y_{1}, y_{2}) = (y_{3}, y_{4}, y_{3})^{-1}.
$$
By Lemma \ref{2.1}, $\tau_{2}(H) = H^{\prime}$, and so, $H/H^{\prime}$
is free abelian of rank $4$. By (17), for $\alpha \in
{\mathbb{Z}}$ and $u, v \in H$, $(u, v)^{\alpha} = \alpha (u, v)$.
By (18), we obtain
$$
(u, v)^{\alpha} = \alpha [u, v] + \frac{\alpha}{2} [u, v, v] +
\frac{\alpha}{2} [u, v, u]. \eqno (19)
$$
Suppose that
$$
\prod_{i=1}^{2} (y_{3}, y_{i})^{\alpha_{i}} \prod_{j=1}^{3}
(y_{4}, y_{j})^{\beta_{j}} \in \gamma_{3}(H),
$$
where $\alpha_{i}, \beta_{j} \in {\mathbb{Z}}$, $i = 1, 2$ and $j
= 1, 2, 3$. Since $\gamma_{3}(H)$ is generated by the elements of
the form $(y_{i}, y_{j}, y_{k})$ with $i, j, k \in \{1, \ldots,
4\}$, we obtain from (19) and the form of elements of $I$ that
$$
a [x_{1}, x_{2}] + \sum_{i=1}^{2} \alpha_{i} [x_{3}, x_{i}] +
\sum_{j=1}^{3} \beta_{j} [x_{4}, x_{j}] \in \gamma_{3}(L_{4,4}),
$$
for $a \in {\mathbb{Q}}$. Since $L_{4,3}$ is a free nilpotent Lie
algebra, we have $\alpha_{i} = 0$ ($i = 1, 2$) and $\beta_{j} = 0$
($j = 1, 2, 3$). Therefore $H^{\prime}/\gamma_{3}(H)$ is
torsion-free. Finally, it is easily shown that $\gamma_{3}(H)$ is
torsion-free. Thus $H$ is a Magnus group. Suppose to a contrary
that $H$ is relatively free. Then, by Theorem A (I), ${\cal
L}_{\mathbb{Q}}(H)$ is relatively free. Since both $R$ and $\exp
{\cal L}_{\mathbb{Q}}(H)$ are Mal'cev completions of $H$, we
obtain $R \cong \exp {\cal L}(H)$ as groups and so, $L \cong {\cal
L}_{\mathbb{Q}}(H)$ as Lie algebras. That is, $L$ is relatively
free which a contradiction, and so, $H$ is not relatively free.
Suppose that there exists a Lie algebra isomorphism $\zeta$ from
$L_{{\mathbb{Q}}}(H)$ into $L$. By Lemma \ref{3.2},
$L_{{\mathbb{Q}}}(H)$ is generated by the set $\{y_{1}H^{\prime},
\ldots, y_{4}H^{\prime}\}$. Write $\zeta(y_{i}H^{\prime}) =
z_{i}$, $i = 1, \ldots, 4$. Since
$$
[y_{1}H^{\prime}, y_{2}H^{\prime}] = (y_{1}, y_{2})\gamma_{3}(H) =
0~~{\rm in}~~L_{{\mathbb{Q}}}(H),
$$
we have $[z_{1}, z_{2}] = 0$ in $L$. For $i = 1, 2$,
$$
z_{i} = \sum_{j=1}^{4} \alpha_{ji} y_{j} + \sum_{\substack{\kappa \neq 2\\
1 \leq \ell < \kappa \leq 4}} \beta_{\kappa \ell i}
[y_{\kappa}, y_{\ell}] + v_{i},
$$
where $v_{i} \in \gamma_{3}(L)$. Since $\{z_{1}+L^{\prime},
z_{2}+L^{\prime}\}$ is linearly independent, we obtain
$\alpha_{\mu 1}\alpha_{\nu 2} - \alpha_{\nu 1} \alpha_{\mu 2} \neq
0$ for some $\mu, \nu \in \{1, \ldots, 4\}$ with $\mu \neq \nu$.
Since $[y_{1}, y_{2}] + [y_{3}, y_{4}, y_{3}] = 0$, and working
with a basis of $L$ consisting of "basic commutators", we obtain
the following equations
$$
\begin{array}{ccl}
 \alpha_{11}\alpha_{22} - \alpha_{12} \alpha_{21} & = & -\alpha_{32} \beta_{431} + \alpha_{31} \beta_{432} \\
 \alpha_{11}\alpha_{32} - \alpha_{12} \alpha_{31} & = & 0 \\
 \alpha_{11}\alpha_{42} - \alpha_{12} \alpha_{41} & = & 0 \\
 \alpha_{21}\alpha_{32} - \alpha_{22} \alpha_{31} & = & 0 \\
 \alpha_{21}\alpha_{42} - \alpha_{22} \alpha_{41} & = & 0 \\
\alpha_{31}\alpha_{42} - \alpha_{32} \alpha_{41} & = & 0 \\
\alpha_{41}\beta_{432} - \alpha_{42} \beta_{431} & = & 0
\end{array}
$$
Thus $\alpha_{11} \alpha_{22} - \alpha_{12} \alpha_{21} \neq 0$.
Therefore $\alpha_{11} \alpha_{22} \neq 0$ or $\alpha_{12}
\alpha_{21} \neq 0$. Suppose that $\alpha_{11} \alpha_{22} \neq
0$. If $\alpha_{12} \alpha_{31} \neq 0$ then $\alpha_{11}
\alpha_{32} \neq 0$ and so $\frac{\alpha_{11}}{\alpha_{12}} =
\frac{\alpha_{31}}{\alpha_{32}}$. Furthermore
$\frac{\alpha_{21}}{\alpha_{22}} =
\frac{\alpha_{31}}{\alpha_{32}}$. Hence $\alpha_{11} \alpha_{22} -
\alpha_{12} \alpha_{21} = 0$ which is a contradiction. If
$\alpha_{12} \alpha_{31} = 0$ then $\alpha_{12} = 0$ or
$\alpha_{31} = 0$. Since $\alpha_{11} \neq 0$ we get $\alpha_{32}
= 0$. Then $\alpha_{22} \alpha_{31} = 0$ and so $\alpha_{13} =
\alpha_{32} = 0$ which is a contradiction. Similar arguments may
be applied if $\alpha_{12} \alpha_{21} \neq 0$. Therefore we
obtain $L$ is not isomorphic to $L_{{\mathbb{Q}}}(H)$ as Lie
algebras.

\vskip .120 in

\noindent\emph{Faculty of Sciences, Department of Mathematics, Aristotle
University of Thessaloniki, GR 541 24, Thessaloniki, Greece},
e-mail: kkofinas@math.auth.gr

\bigskip

\noindent\emph{Department of Mathematics, University of the Aegean, GR 832
00, Karlovassi, Samos, Greece}, e-mail: vmet@aegean.gr

\bigskip

\noindent\emph{Faculty of Sciences, Department of Mathematics, Aristotle
University of Thessaloniki, GR 541 24, Thessaloniki, Greece},
e-mail: apapist@math.auth.gr


\addcontentsline{toc}{section}{References}
\begin{thebibliography}{99}

\bibitem{as} R.K. Amayo and I. Stewart, Infinite-dimensional Lie algebras, (Noordhoff International Publishing, 1974).

\bibitem{andreev} K.K. Andreev, Nilpotent groups and Lie algebras,
{\it Algebra i Logika} {\bf 7} (1968), 4--14.

\bibitem{bahturin}  Yu.A. Bahturin, Identical relations in Lie
algebras, (Nauka, Moscow, 1985), (Russian). Translation: (VNU
Science Press, Utrecht, 1987).

\bibitem{baumslag} G. Baumslag, Lectures notes on nilpotent groups, C.B.M.S. Regional
Conference Series 2, American Mathematical Society, Providence
1971.

\bibitem{bh} M.R. Bridson and A. Haefliger, Metric spaces of non-positive curvature. Springer-Verlag Berlin, 1999.
Die Grundlehren der mathematischen Wissenschaften, Band 319.

\bibitem{bg} R.M. Bryant, J.R.J. Groves, Algebraic groups of
automorphisms of nilpotent groups and Lie algebras, {\it J. London
Math. Soc.} {\bf 33} (1986), 453--466.

\bibitem{fm} B. Farb and L. Mosher, Problems on the geometry of finitely generated
solvable groups. In: Crystallographic groups and their generalizations (Kortrijk, 1999),
121--134, {\it Contemp. Math.} {\bf 262}, Amer. Math. Soc., Providence, RI, 2000.


\bibitem{hirsch} K.A. Hirsch, Infinite soluble groups II, {\it Proc. London Math. Soc.}
{\bf 44} (1938), 336--365.

\bibitem{hb} B. Huppert, N. Blackburn, Finite groups II,
Springer, 1982.

\bibitem{jacobson}  N. Jacobson, Lie algebras, (Wiley, New York, 1962).

\bibitem{jennings}  S.A. Jennings, The group ring of a class of infinite nilpotent groups, {\it Canad. J. Math.} {\bf 7} (1955), 169--187.

\bibitem{kp}  C.E. Kofinas and A.I. Papistas, Automorphisms of relatively free nilpotent Lie algebras, (submitted).

\bibitem{kurosh} A.G. Kurosh, The theory of groups, Vol. II, Chelsea, 1956.

\bibitem{luck} W. L\"uck, Survey on Geometric Group Theory, arXiv:0806.3771.

\bibitem{mks} W. Magnus, A. Karrass and D. Solitar, Combinatorial
Group Theory, (Wiley, New York, 1966).

\bibitem{neumann}  H. Neumann, Varieties of groups, Ergebnisse der mathematik and ihner Grenzgebiete, band 37 (Springer-Verlag, 1967).

\bibitem{pansu} P. Pansu, M\'etrique de Carnot-Carath\'eodory et quasiisometries des espaces sym\'etrique de rang un, {\it Ann. Math.} {\bf 129} (1989), 1--60.

\bibitem{segal}  D. Segal, Polycyclic groups, (Cambridge University Press, 1983).

\bibitem{shmelkin} A.L. Shmel'kin, Free poly-nilpotent groups,
{\it Izv. Akad. Nauk SSSR, Ser. Mat.} {\bf 28} (1964), 91-122.
\end{thebibliography}
\end{document}